\documentclass[english]{amsart}
\usepackage{amsmath}

\usepackage{amssymb}
\usepackage{svn}

\usepackage{amscd}
\usepackage[all,cmtip,line]{xy}

\SVN $LastChangedDate: 2012-10-30 15:27:28 +0000 (Tue, 30 Oct 2012) $
\SVN $LastChangedRevision: 28 $

\def\AA{\mathbb{A}}

\def\CC{\mathbb{C}}

\def\FF{\mathbb{F}}

\def\QQ{\mathbb{Q}}
\def\RR{\mathbb{R}}

\def\ZZ{\mathbb{Z}}

\newcommand{\uA}{{\underline{A}}}

\newcommand{\calC}{{\mathcal{C}}}

\newcommand{\onto}{\twoheadrightarrow}

\newlength{\ownl}

\newcommand{\End}{{\operatorname{End}\,}}

\newcommand{\Fil}{{\operatorname{Fil}\,}}

\newcommand{\Frob}{{\operatorname{Frob}}}

\newcommand{\gr}{{\operatorname{gr}\,}}

\newcommand{\Hom}{{\operatorname{Hom}\,}}

\newcommand{\Res}{{\operatorname{Res}}}

\newcommand{\Spec}{{\operatorname{Spec}\,}}
\newcommand{\Spf}{{\operatorname{Spf}\,}}

\newcommand{\GL}{\operatorname{GL}}

\newcommand{\SL}{\operatorname{SL}}

\newcommand{\cris}{{\operatorname{cris}}}

\newcommand{\tor}{{\operatorname{tor}}}

\newcommand{\A}{{\mathbb{A}}}

\newcommand{\C}{{\mathbb{C}}}
\newcommand{\D}{{\mathbb{D}}}

\newcommand{\F}{{\mathbb{F}}}
\newcommand{\G}{{\mathbb{G}}}

\newcommand{\Q}{{\mathbb{Q}}}

\newcommand{\Z}{{\mathbb{Z}}}

\newcommand{\CA}{{\mathcal{A}}}

\newcommand{\CE}{{\mathcal{E}}}
\newcommand{\CF}{{\mathcal{F}}}
\newcommand{\CG}{{\mathcal{G}}}
\newcommand{\CH}{{\mathcal{H}}}
\newcommand{\CI}{{\mathcal{I}}}
\newcommand{\CJ}{{\mathcal{J}}}

\newcommand{\CL}{{\mathcal{L}}}
\newcommand{\CM}{{\mathcal{M}}}
\newcommand{\CN}{{\mathcal{N}}}
\newcommand{\CO}{{\mathcal{O}}}
\newcommand{\CP}{{\mathcal{P}}}

\newcommand{\CS}{{\mathcal{S}}}

\newcommand{\gI}{{\mathfrak{I}}}

\newcommand{\gc}{{\mathfrak{c}}}
\newcommand{\gd}{{\mathfrak{d}}}

\newcommand{\gm}{{\mathfrak{m}}}
\newcommand{\gn}{{\mathfrak{n}}}

\newcommand{\gp}{{\mathfrak{p}}}
\newcommand{\gq}{{\mathfrak{q}}}

\newcommand{\Qbar}{{\overline{\Q}}}

\DeclareMathOperator{\lcm}{lcm}

\newcommand{\Qpbar}{\overline{\Q}_p}

\newcommand{\Fpbar}{\overline{\F}_p}

\newcommand{\Ver}{\operatorname{Ver}}

\def\smallmat#1#2#3#4{\bigl(\begin{smallmatrix}{#1}&{#2}\\{#3}&{#4}\end{smallmatrix}\bigr)}

\usepackage{hyperref}

\newcommand{\Ha}{{\mathrm{Ha}}}
\newcommand{\Ig}{{\mathrm{Ig}}}

\newcommand{\be}{{\mathbf{e}}}
\newcommand{\f}{\mathbf{f}}
\newcommand{\bh}{\mathbf{h}}
\newcommand{\bk}{\mathbf{k}}
\newcommand{\bl}{\mathbf{l}}
\newcommand{\bm}{\mathbf{m}}
\newcommand{\bn}{\mathbf{n}}
\newcommand{\bt}{\mathbf{t}}

\newcommand{\Nm}{\mathrm{Nm}}
\newcommand{\uhp}{\mathfrak{H}}
\newcommand{\Shom}{{\mathcal{H}om}}

\newcommand{\Lie}{{\mathcal{L}ie}}
\newcommand{\crys}{\mathrm{crys}}
\newcommand{\tcalC}{{\widetilde{\calC}}}

\newcommand{\ord}{{\operatorname{ord}}}
\newcommand{\tord}{{\operatorname{tord}}}

\newcommand{\dr}{{\operatorname{dR}}}
\newcommand{\tot}{{\operatorname{tot}}}
\newcommand{\im}{{\operatorname{im}}}

\newcommand{\ol}{\overline}
\newcommand{\os}{{\ol{s}}}
\newcommand{\wS}{{\widehat{S}}}

\theoremstyle{plain}
\newtheorem{theorem}{Theorem}[subsection]
\newtheorem{proposition}[theorem]{Proposition}
\newtheorem{corollary}[theorem]{Corollary}
\newtheorem{lemma}[theorem]{Lemma}

\newtheorem{ithm}{Theorem}

\theoremstyle{definition}
\newtheorem{remark}[theorem]{Remark}
\newtheorem{definition}[theorem]{Def\/inition}

\newcommand{\begpf}{\noindent{\bf Proof.}\enspace}
\newcommand{\epf}{{\ifhmode\unskip\nobreak\hfil\penalty50 \hskip1em
\else\nobreak\fi \nobreak\mbox{}\hfil\mbox{$\square$} \parfillskip=0pt
\finalhyphendemerits=0 \par\vskip5pt}}

\begin{document}

\title[Weight-shifting operators on Hilbert modular forms]{Geometric weight-shifting operators on Hilbert modular forms in characteristic $p$}
\author{Fred Diamond}
\email{fred.diamond@kcl.ac.uk}
\address{Department of Mathematics,
King's College London, WC2R 2LS, UK}

\maketitle

\date{November 2020}

\begin{abstract}
We carry out a thorough study of weight-shifting operators on Hilbert modular forms in characteristic $p$,
generalizing the author's prior work with Sasaki to the case where $p$ is ramified in the totally real field.
In particular we use the partial Hasse invariants and Kodaira--Spencer filtrations defined by Reduzzi and Xiao
to improve on Andreatta and Goren's construction of partial $\Theta$-operators, obtaining ones whose effect
on weights is optimal from the point of view of geometric Serre weight conjectures.  Furthermore we describe
the kernels of partial $\Theta$-operators in terms of images of geometrically constructed partial Frobenius operators.
Finally we apply our results to prove a partial positivity result for minimal weights of mod $p$ Hilbert modular forms.
\end{abstract}


\section{Introduction}
The study of weight-shifting operations on modular forms has a rich and fruitful history.
Besides those naively obtained from the graded algebra structure on the space
of classical modular forms of all weights, there is a deeper construction due to
Ramanujan~\cite{ram} which shifts the weight by two using differentiation, leading
to a more general theory of Maass--Shimura operators.
Analogous weight-shifting operations in characteristic $p$, first studied by Swinnerton-Dyer
and Serre~\cite{serre}, take on special significance in the context of congruences between
modular forms and the implications for associated Galois representations.  In particular one
has the following linear maps from the space of mod $p$ modular forms of weight $k$
and some fixed level $N$ prime to $p$:
\begin{itemize}
\item multiplication by a Hasse invariant $H$, to forms of weight $k+p-1$;
\item a differential operator $\Theta$, to forms of weight $k+p+1$;
\item a linearized $p$-power map $V$, to forms of weight $pk$.
\end{itemize}
These maps all have simple descriptions in terms of associated $q$-expansions: 
if $f$ has $q$-expansion $\sum a_n q^n$, then that of $Hf$ (resp.~$\Theta f$, $V f$)
is  $\sum a_n q^n$ (resp.~$\sum n a_n q^n$, $\sum a_n q^{pn}$).  Following the work
of Swinnerton-Dyer and Serre, there were further significant developments to the theory
due to Katz~\cite{katz1, katz2} (interpreting the constructions more geometrically), 
Jochnowitz~\cite{joch1, joch2} (on the weight filtration and Tate's $\Theta$-cycles)
and Gross~\cite{gross} (in the study of companion forms), providing crucial ingredients
for Edixhoven's proof of the weight part of Serre's Conjecture in~\cite{edix}.

Suppose now that $F$ is a totally real field of degree $d = [F:\Q]$ and consider spaces of
Hilbert modular forms of weight $\bk \in \ZZ^\Sigma$ and fixed level prime to $p$, where $\Sigma$
denotes the set of embeddings $\{\sigma: F \hookrightarrow \Qpbar\}$.  For such spaces of $p$-adic modular
forms, Katz \cite{Katz3} constructed a family of commuting differential operators $\Theta_\sigma$, indexed by the
$d$-embeddings $\sigma \in \Sigma$.  The theory was further developed by Andreatta
and Goren~\cite{AG} who, building on Katz's work and Goren's definition of partial Hasse invariants in \cite{Go1, Go2}
(if $p$ is unramified in $F$), defined partial $\Theta$-operators on spaces of mod $p$ Hilbert modular forms.

Under the assumption that $p$ is unramified in $F$, some aspects of the construction
of partial $\Theta$-operators in \cite{AG} were simplified in \cite{DS}, which also went on to define
partial Frobenius operators (generalizing $V$) geometrically and use their image to describe kernels of
partial $\Theta$-operators.  When $p$ is ramified in $F$, the effectiveness of the
approach in \cite{AG} was limited by the singularities of the available (Deligne--Pappas)
model for the Hilbert modular variety.  Since then however, a smooth integral
model was constructed by Pappas and Rapoport~\cite{PR}, and the theory of partial
Hasse invariants was further developed in this context by Reduzzi and Xiao in \cite{RX}.
The theory of partial $\Theta$-operators was revisited in that light by Deo, Dimitrov and Wiese in \cite{DDW}, where
they closely follow \cite{AG}.   Here we instead exploit the observations and techniques introduced in \cite{DS},
applying them directly to the special fibre of the Pappas--Rapoport model to construct and relate partial $\Theta$
and Frobenius operators.  In particular this eliminates extraneous multiples of partial Hasse invariants that
appear in \cite{DDW}, and yields results whose implications for minimal weights are 
motivated by the forthcoming generalization to the ramified case of the geometric Serre weight
conjectures of \cite{DS}.  The main contributions of this paper may be summarized as follows:
\begin{itemize}
\item a construction of operators $\Theta_\tau$ with optimal effect on weight (Theorem~\ref{thm:theta});
\item a geometric construction of partial Frobenius operators $V_\gp$ (see \S\ref{sec:V});
\item a description of the kernel of $\Theta_\tau$ in terms of the image of $V_\gp$ (see \S\ref{sec:ker});
\item an application to positivity of minimal weights (Theorem~\ref{thm:ptwt0}).
\end{itemize}

We should emphasize that the focus of this paper is entirely on Hilbert modular forms in characteristic $p$.
There is also a rich theory of $\Theta$-operators on $p$-adic automorphic forms which has seen major
progress recently in the work of de Shalit and Goren~\cite{dSG2}, and Eischen, Fintzen, Mantovan and
Varma~\cite{EFMV}, which in turn has implications in the characteristic $p$ setting~\cite{dSG1, dSG2, EFGMM, EM}.
Another advance in characteristic $p$ has been Yamauchi's construction~\cite{yam} of $\Theta$-operators for 
mod $p$ Siegel modular forms of degree two.  We remark however that all of the work just mentioned
only considers automorphic forms on reductive groups which are unramified at $p$; the novelty of this
paper is largely in the treatment of ramification at $p$. 

\bigskip

We now describe the contents in more detail.

We first set out some basic notation and constructions in \S\ref{sec:notation}.  In particular we
fix a prime $p$,  a totally real field $F$ of degree\footnote{Including the case $F = \QQ$ would
introduce different complications in the treatment of cusps and provide no new results.}
$d = [F:\QQ] > 1$, and let $\CO_F$ denote the ring of integers of $F$ and $S_p$ the set of prime ideals of $\CO_F$ dividing $p$.
For each $\gp \in S_p$, let $\Sigma_{\gp,0}$ denote the set of $f_\gp$ embeddings $\CO_F/\gp \to \Fpbar$
and $\Sigma_\gp$ the set of $e_\gp f_\gp$ embeddings $F_\gp \to \Qpbar$, where $f_\gp$ (resp.~$e_\gp$)
is the residual (resp.~ramification) degree of $\gp$.  We let
$$\begin{array}{rcccl} \Sigma_0 & = & \coprod_{\gp \in S_p} \Sigma_{\gp,0} & = & \{\,\tau_{\gp,i}\,|\, \gp \in S_p, i \in \ZZ/f_\gp\ZZ\,\}\\
\mbox{and} \qquad \Sigma & = & \coprod_{\gp \in S_p} \Sigma_\gp & =  & \{\,\theta_{\gp,i,j}\,|\, \gp \in S_p, i \in \ZZ/f_\gp\ZZ,  j = 1,\ldots,e_\gp\,\},\end{array}$$
where each $\tau_{\gp,0} \in \Sigma_{\gp,0}$ is chosen arbitrarily, $\tau_{\gp,i} = \tau_{\gp,0}^{p^i}$ and 
$\theta_{\gp,i,1},\ldots,\theta_{\gp,i,e_\gp}$ is any ordering of the lifts of $\tau_{\gp,i}$ to $\Sigma_\gp$.
We also define a ``right-shift'' permutation $\sigma$ of $\Sigma$ by
$$\sigma(\theta_{\gp,i,j}) = \left\{ \begin{array}{ll}  \theta_{\gp,i,j-1},&\mbox{if $j > 1$};\\   \theta_{\gp,i-1,e_\gp},&\mbox{if $j=1$.}\end{array}\right.$$

In \S\ref{sec:models} we recall the definition of the Pappas--Rapoport model $Y_U$ for the Hilbert modular variety of level $U$,
where $U$ is any sufficiently small open compact subgroup of $\GL_2(\AA_{F,\f})$ of level prime to $p$.  This may be
viewed as a coarse moduli space for Hilbert-Blumenthal abelian varieties with additional structure, where this additional
structure includes a suitable collection of filtrations on direct summands of its sheaf of invariant differentials.   The scheme
$Y_U$ is then smooth of relative dimension $d$ over $\CO$, where $\CO$ is the ring of integers of a finite extension
of $\Q_p$ in $\Qpbar$.  Since the main results of the paper concern Hilbert modular forms in characteristic $p$, we will
restrict our attention to this setting for the remainder of the Introduction, and let $\ol{Y}_U = Y_{U,\FF}$ where $\FF$
is the residue field of $\CO$.

In \S\S\ref{sec:pairings}--\ref{sec:bundles} we construct the automorphic line bundles $\ol{\CA}_{\bk,\bl}$ on $\ol{Y}_{U}$
for all $\bk$, $\bl \in \ZZ^\Sigma$ and sufficiently small $U$ (of level prime to $p$), and  
define the space of Hilbert modular forms of weight $(\bk,\bl)$ and level $U$ over $\FF$ to be
$$M_{\bk,\bl}(U;\FF) = H^0(\ol{Y}_{U},\ol{\CA}_{\bk,\bl}).$$
The spaces are equipped with a natural Hecke action making 
$$M_{\bk,\bl}(\FF) := \varinjlim_U M_{\bk,\bl}(U;\FF)$$
a smooth admissible representation of $\GL_2(\AA_{F,\f})$ over $\FF$.
A key point, as already observed in \cite{DS} in the unramified case, is that the parity condition on $\bk$ imposed
in the definition of Hilbert modular forms in characteristic zero (for the group $\Res_{F/\QQ} \GL_2$) disappears in characteristic $p$.
We remark also that the effect of the weight parameter $\bl$ (in characteristic $p$) is to introduce twists by torsion bundles
that make various constructions, in particular that of partial $\Theta$-operators, compatible with the natural Hecke action.

 In \S\S\ref{sec:KS}--\ref{sec:stratification} we recall results of Reduzzi and Xiao~\cite{RX} that will underpin 
 our construction of partial $\Theta$-operators.  Firstly there is a natural Kodaira--Spencer filtration on direct
 summands of $\Omega^1_{\ol{Y}_U/\FF}$ whose graded pieces are isomorphic to
 the automorphic line bundles $\ol{\CA}_{2\be_\theta,-\be_\theta}$ (where $\be_\theta$ denotes the basis element
 of $\ZZ^\Sigma$ indexed by $\theta$).  Secondly for each $\theta = \theta_{\gp,i,j} \in \Sigma$, there is a partial Hasse invariant 
 $$H_\theta \in M_{\bh_\theta,\bf{0}}(U;\FF),\qquad\mbox{where \,\,$ \bh_\theta = n_\theta \be_{\sigma^{-1}\theta} - \be_\theta$}$$
 with $n_\theta = p$ if $\theta_{\gp,i,1}$ for some $\gp,i$, and $n_\theta = 1$ otherwise.  Multiplication by the partial
 Hasse invariant $H_\theta$ thus defines a map
 $$\cdot H_\theta:  M_{\bk,\bl}(U;\FF)  \longrightarrow   M_{\bk+\bh_\theta,\bl}(U;\FF)$$
which is easily seen to be Hecke-equivariant.
 We also analogously define invariants $G_\theta \in M_{\bf{0},\bh_\theta}(U;\FF)$ which trivialize the bundles 
$\ol{\CA}_{\bf{0},\bh_\theta}$.  The partial Hasse invariants $H_\theta$ refine the ones
defined by Andreatta and Goren in~\cite{AG} and give rise to a natural stratification on $\ol{Y}_U$
and a notion of minimal weight $\bk_{\min}(f)$ for non-zero $f \in M_{\bk,\bl}(U;\FF)$, which the main
result of \cite{DK2} shows lies in a certain cone $\Xi^{\min} \subset \ZZ_{\ge 0}^\Sigma$ (see (\ref{eqn:Ximin})).

 We then follow the approach of \cite[\S8]{DS} to define partial $\Theta$-operators in~\S\ref{Theta}.
 For each $\tau = \tau_{\gp,i} \in \Sigma_0$, this gives a Hecke-equivariant operator
 $$\Theta_\tau:  M_{\bk,\bl}(U;\FF) \longrightarrow M_{\bk + \bh_\theta + 2\be_\theta ,\bl - \be_\theta}(U;\FF)$$
 where $\theta = \theta_{\gp,i,e_\gp}$.  Note in particular that {\em if $\gp$ is ramified, then the shift\footnote{This precise shift is
 predictable from the point of view of forthcoming work with Sasaki generalizing the geometric Serre weight conjectures
 of \cite{DS} to the ramified case.} in the weight parameter $\bk$
 is by $\be_{\sigma^{-1}\theta} + \be_{\theta}$}.   The idea of the construction, inspired by the one in \cite[\S12]{AG},
 is to divide by fundamental Hasse invariants to get a rational function on the Igusa cover of $\ol{Y}_U$, differentiate,
 project to the top graded piece of the $\tau$-component of the Kodaira--Spencer filtration, and finally multiply by
 fundamental Hasse invariants to descend to $\ol{Y}_U$ and eliminate poles.  The argument also gives a 
 direct (albeit local) definition of the $\Theta$-operator without reference to the Igusa cover in (\ref{eqn:altdeftheta}),
 and establishes the following result (Theorem~\ref{thm:theta}) generalizing \cite[Thm.~8.2.2]{DS}:
 \begin{ithm} \label{ithm:theta} Let $\tau = \tau_{\gp,i}$ and $\theta = \theta_{\gp,i,e_\gp}$.
 Then $\Theta_{\tau}(f)$ is divisible by $H_{\theta}$ if and only if either $f$ is divisible by $H_{\theta}$ or $p|k_{\theta}$.
\end{ithm}

We turn to the construction of partial Frobenius operators $V_\gp$ in~\S\ref{Frob}.   This essentially generalizes a definition
in \cite[\S9.8]{DS}, but requires significantly more work to actualize if $\gp$ is ramified.   We do this using Dieudonn\'e theory to define
a partial Frobenius endomorphism $\Phi_\gp$ of $\ol{Y}_U$ and an isomorphism $\Phi_\gp^*\ol{\CA}_{\bk,\bl} \cong \ol{\CA}_{\bk'',\bl''}$, where
$$\bk'' = \bk + \sum_{\theta \in \Sigma_\gp} k_\theta \bh_\theta\quad\mbox{and} \quad \bl'' = \bl + \sum_{\theta \in \Sigma_\gp} l_\theta \bh_\theta,$$
in order to obtain, for $\gp \in S_p$, commuting Hecke-equivariant operators
$$V_\gp: M_{\bk,\bl}(U;\FF) \longrightarrow M_{\bk'',\bl''}(U;\FF).$$

We will use $q$-expansions to relate the kernel of $\Theta_\tau$ to the image of $V_\gp$ for $\tau \in \Sigma_{\gp,0}$,
so we recall the theory in \S\ref{q1}.  This is a straightforward adaptation to our setting of results and methods developed
in \cite{rap, chai, dim, dem}.  In \S\ref{q2} we compute the (constant) $q$-expansions of the invariants $H_\theta$ and $G_\theta$
at each cusp of $\ol{Y}_U$, and we obtain formulas generalizing the classical ones for the effect of the operators $\Theta_\tau$
and $V_\gp$ on all $q$-expansions.  In particular this shows that the operators $\Theta_\tau$ for varying $\tau$ commute.

In \S\ref{sec:ker} we turn our attention to the description of the kernel of $\Theta_\tau$.
The $q$-expansion formulas also show that $\Theta_\tau \circ V_\gp = 0$ if $\tau \in \Sigma_{\gp,0}$, and that
$\ker(\Theta_\tau)$ is the same for all $\tau \in \Sigma_{\gp,0}$.  Theorem~\ref{ithm:theta} then reduces the study of
the kernel to the case of weights of the form $(\bk'',\bl'')$ where $\bk'',\bl''$ are as in the definition of $V_\gp$,
for which the argument proving~\cite[Thm.~9.8.2]{DS} gives the following:\footnote{This is a slight reformulation of Theorem~\ref{thm:kertheta}.}
\begin{ithm}  \label{ithm:ker}  If $\bk,\bl \in \ZZ^\Sigma$ and $\tau = \tau_{\gp,i}$ and $\theta = \theta_{\gp,i,e_\gp}$, 
then the sequence
$$0 \longrightarrow M_{\bk,\bl}(U,\FF) \stackrel{V_\gp}{\longrightarrow}  M_{\bk'',\bl''}(U;\FF) 
            \stackrel{\Theta_\tau}{\longrightarrow}  M_{\bk'' + \bh_\theta + 2\be_\theta,\bl'' - \be_\theta}(U;\FF)$$
is exact.
\end{ithm}
 
Before discussing the application to positivity of minimal weights, we remark that a less precise relation
among the weight-shifting operations can be neatly encapsulated in terms of the
algebra of modular forms of all weights
$$M_\tot(U;\FF) := \bigoplus_{\bk,\bl \in \ZZ^\Sigma}  M_{\bk,\bl}(U;\FF),$$
or even its direct limit $M_\tot(\FF) := \varinjlim_U M_\tot(U;\FF)$ (over all sufficiently small levels prime to $p$).
It follows from its definition that the operator $V_\gp$ (resp.~$\Theta_\tau$) on
the direct sum is an $\FF$-algebra endomorphism
(resp.~$\FF$-derivation) of $M_\tot(\FF)$.  One also finds that $V_\gp$ maps the ideal
$$\gI = \langle H_\theta' - 1, G_\theta' - 1\rangle_{\theta \in \Sigma}$$
to itself.\footnote{The $H_\theta'$ and $G_\theta'$ are slight modifications of the $H_\theta$ and $G_\theta$
obtained by rescaling those for which $j=1$; see \S\ref{sec:exact}.}  Furthermore 
$\Theta_\tau(H_\theta') = \Theta_\tau(G_\theta') = 0$ for all $\theta \in \Sigma$, so 
$V_\gp$ (resp.~$\Theta_\tau$) induces an $\FF$-algebra endomorphism (resp.~$\FF$-derivation)
of the quotient $M_\tot(\FF)/\gI$.   We then have the following consequence of 
Theorem~\ref{ithm:ker} (see Theorem~\ref{thm:exact}):
\begin{ithm}  If $\tau \in \Sigma_{\gp,0}$, then the sequence
$$0 \longrightarrow M_\tot(\FF)/\gI \stackrel{V_\gp}{\longrightarrow}
M_\tot(\FF)/\gI \stackrel{\Theta_\tau}{\longrightarrow}M_\tot(\FF)/\gI$$
is exact.
\end{ithm}

In \S\ref{sec:ptwt0} we apply our results to refine the main result of \cite{DK2}, which we recall states that minimal
weights of non-zero forms always lie in $ \Xi^{\min}$.   The geometric Serre weight conjectures of \cite{DS} (and its forthcoming
generalization to the ramified case) predict that if $f$ is a mod $p$ Hecke eigenform which is non-Eisenstein
(in the sense that the associated Galois representation is irreducible), then $\bk_{\min}(f)$ should be totally positive.
We use Theorems~\ref{ithm:theta} and~\ref{ithm:ker} to prove a partial result in this direction (Theorem~\ref{thm:ptwt0}):
\begin{ithm}  \label{ithm:ptwt0} Suppose that $\gp \in S_p$ is such that $F_\gp \neq \Q_p$ and $p^{f_\gp} > 3$.
Suppose that $f \in M_{\bk,\bl}(U;\F)$ is non-zero and $\bk = \bk_{\min}(f)$.  If $k_\theta = 0$ for some $\theta \in \Sigma_\gp$,
then $\bk = \bf{0}$.
\end{ithm}
Since the Hecke action on forms of weight $(\bf{0},\bl)$ is Eisenstein (see Proposition~\ref{prop:wt0}), the
theorem implies the total positivity of minimal weights of non-Eisenstein eigenforms in many situations, for
example if $p>3$ and there are no primes $\gp \in S_p$ such that $F_\gp = \Q_p$.  We remark that the
hypothesis $F_\gp \neq \Q_p$ cannot be removed from Theorem~\ref{ithm:ptwt0}:  if $\Sigma_\gp = \{\theta\}$,
then there are non-zero forms whose minimal weight $\bk$ satisfies $k_\theta = 0$ and $k_{\theta'} > 0$
for some $\theta' \neq \theta$.   However forthcoming work with Kassaei will show that the Hecke action
on such forms is Eisenstein; like in~\cite{DK, DK2}, the case of split primes seems to  require
a completely different method.   Unfortunately the case  of $p \le 3$, $f_\gp = 1$, $e_\gp > 1$ slips
through the crack between the two methods.  We do not know whether Theorem~\ref{ithm:ptwt0}
should hold in this case, but we still at least conjecture the failure is Eisenstein.

 \bigskip

\noindent {\bf Acknowledgements:}  The author would like to thank Payman Kassaei and Shu Sasaki for numerous
helpful conversations, particularly in the context of collaboration on forthcoming work that motivates
the results in this paper.  It is also a pleasure to acknowledge the evident debt this work owes to 
that of Reduzzi and Xiao in \cite{RX}.  The author is also grateful to Deo, Dimitrov 
and Wiese, whose interest in partial $\Theta$-operators (and their kernels) encouraged him to explore 
the ramifications of the approach in \cite{DS}, and to the referee for pointing out a number of typos and
useful references.

\section{Preliminaries} \label{prelim}
\subsection{Embeddings and decompositions} \label{sec:notation}
We first set out notation and conventions for various constructions associated to the set of embeddings
of a totally real field $F$, which together with a prime $p$, will be fixed throughout the paper.

We assume that $F$ has degree $d=[F:\QQ] > 1$, let $\CO_F$ denote its ring of integers, $\gd$ its different, and 
$\Sigma$ the set of embeddings $F \to \Qbar$, where $\Qbar$ is the algebraic closure of $\Q$ in $\C$.  

We also fix an embedding $\Qbar \to \Qpbar$.  We let $S_p$ denote the set of primes
of $\CO_F$ dividing $p$, and identify 
$\Sigma$ with $\coprod_{\gp \in S_p} \Sigma_{\gp}$ under the natural bijection, where 
$\Sigma_{\gp}$ denotes the set of embeddings $F_{\gp} \to \Qpbar$.

For each $\gp \in S_p$, we let $F_{\gp,0}$ denote the maximal unramified subextension of $F_\gp$,
which we identify with the field of fractions of $W(\CO_F/\gp)$.   We also let
$f_\gp$ denote the residue degree $[F_{\gp,0}:\Q_p]$,  $e_\gp$ the ramification index $[F_\gp:F_{\gp,0}]$,
and  $\Sigma_{\gp,0}$ the set of embeddings $F_{\gp,0} \to \Qpbar$, which we may identify with the
set of embeddings $\CO_F/\gp \to \Fpbar$, or homomorphisms $W(\CO_F/\gp) \to W(\Fpbar)$.
For each $\gp \in S_p$, we fix a choice of embedding $\tau_{\gp,0}  \in \Sigma_{\gp,0}$,
and for $i \in \Z/f_\gp\Z$, we let $\tau_{\gp,i} = \phi^i \circ \tau_{\gp,0}$ where $\phi$
is the Frobenius automorphism of $\Fpbar$ (or $W(\Fpbar)$ or its field of fractions), so that
$\Sigma_{\gp,0} = \{\tau_{\gp,1},\tau_{\gp,2},\ldots,\tau_{\gp,f_{\gp}}\}$.  We also
let $\Sigma_0 = \coprod_{\gp \in S_p} \Sigma_{\gp,0}$.  Letting $\gq = \prod_{\gp \in S_p}\gp$ denote the
radical of $p$ in $\CO_F$, note that $\Sigma_0$ may also be identified with the set of ring homomorphisms
$\CO_F/\gq \to \Fpbar$ (or indeed $\CO_F \to \Fpbar$).

For each 
$\tau = \tau_{\gp,i} \in \Sigma_0$, we let $\Sigma_\tau \subset \Sigma_\gp$ denote the set of embeddings restricting
to $\tau$, for which we choose an ordering $\theta_{\gp,i,1},\theta_{\gp,i,2},\ldots,\theta_{\gp,i,e_\gp}$, so that
$$\Sigma = \coprod_{\tau \in \Sigma_0} \Sigma_\tau = \{\,\theta_{\gp,i,j}\,|\,\gp \in S_p, i \in \ZZ/f_\gp\ZZ, 1 \le j \le e_\gp\,\}.$$
We also define a permutation $\sigma$ of $\Sigma$ whose restriction to each $\Sigma_\gp$ is the $e_\gp f _\gp$-cycle
corresponding to the right shift of indices with respect to the lexicographic ordering, i.e.,
$$\begin{array}{ccccccccc}
 (1,1)&  \mapsto &(1,2) &\mapsto &\cdots &\mapsto &(1,e_\gp) & \mapsto &\\
 (2,1) & \mapsto & (2,2) & \mapsto & \cdots &\mapsto &(2,e_\gp) &\mapsto & \\
 &&&&\vdots&&&&\\
 (f_\gp,1)& \mapsto&(f_\gp,2)&\mapsto & \cdots &\mapsto &(f_\gp,e_\gp) &\mapsto &(1,1) .\end{array}$$

Let $E \subset \Qbar$ be a number field containing the image of $\theta$ for all $\theta \in \Sigma$,
let $\CO$ be the completion of $\CO_E$ at the prime determined by the choice of
$\Qbar \to \Qpbar$, and let $\FF$ be its residue field.  For any $\CO_{F,p} = \CO_F \otimes \ZZ_p$-module $M$, we write 
$M = \bigoplus_{\gp \in S_p} M_\gp$ for the decomposition obtained from that of
$$\CO_{F,p} \cong \prod_{\gp \in S_p}  \CO_{F,\gp}.$$
Similarly for any $W(\CO_F/\gp) \otimes_{\ZZ_p} \CO$-module $N$, we have a
decomposition $N = \bigoplus_{\tau \in \Sigma_{\gp,0}} N_\tau$ obtained from
$$W(\CO_F/\gp) \otimes_{\ZZ_p} \CO  \cong \prod_{\tau \in \Sigma_{\gp,0}} \CO.$$
In particular, for any $\CO_F\otimes \CO$-module $M$, we have the decomposition
$$M = \bigoplus_{\gp \in S_p} M_\gp = \bigoplus_{\tau \in \Sigma_0} M_\tau,$$
where we simply write $M_\tau$ for $M_{\gp,\tau}$.  We also write $M_{\gp,i}$ for
$M_\tau$ if $\tau = \tau_{\gp,i}$; thus $M_{\gp,i}$ is the summand of 
the $\CO_{F,\gp}\otimes_{\ZZ_p} \CO$-module $M_\gp$ on which $W(\CO_F/\gp)$
acts via $\tau_{\gp,i}$.

We also fix a choice of uniformizer $\varpi_\gp$ for each $\gp \in S_p$.  We let $f_\gp(u)$
denote the minimal polynomial of $\varpi_\gp$ over $W(\CO_F/\gp)$, and let $f_\tau$ denote its image
in $\CO[u]$ for each $\tau \in \Sigma_{\gp,0}$; thus $u \mapsto \varpi_\gp \otimes 1$ induces an isomorphism 
$$\CO[u]/(f_\tau(u)) \stackrel{\sim}{\longrightarrow} \CO_{F,\gp} \otimes_{W(\CO_F/\gp),\tau} \CO.$$
Furthermore we have $f_\tau(u) = \prod_{\theta \in \Sigma_\tau}  (u - \theta(\varpi_\gp))$, and we
define elements
\begin{equation} \label{eqn:st} \begin{array}{rcc}
s_{\tau,j} &=& (u-\theta_{\gp,i,1}(\varpi_\gp))\cdots(u-\theta_{\gp,i,j}(\varpi_\gp))\\
\mbox{and}\quad t_{\tau,j} &=& (u-\theta_{\gp,i,j+1}(\varpi_\gp))\cdots(u-\theta_{\gp,i,e_\gp}(\varpi_\gp))\end{array}\end{equation}
of $\CO[u]/(f_\tau(u))$ for $j=0,\ldots,e_\gp$ (with the obvious convention that $s_{\tau,0} = t_{\tau,e_\gp} = 1$).
Note that each of the ideals $(s_{\tau,j})$ and $(t_{\tau,j})$ is the other's annihilator; furthermore the quotients
of $\CO[u]/(f_\tau(u))$ by these ideals are free over $\CO$, and the corresponding
ideals in $\CO_{F,\gp} \otimes_{W(\CO_F/\gp),\tau} \CO$ may be described as kernels of projection maps
to products of copies of $\CO$, hence depend only on $j$ and the ordering of embeddings, and not on the
choice of uniformizer $\varpi_\gp$.

For an invertible $\CO_F$-module $L$ and an embedding $\theta = \theta_{\gp,i,j} \in \Sigma_\tau$,
we define $L_\theta$ to be the free rank one $\CO$-module
 \begin{equation} \label{eqn:graded} L_\theta = t_{\tau,j}(L\otimes\CO)_\tau  \otimes_{\CO[u],\theta} \CO. \end{equation}
Note that $L_\theta$ is {\em not} to be identified with $L\otimes_{\CO_F,\theta} \CO$; rather there is a canonical map
$L_\theta \to L\otimes_{\CO_F,\theta} \CO$ which is an isomorphism if and only if $j = e_\gp$.
If $L$ and $L'$ are invertible $\CO_F$-modules, we will write $LL'$ for $L\otimes_{\CO_F}L'$ and
$L^{-1}$ for $\Hom_{\CO_F}(L,\CO_F)$.  Note that there are natural maps $L_\theta \otimes_{\CO} L'_\theta \to (LL')_\theta$
and $(L^{-1})_\theta \to \Hom_{\CO}(L_\theta,\CO)$, but again these are isomorphisms if and only if $j = e_\gp$.



 \subsection{Pappas--Rapoport models} \label{PR}
\label{sec:models}
In this section we recall the description of the Hilbert modular variety as a coarse moduli space for abelian varieties
with additional structure, along with the construction by Pappas and Rapoport of a smooth integral model 
(see~\cite{PR} and \cite{Shu}).

Let $G = \Res_{F/\QQ} \GL_2$ and let $U$ be an open compact subgroup of $\GL_2(\widehat{\CO}_F) \subset \GL_2(\A_{F,\f})=
G(\A_\f)$ of the form $U_pU^p$, where $U_p = \GL_2(\CO_{F,p})$ and $U^p \subset \GL_2(\A_{F,\f}^{(p)})$ is sufficiently small,
in a sense to be specified below.

We consider the functor which associates, to a locally Noetherian $\CO$-scheme $S$, the set of isomorphism classes of data
$(A,\iota,\lambda,\eta,\CF^\bullet)$, where:
\begin{itemize}
\item $s:A \to S$ is an abelian scheme of relative dimension $d$;
\item $\iota: \CO_F\to \End_S(A)$ is an embedding such that $(s_*\Omega^1_{A/S})_\gp$
is, locally on $S$, free of rank $e_\gp$ over $W(\CO_F/\gp)\otimes_{\ZZ_p} \CO_S $
for each $\gp \in S_p$;
\item $\lambda$ is an $\CO_F$-linear quasi-polarization of $A$ such that for each connected component
$S_i$ of $S$, $\lambda$ induces an isomorphism $\gc_i\gd \otimes_{\CO_F}  A_{S_i} \to A_{S_i}^\vee$
for some fractional ideal $\gc_i$ of $F$ prime to $p$;
\item $\eta$ is a level $U^p$ structure on $A$, 
i.e., for a choice of geometric point $\overline{s}_i$ on each connected component $S_i$ of $S$,
the data of a $\pi_1(S_i,\overline{s}_i)$-invariant $U^p$-orbit of
$\widehat{\CO}_F^{(p)} = \CO_F\otimes \widehat{\Z}^{(p)}$-linear isomorphisms\footnote{Note the
conventions in place with respect to the different, which are motivated by the point of view that
we wish to systematically trivialize modules defined by cohomological constructions.}
$$\eta_i :(\widehat{\CO}_F^{(p)})^2 \to \gd \otimes_{\CO_F} T^{(p)}(A_{\overline{s}_i}),$$
where $T^{(p)}$ denotes the product over $\ell \neq p$ of the $\ell$-adic Tate modules,
and $g \in U^p$ acts on $\eta_i$ by pre-composing with right multiplication by $g^{-1}$;
\item $\CF^\bullet$ is a collection of Pappas--Rapoport filtrations, i.e., for
each $\tau = \tau_{\gp,i} \in \Sigma_0$, an increasing filtration of 
$\CO_{F,\gp} \otimes_{W(\CO_F/\gp),\tau} \CO_S$-modules 
$$0 = \CF_\tau^{(0)} \subset \CF_\tau^{(1)} \subset \cdots 
   \subset \CF_\tau^{(e_\gp - 1)} \subset \CF_\tau^{(e_\gp)} 
   = (s_*\Omega_{A/S}^1)_\tau$$
such that for $j=1,\ldots,e_{\gp}$, the quotient
$${\CL}_{\gp,i,j}  :=  \CF_\tau^{(j)}/\CF_\tau^{(j-1)}$$
is a line bundle on $S$ on which $\CO_F$ acts via $\theta_{\gp,i,j}$.
\end{itemize}

The proof of \cite[Lemma~2.4.1]{DS} does not assume $p$ is unramified in $F$, and shows that if
$U^p$ is sufficiently small and $\alpha$ is an automorphism of a triple $(A,\iota,\eta)$ over a
connected scheme $S$, then $\alpha = \iota(\mu)$ for some $\mu \in U \cap \CO_F^\times$.
If we assume further that $-1 \not\in U\cap \CO_F^\times$, then it follows from standard arguments
that the functor above is representable by an infinite disjoint union of quasi-projective schemes over
$\CO$, which we denote by $\widetilde{Y}_U$, and the argument in the proof of \cite[Prop.~6]{Shu}
shows that $\widetilde{Y}_U$ is smooth of relative dimension $d$ over $\CO$.  Furthermore defining an action of
$\CO_{F,(p),+}^\times$ on $\widetilde{Y}_U$ by 
$$\nu\cdot(A,\iota,\lambda,\eta,\CF^\bullet) = (A,\iota,\nu\lambda,\eta,\CF^\bullet)$$
(as in \cite[\S2.1.3]{DKS}), we see that the resulting action of $\CO_{F,(p),+}^\times/(U\cap\CO_F^\times)^2$ is free
and the quotient is representable by a smooth quasi-projective scheme over $\CO$, which we denote by $Y_U$.

We also have a natural right action of $\GL_2(\A_{F,\f}^{(p)})$ on the inverse system of schemes $Y_U$
induced by pre-composing the level structure $\eta$ with right-multiplication by $g^{-1}$.  More precisely suppose that
$U_1$ and $U_2$ are as above (with $U_1^p$ and $U_2^p$ sufficiently small) and $g \in \GL_2(\A_{F,\f}^{(p)})$
is such that $g^{-1}U_1g \subset U_2$.  Letting $(A,\iota,\lambda,\eta,\CF^\bullet)$ denote the universal object
over $\widetilde{Y}_{U_1}$, there is a prime-to-$p$ quasi-isogeny $A \to A'$ of abelian varieties
with $\CO_F$-action inducing isomorphisms $\gd\otimes_{\CO_F} T^{(p)}(A'_{\overline{s}_i}) \cong \eta_i((\widehat{\CO}_{F}^{(p)})^2g^{-1})$
for each $i$, from which we obtain a level $U_2$-structure $\eta' = \eta\circ r_{g^{-1}}$ on $A'$ (where $r_{g^{-1}}$ denotes
right multiplication by $g^{-1}$).  Together with the other
data inherited from $A$, we obtain an object $(A',\iota',\lambda',\eta',\CF'^\bullet)$ corresponding to a 
morphism $\widetilde{\rho}_g: \widetilde{Y}_{U_1} \to \widetilde{Y}_{U_2}$ and descending to a morphism 
$Y_{U_1} \to Y_{U_2}$ which we denote $\rho_g$.  These morphisms satisfy the evident compatibility
$\rho_{g_2}\circ \rho_{g_1} = \rho_{g_1g_2}$ whenever $g_1^{-1}U_1g_1 \subset U_2$
and  $g_2^{-1}U_2g_2 \subset U_3$. 

Finally we remark that the schemes $Y_U$ define smooth integral models
over $\CO$ for the Hilbert modular varieties associated to the group $G$ (with the usual choice of Shimura datum),
and their generic fibres and resulting $\GL_2(\A_{F,\f}^{(p)})$-action may be identified with those obtained from a
system of canonical models.  In particular for any $\CO \to \CC$, we have isomorphisms
$$Y_U(\CC) \cong \GL_2(F)_+\backslash (\uhp^\Sigma \times \GL_2(\A_{F,\f})/U) 
 \cong  \GL_2(\CO_{F,(p)})_+\backslash (\uhp^\Sigma \times \GL_2(\A^{(p)}_{F,\f})/U^p)$$
compatible with the right action of $\GL_2(\A_{F,\f}^{(p)})$ on the inverse system,
and inducing a bijection between the set of geometric components of $Y_U$ and
$$\A_{F,\f}^\times/F_+^\times\det(U) \cong (\A_{F,\f}^{(p)})^\times/\CO_{F,(p),+}^\times\det(U^p)$$
These isomorphisms arise in turn from ones of the form
$$\widetilde{Y}_U(\CC) \cong  \SL_2(\CO_{F,(p)})\backslash (\uhp^\Sigma \times \GL_2(\A^{(p)}_{F,\f})/U^p),$$
under which the set of geometric components of $\widetilde{Y}_U$ is described by
$(\A_{F,\f}^{(p)})^\times/\det(U^p)$.

\section{Automorphic bundles}  \label{bundles}
\subsection{Pairings and duality} \label{sec:pairings}
Before introducing the line bundles whose sections define the automorphic forms of interest in the paper, we
present a  plethora of perfect pairings provided by Poincar\'e duality.

We fix a sufficiently small $U$ as in \S\ref{sec:models} and
consider the de Rham cohomology sheaves $\CH^1_\dr(A/S) = \RR^1s_*\Omega^\bullet_{A/S}$ on
the universal abelian scheme $A$ over $S = \widetilde{Y}_U$.  Recall that these sheaves are locally
free of rank two over $\CO_F\otimes \CO_S$.  Furthermore Poincar\'e duality and the polarization $\lambda$ induce an
$\CO_F\otimes \CO_S$-linear isomorphism
$$\begin{array}{ccccc}
\CH^1_\dr(A/S)& \stackrel{\sim}{\longrightarrow} &  \Shom_{\CO_S}(\CH^1_\dr(A^\vee/S),\CO_S)&
  \stackrel{\sim}{\longleftarrow}& \Shom_{\CO_S}(\CH^1_\dr(\gc\gd \otimes_{\CO_F} A)/S),\CO_S)\\
 &&&& \wr\parallel\\
 &&&& \Shom_{\CO_S}(\gd^{-1} \otimes_{\CO_F}\CH^1_\dr(A/S),\CO_S)\end{array}$$
 (where $\gc$ depends on the connected component of $S$ and disappears from the last expression since it
 is prime to $p$).  This in turn induces $\CO_{F,\gp} \otimes_{W(\CO_F/\gp),\tau} \CO_S$-linear isomorphisms
\begin{equation} \label{eqn:pairing} \CH^1_\dr(A/S)_\tau \cong \Shom_{\CO_S}(\gd_\gp^{-1} \otimes_{\CO_{F,\gp}} \CH^1_\dr(A/S)_\tau,\CO_S) \end{equation}
which we view as defining a perfect $\CO_S$-bilinear pairing 
$\langle \cdot,\cdot \rangle_\tau^0$ between $\CH_\tau := \CH^1_\dr(A/S)_\tau$ and $\gd_{\gp}^{-1}\otimes_{\CO_{F,\gp}} \CH_\tau$ for $\tau \in \Sigma_{\gp,0}$.
Furthermore the pairing is alternating in the sense that $\langle x, c \otimes y \rangle_\tau^0 = - \langle y, c \otimes x\rangle_\tau^0$
on sections.  Alternatively, we may apply the canonical $\CO_F\otimes R$-linear isomorphism
$$\Hom_{\CO_F\otimes R} (M, \CO_F \otimes R) \stackrel{\sim}{\longrightarrow} \Hom_R(\gd^{-1} \otimes_{\CO_F} M,R)$$
induced by the trace for any $\CO_F \otimes R$-module $M$ to obtain an $\CO_F \otimes \CO_S$-linear isomorphism
\begin{equation} \label{eqn:pairing2}\CH^1_\dr(A/S)   
\stackrel{\sim}{\longrightarrow} \Shom_{\CO_F\otimes\CO_S}(\gc^{-1} \otimes_{\CO_F}\CH^1_\dr(A/S),\CO_F \otimes \CO_S),\end{equation}
and hence a perfect alternating $\CO_{F,\gp} \otimes_{W(\CO_F/\gp),\tau} \CO_S$-bilinear pairing $\langle\cdot,\cdot\rangle_\tau$ on $\CH_\tau$.

Note that $\CH_\tau$ is locally free of rank two over $\CO_{F,\gp} \otimes_{W(\CO_F/\gp),\tau} \CO_S$, hence a
vector bundle of rank $2e_\gp$ over $\CO_S$.  Furthermore $(s_*\Omega^1_{A/S})_\tau$ is a subbundle  of $\CH_\tau$
of rank $e_\gp$, but is not locally free over $\CO_{F,\gp} \otimes_{W(\CO_F/\gp),\tau} \CO_S$ if $e_\gp > 1$
(in which case the failure is on a closed subscheme of codimension one), and more generally $\CF_\tau^{(j)}$
is a subbundle of rank $j$ for $j=0,1,\ldots,e_\gp$.

Recall that for $j=0,1,\ldots,e_\gp$, we defined (see (\ref{eqn:st})) elements $s_{\tau,j}$ and $t_{\tau,j}$ 
of the ring $\CO_{F,\gp} \otimes_{W(\CO_F/\gp),\tau} \CO \cong \CO[u]/(f_\tau(u))$, where
$f_\tau$ is the image under $\tau$ of the Eisenstein polynomial associated to our choice of uniformizer $\varpi_\gp$.
In the following, we shall fix $\tau$ and omit the subscripts $\tau$, $\gp$
and $i$ to disencumber the notation; we also write simply $W$ for $W(\CO_F/\gp)$.  
Note that since $\CH$ is locally free over $\CO_{F,\gp} \otimes_W \CO_S$ and
$\CF^{(j)}$ is annihilated by $s_j$, we have that $\CF^{(j)}$ is in fact a subbundle
of $t_j\CH$.

For a subsheaf $\CE \subset \CH$ of $\CO_{F,\gp} \otimes_W \CO_S$-submodules,
we define $\CE^\perp$ to be its orthogonal complement under the pairing $\langle\cdot,\cdot\rangle$,
i.e., the kernel of the morphism:
$$\CH \stackrel{\sim}{\longrightarrow} \Shom_{\CO_{F,\gp}\otimes_W\CO_S}(\CH,\CO_{F,\gp}\otimes_W\CO_S) 
\longrightarrow \Shom_{\CO_{F,\gp}\otimes_W\CO_S}(\CE,\CO_{F,\gp}\otimes_W\CO_S),$$
or equivalently the orthogonal complement of $\gd_\gp^{-1} \otimes_{\CO_{F,\gp}} \CE$ under
the pairing $\langle\cdot,\cdot\rangle^0$.  Note from the latter description
that if $\CE$ is an $\CO_S$-subbundle of $\CH$, then so is $\CE^\perp$.

\begin{lemma}  \label{lem:ortho}  We have the equality $(\CF^{(j)})^\perp  = t_j^{-1} \CF^{(j)}$
for $j=0,1,\ldots,e$,  where $t_j^{-1} \CF^{(j)} = \ker(\CH \stackrel{t_j}{\longrightarrow} t_j \CH \longrightarrow t_j\CH/\CF^{(j)})$
is the preimage sheaf of $\CF^{(j)}$ under $t_j$.
\end{lemma}
\begpf  We prove the lemma by induction on $j$, the case of $j=0$ being obvious.

Suppose then that $1 \le j \le e$ and that $(\CF^{(j-1)})^\perp  = t_{j-1}^{-1} \CF^{(j-1)}$.  Note that $(\CF^{(j)})^\perp$ and $t_j^{-1} \CF^{(j)}$
are both kernels of surjective morphisms from $\CH$ to vector bundles of rank $j$ on $S$, so each is a subbundle of rank $2e-j$, and
hence it suffices to prove the inclusion $t_j^{-1}\CF^{(j)} \subset (\CF^{(j)})^\perp$.  To do so, we may work locally on $S$, and assume
that $\CF^{(j)}(V) = \CF^{(j-1)}(V) \oplus R x_j$ where $V = \Spec R$ is a Noetherian open subscheme of $\widetilde{Y}_U$ such that
$\CH|_V$ is free over $\CO_{F,\gp} \otimes_W \CO_V$ and $x_j \in \CH(V)$ satisfies $(u-\theta_j(\varpi))x_j \in \CF^{(j-1)}(V)$.
In particular $x_j = t_jy_j$ for some $y_j \in t_{j-1}^{-1}\CF^{(j-1)}(V)$, so that
$$t_j^{-1}\CF^{(j)}(V) = t_j^{-1}\CF^{(j-1)}(V) \oplus Ry_j\quad\mbox{and}\quad
   (\CF^{(j)}(V))^\perp = (\CF^{(j-1)}(V))^\perp \cap (Rx_j)^\perp.$$

Note that if $w \in t_j^{-1}\CF^{(j)}(V)$, then 
$$t_{j-1}w = (u-\theta_j(\varpi)) t_j w \in (u-\theta_j(\varpi))\CF^{(j)}(V) \subset \CF^{(j-1)}(V),$$
so that $w \in t_{j-1}^{-1} \CF^{(j-1)}(V) = (\CF^{(j-1)})(V))^\perp$.  Furthermore if
$w \in  t_j^{-1}\CF^{(j-1)}(V)$, then
$$\langle  w,x_j \rangle = \langle w, t_jy_j \rangle = \langle t_jw, y_j \rangle$$
since $t_jw \in \CF^{(j-1)}$ and $y_j \in  t_{j-1}^{-1}\CF^{(j-1)}(V) =  (\CF^{(j-1)}(V))^\perp$.
Finally since the pairing is alternating, we have
$$\langle y_j , x_j \rangle = \langle y_j,t_jy_j \rangle = \langle t_jy_j,y_j \rangle = \langle x_j,y_j \rangle = - \langle y_j,x_j \rangle,$$
which implies that $\langle y_j,x_j \rangle = 0$ (since $\widetilde{Y}_U$ is flat over $\ZZ_2$ if $p=2$).

We have now shown that $t_j^{-1}\CF^{(j)}(V) \subset (\CF^{(j-1)}(V))^\perp \cap (Rx_j)^\perp = (\CF^{(j)}(V))^\perp$, as required.
\epf

We now define $\CG^{(j)} = (u-\theta_j(\varpi))^{-1}\CF^{(j-1)}$ for $j=1,\ldots,e$.  Thus $\CG^{(j)}$ is a rank $j+1$
subbundle of $\CH$, and we have inclusions of subbundles $\CF^{(j-1)} \subset \CF^{(j)} \subset \CG^{(j)}$, so that
${\CL}_j := \CF^{(j)}/\CF^{(j-1)}$ is a rank one subbundle of the rank two vector bundle
${\CP}_j := \CG^{(j)}/\CF^{(j-1)}$.  Furthermore all the inclusions are morphisms of 
$\CO_{F,\gp}\otimes_W \CO_S$-modules, and $\CO_F$ acts on ${\CP}_j$ via $\theta_j$.

Note that $\CG^{(j)}$ is annihilated by $s_j$, so that $\CG^{(j)} \subset t_j\CH$, and
we have $t_j^{-1}\CG^{(j)} = t_{j-1}^{-1}\CF^{(j-1)} = (\CF^{j-1})^\perp$ by Lemma~\ref{lem:ortho},
from which it follows also that 
$$t_j^{-1}\CF^{(j-1)} = t_j^{-1}(t_j^{-1}\CG^{(j)})^\perp = (\CG^{(j)})^\perp.$$
(The last equality can be seen by arguing locally on sections, or by noting that the diagram
$$\xymatrix{ \CH \ar[r]\ar[d] & \Shom_{\CO_{F,\gp}\otimes_W\CO_S} (\CH,\CO_{F,\gp}\otimes_W\CO_S)  \ar[r]\ar[d] & 
\Shom_{\CO_{F,\gp}\otimes_W\CO_S}(\CG^{(j)},\CO_{F,\gp}\otimes_W\CO_S) \ar[d] \\
 \CH \ar[r]& \Shom_{\CO_{F,\gp}\otimes_W\CO_S} (\CH,\CO_{F,\gp}\otimes_W\CO_S)  \ar[r]& 
 \Shom_{\CO_{F,\gp}\otimes_W\CO_S}(t_j^{-1}\CG^{(j)},\CO_{F,\gp}\otimes_W\CO_S)}$$
 commutes, where the left horizontal morphisms are defined by the pairing, the right by restriction, and 
 all the vertical morphisms by $t_j$.  The kernel of the composite along the top is $(\CG^{(j)})^\perp$, whereas
  $t_j^{-1}(t_j^{-1}\CG^{(j)})^\perp$ is the kernel of the composite along the left and bottom.  Since
  $t_j:t_j^{-1}\CG^{(j)} \to \CG^{(j)}$ is a surjective morphism of vector bundles, the leftmost vertical arrow
  is injective, so these kernels coincide.)
Therefore multiplication by $t_j$ defines an isomorphism $(\CF^{(j-1)})^\perp/(\CG^{(j)})^\perp \stackrel{\sim}{\longrightarrow}
\CG^{(j)}/\CF^{(j-1)}$, and composing its inverse with the isomorphism
$$(\CF^{(j-1)})^\perp/(\CG^{(j)})^\perp  \stackrel{\sim}{\longrightarrow} \Shom_{\CO_{F,\gp}\otimes_W\CO_S}(\CG^{(j)}/\CF^{(j-1)},\CO_{F,\gp}\otimes_W\CO_S)$$
induced by the pairing on $\CH$, we obtain an alternating $(\CO_{F,\gp} \otimes_W \CO_S)$-valued pairing $\langle \cdot , \cdot \rangle_j$
on ${\CP}_j := \CG^{(j)}/\CF^{(j-1)}$, whose description on sections is given in terms of (\ref{eqn:pairing2}) by
$$\langle t_j x, t_j y \rangle_j = \langle x, t_j y \rangle = \langle t_j x , y \rangle.$$
Note that since $\CO_{F,\gp}$ acts via $\tau_j$ on ${\CP}_j$, we in fact have the identification
$$\Shom_{\CO_{F,\gp}\otimes_W\CO_S}({\CP}_j,\CO_{F,\gp}\otimes_W\CO_S) =
\Shom_{\CO_S}({\CP}_j,\CI_j)$$
where $\CI_j$ is the sheaf of ideals, and trivial rank one $\CO_S$-subbundle, of $\CO_{F,\gp}\otimes_W\CO_S$ generated
by the global section $\prod_{j'\neq j}(u-\tau_{j'}(\varpi_\gp))$.   We thus obtain a trivialization of $\wedge^2_{\CO_S} {\CP}_j$
corresponding to a perfect $\CO_S$-valued pairing $\langle \cdot , \cdot \rangle_j^0$, which an unravelling of
definitions shows is given in terms of the original pairing $\langle \cdot , \cdot \rangle^0$ of (\ref{eqn:pairing}) by the formula
$$\langle t_j x , y \rangle_j^0 = \langle  x, f'(\varpi_{\pi_\gp})^{-1} \otimes y \rangle^0$$
on sections (where $f'$ is the derivative of the Eisenstein polynomial $f$).

\subsection{Automorphic line bundles} \label{sec:bundles}
Recall that the Pappas--Rapoport model $\widetilde{Y}_U$ is equipped with line bundles\footnote{For the moment, we continue
to suppress the fixed $\tau = \tau_{\gp,i} \in \Sigma_0$ from the notation.}  $\CL_j$, which we described in \S\ref{sec:pairings} as sub-bundles of
the rank two vector bundles $\CP_j$.  It is natural and convenient to consider also the twists of $\CL_j$ by powers of the determinant
bundle of $\CP_j$:
$${\CN}_j = \wedge^2_{\CO_S}{\CP}_j \cong {\CL}_j \otimes_{\CO_S} {\CM}_j,$$
where ${\CM}_j$ is the line bundle ${\CP}_j/{\CL}_j$.
Note that the pairing $\langle \cdot , \cdot \rangle^0_j$ defines an isomorphism ${\CM}_j\stackrel{\sim}{\longrightarrow} {\CL}_j^{-1}$
and a trivialization $\CO_S \stackrel{\sim}{\longrightarrow} {\CN}_j$ (which depends on the choice of $\varpi$).

As we will now consider these bundles for varying $\tau$, we resume writing the indicative subscripts; thus for
$\tau = \tau_{\gp,i}$, we will denote $\CG^{(j)}$ by $\CG_\tau^{(j)}$, ${\CP}_j$ by ${\CP}_{\gp,i,j}$,
and similarly for ${\CM}_j$ and ${\CN}_j$.  We also freely replace the subscript ``$\gp,i,j$''
by $\theta$, where $\theta = \theta_{\gp,i,j}$, so that for each $\theta \in \Sigma$, we have now defined a rank
two vector bundle ${\CP}_\theta$ and line bundles ${\CL}_\theta$, ${\CM}_\theta$, 
${\CN}_\theta$ on $S = \widetilde{Y}_U$, along with exact sequences
\begin{equation} \label{eqn:Hodge1}
 0 \longrightarrow {\CL}_\theta \longrightarrow {\CP}_\theta \longrightarrow {\CM}_\theta \longrightarrow 0 \end{equation}
and a trivialization of ${\CN}_\theta = \wedge^2_{\CO_S} {\CP}_\theta \cong \CL_\theta \otimes_{\CO_S} {\CM}_\theta$.
Furthermore the bundles ${\CP}_\theta$, ${\CL}_\theta$ and ${\CM}_\theta$ are $\CO_F\otimes \CO_S$-subquotients
of $\CH^1_\dr(A/S)$ on which $\CO_F$ acts via $\theta$.

Recall that we have an action of $\CO_{F,(p),+}^\times$ on $\widetilde{Y}_U$ defined by multiplication on the 
quasi-polarization.  In particular if $\nu \in \CO_{F,(p),+}^\times$, then the identification of $\nu^*A$ with $A$
induces an $\CO_F\otimes \CO_S$-linear isomorphism $\nu^*\CH^1_\dr(A/S) \cong \CH^1_\dr(A/S)$ under
which $\nu^*\CF^\bullet$ corresponds to $\CF^\bullet$, and we thus obtain isomorphisms 
$\alpha_\nu: \nu^*{\CP}_\theta \stackrel{\sim}{\longrightarrow} {\CP}_\theta$
compatible with (\ref{eqn:Hodge1}) and satisfying $\alpha_{\nu'\nu} = \alpha_\nu\circ \nu^*(\alpha_{\nu'})$ (for $\nu,\nu' \in \CO_{F,(p),+}^\times$).
Recall also that the action of $\CO_{F,(p),+}^\times$ on $\widetilde{Y}_U$ factors through $\CO_{F,(p),+}^\times/(U\cap \CO_F^\times)^2$,
the isomorphism $\uA \stackrel{\sim}{\longrightarrow} \nu^*\uA$ being defined by $\iota(\mu^{-1})$ if $\nu = \mu^2$ for $\mu \in U\cap \CO_F^\times$,
and one finds that the automorphism of ${\CP}_\theta$ obtained from $\alpha_\nu$ is multiplication by $\theta(\mu)$, so the
natural action of $\CO_{F,(p),+}^\times$ on the bundles fails to define descent data with respect to the cover $\widetilde{Y}_U \to Y_U$.
We do however obtain descent data after taking suitable tensor products or base-changes of these bundles, which we now consider.

For any $\CO$-algebra $R$, we will use $\cdot_R$ to denote the base-change to $R$ of an $\CO$-scheme $X$, as well as the pull-back to $X_R$ of a
quasi-coherent sheaf on $X$.   Let $\{\,\be_\theta\,|\,\theta\in \Sigma\,\}$ denote the standard basis of $\ZZ^\Sigma$.
For $\bk = \sum k_\theta\be_{\theta}$ and $\bl = \sum l_\theta\be_\theta \in \ZZ^\Sigma$, we define the line bundle
$$\widetilde{\CA}_{\bk,\bl} =  \bigotimes_{\theta\in \Sigma} \left( {\CL}_\theta^{\otimes k_\theta} \otimes {\CN}_\theta^{\otimes l_\theta} \right)
    \cong  \bigotimes_{\theta\in \Sigma} \left( {\CL}_\theta^{\otimes k_\theta + l_\theta} \otimes {\CM}_\theta^{\otimes l_\theta} \right)$$
on $S = \widetilde{Y}_U$, where all tensor products are over $\CO_S$.  For $\bn = \sum n_\theta\be_\theta \in \ZZ^\Sigma$, we
let $\chi_{\bn}:\CO_F^\times \to \CO^\times$ denote the character defined by $\chi_{\bn}(\mu) = \prod_{\theta} \theta(\mu)^{n_\theta}$,
and we let $\chi_{\bn,R}$ denote the associated $R^\times$-valued character.  If $\bk$, $\bl$, $R$ and $U$ are such that
$\chi_{\bk+2\bl,R}$ is trivial on $\CO_F^\times \cap U$, then the action of $\CO_{F,(p),+}^\times$
on $\widetilde{\CA}_{\bk,\bl,R}$ (over its action on $\widetilde{Y}_{U,R}$) factors through $\CO_{F,(p),+}^\times/(U\cap \CO_F^\times)^2$ and hence
defines descent data, in which case we denote the resulting line bundle on $Y_{U,R}$ by $\CA_{\bk,\bl,R}$.  
\begin{definition}\label{def:HMFs}  For $\bk$, $\bl$, $U$ and $R$ as above, we call $\CA_{\bk,\bl,R}$ the {\em automorphic line bundle}
of weight $(\bk,\bl)$ on $Y_{U,R}$, and we define the space of {\em Hilbert modular forms} of weight $(\bk,\bl)$ and level $U$ with
coefficients in $R$ to be
$$M_{\bk,\bl}(U;R) : =  H^0(Y_{U,R} ,\CA_{\bk,\bl,R}).$$
\end{definition}
We note some general situations in which this space is defined:
\begin{itemize}
\item The paritious setting: if $w=k_\theta + 2l_\theta$ is independent of $\theta$, then $\chi_{\bk+2\bl}(\mu) = \Nm_{F/\QQ}(\mu)^w = 1$
 for all $\mu \in U \cap \CO_F^\times$ (assuming only $U$ is small enough that $U\cap \CO_F^\times$ has no elements of norm $-1$ if $w$ is odd).
\item The mod $p$ setting:  if $R$ is any $\FF = \CO/\gm_\CO$-algebra and $U$ is sufficiently small that $\mu \equiv 1 \bmod \gp$ for all $\mu \in U \cap \CO_F^\times$
and $\gp \in S_p$, then $\theta(\mu) \equiv 1 \bmod \gm_\CO$ for all $\theta \in \Sigma$, so $\chi_{\bk+2\bl,\FF}$ is trivial on $U \cap \CO_F^\times$,
and hence so is $\chi_{\bk+2\bl,R}$.
\item The torsion setting:  if $R$ is an $\CO/p^N\CO$-algebra and $U$ is sufficiently small that $\mu \equiv 1 \bmod p^N\CO_F$ for all $\mu \in U \cap \CO_F^\times$,
so $\chi_{\bk+2\bl,\CO/p^N\CO}$ is trivial on $U \cap \CO_F^\times$, and hence so is $\chi_{\bk+2\bl,R}$.
\end{itemize}

We also have a natural left action of $\GL_2(\A_{F,\f}^{(p)})$ on the direct limit over $U$ of the spaces $M_{\bk,\bl}(U;R)$.
More precisely suppose that $U_1$ and $U_2$ are as above and $g \in \GL_2(\A_{F,\f}^{(p)})$
is such that $g^{-1}U_1g \subset U_2$, in which case recall that in \S\ref{sec:models} 
we defined a morphism $\widetilde{\rho}_g: \widetilde{Y}_{U_1} \to \widetilde{Y}_{U_2}$
descending to a morphism $\rho_g:Y_{U_1} \to Y_{U_2}$.   Furthermore the morphism $\widetilde{\rho}_g$ is obtained from a
prime-to-$p$ quasi-isogeny $A \to A'$ where $A$ is the universal abelian scheme over $\widetilde{Y}_{U_1}$ and $A'$ is the pull-back
of the universal abelian scheme.  We thus obtain isomorphisms $\widetilde{\rho}_g^*{\CP}_{\theta,2} \to {\CP}_{\theta,1}$
compatible with (\ref{eqn:Hodge1}) and the action of $\CO_{F,(p),+}^\times$ (augmenting the notation for the automorphic bundles on
$\widetilde{Y}_{U_i}$ and $Y_{U_i,R}$ with the subscript $i$).  Note that $U_1 \cap \CO_F^\times \subset U_2 \cap \CO_F^\times$, so if
$\widetilde{\CA}_{\bk,\bl,R,2}$ descends to $Y_{U_2,R}$, then $\widetilde{\CA}_{\bk,\bl,R,1}$ descends to $Y_{U_1,R}$, and we obtain
an isomorphism $\rho_g^*\CA_{\bk,\bl,R,2}  \cong \CA_{\bk,\bl,R,1}$.  We then define $[g] = [g]_{U_1,U_2}: M_{\bk,\bl}(U_2;R) \to M_{\bk,\bl}(U_1;R)$
as the composite
$$H^0(Y_{U_2}, \CA_{\bk,\bl,R,2})  \stackrel{\rho_g^*}{\longrightarrow} H^0(Y_{U_1}, \rho_g^*\CA_{\bk,\bl,R,2})
    \stackrel{\sim}{\longrightarrow} H^0(Y_{U_1}, \CA_{\bk,\bl,R,1}).$$
These maps satisfy the obvious compatibility, namely that 
$$[g_1]_{U_1,U_2} \circ [g_2]_{U_2,U_3} = [g_1g_2]_{U_1,U_3}$$
whenever $g_1^{-1}U_1g_1 \subset U_2$ and  $g_2^{-1}U_2g_2 \subset U_3$, and hence define an action of 
 $\GL_2(\A_{F,\f}^{(p)})$ on 
 \begin{equation}\label{eqn:inflevel}  M_{\bk,\bl}(R) := \varinjlim_U M_{\bk,\bl}(U;R)\end{equation}
 (where the limit is over sufficiently small open compact $U$
 with respect to the maps $[1]_{U_1,U_2}$).  For paritious $\bk$, $\bl$ and any choice of $\CO \to \CC$, we may identify the spaces
 $M_{\bk,\bl}(U;\CC)$ with those of holomorphic Hilbert modular forms, compatibly with the usual action 
 (up to normalization by a factor of $||\det g||$ depending on conventions).

Finally we remark that the action of $\nu \in \CO_{F,(p),+}^\times$ on the trivialization of ${\CN}_\theta$ is given by multiplication
by $\theta(\nu)$, so their products do not descend to trivializations of line bundles on $Y_{U,R}$.  However since the stabilizer
of each geometric connected component of $\widetilde{Y}_U$ is $\CO_{F,+}^\times \cap \det(U)$, we can obtain a (non-canonical) trivialization 
of $\CA_{{\mathbf{0}},\bl,R}$ as in \cite[Prop.~3.6.1]{DS}, provided $\chi_{\bl,R}$ is trivial on $\CO_{F,+}^\times \cap \det(U)$
and the geometric connected components of $Y_U$ are defined over $R$.  
Furthermore the same argument as in the proof of \cite[Lemma~4.5.1]{DS} shows the following:
\begin{proposition}  \label{prop:wt0}  If $p^N R = 0$, then the action of $\GL_2(\A_{F,\f}^{(p)})$ on  $M_{\bf{0},\bl}(R)$
factors through $\det: \GL_2(\A_{F,\f}^{(p)}) \to (\A_{F,\f}^{(p)})^\times$; furthermore as a representation of $(\A_{F,\f}^{(p)})^\times$,
$M_{\bf{0},\bl}(R)$ is isomorphic to the smooth induction of the character $\CO_{F,(p),+}^\times \to R^\times$ defined by
$\alpha \mapsto \prod_\theta \theta(\alpha)^{l_\theta}$.
\end{proposition}

\subsection{The Kodaira--Spencer filtration}  \label{sec:KS}
In this section we define a filtration on $\Omega^1_{Y_U/\CO}$ whose pieces are described by automorphic bundles
with weight components $k_\theta = 2$, $l_\theta = -1$.  The construction of the filtration is 
due to Reduzzi and Xiao (see~\cite[\S2.8]{RX}), but their presentation is complicated by the fact they wish to prove
smoothness simultaneously, and it obscures the fact that the bundles we denoted $\CG^{(j)}$ automatically
satisfy the orthogonality condition appearing in the definition of their counterparts in \cite{RX}.
We will show below that, with smoothness already established, one can give a more direct
conceptual description of the filtration and its properties.\footnote{The simultaneous treatment in \cite{RX} 
seems natural in view of the inherent overlap in the analysis of deformations needed for both results.
However the decision not to appeal to the results in \cite{Vol} and \cite{Shu} also makes reference to a perceived
minor gap in the proof of \cite[Prop~2.11]{Vol}; we found no such gap, nor is that result even needed, but we remark
that we made implicit use of \cite[Cor.~2.10]{Vol} when invoking the proof of \cite[Prop.~6]{Shu} to conclude
that $\widetilde{Y}_U$ is smooth over $\CO$.}
Furthermore in the case $p=2$, the argument in \cite{RX} appeals to a very general flatness assertion for divided power
envelopes for which we could not find a proof or reference, so it is not used here.  

\begin{theorem}[Reduzzi--Xiao] \label{thm:KSfil} There exists a decomposition 
$$\Omega^1_{Y_U/\CO} = \bigoplus_{\gp \in S_p} \bigoplus_{i = 1}^{f_\gp} \Omega^1_{Y_U/\CO,\gp,i},$$
together with an increasing filtration
$$\begin{array}{rcccl}
0 = \Fil^0 (\Omega^1_{Y_U/\CO,\gp,i} ) &\subset& \Fil^1 (\Omega^1_{Y_U/\CO,\gp,i}) &\subset& \cdots\\
&\subset & \Fil^{e_\gp - 1}( \Omega^1_{Y_U/\CO,\gp,i} )&
  \subset & \Fil^{e_\gp}( \Omega^1_{Y_U/\CO,\gp,i} ) =  \Omega^1_{Y_U/\CO,\gp,i}\end{array}$$
for each $\gp \in S_p$ and $i=1,\ldots,f_\gp$,  such that for each $j=1,\ldots,e_\gp$, $\gr^j (\Omega^1_{Y_U/\CO,\gp,i} )$ is
isomorphic to the automorphic bundle $\CA_{2\be_\theta,-\be_\theta,\CO}$, where $\theta =  \theta_{\gp,i,j}$.
\end{theorem}
\begpf  As usual, we first prove the analogous result for $S: = \widetilde{Y}_U$ and then descend to $Y_U$.

We let $\delta_0: S \stackrel{\sim}{\longrightarrow} \Delta \hookrightarrow Z_0$ denote the first infinitesimal thickening of the
diagonal embedding, and we view $\Omega^1_{S/\CO}$ as $\delta_0^*\CI$, where $\CI$ denotes
the sheaf of ideals defining $\Delta$ in $Z_0$.   Letting $s:A \to S$ denote the universal abelian scheme, the
transition maps for the crystal $R^1s_{\cris,*}\CO_{A/\Z_p}$ and canonical isomorphisms with de Rham
cohomology yield an $\CO_F \otimes \CO_{Z_0}$-linear isomorphism
$$\alpha:  p_0^*\CH^1_\dr(A/S)  \stackrel{\sim}{\longrightarrow}  q_0^*\CH^1_\dr(A/S)$$
extending the identity on $S \cong \Delta$, where $p_0,q_0: Z_0 \to S$ are the two projection maps
$Z_0 \to S$.  Since $\alpha$ is $\CO_F$-linear, it follows from the definition of $\CP_{\tau,1} = \CG_\tau^{(1)}$
that $\alpha$ restricts to an isomorphism
$$\alpha_{\tau,1}:  p_0^*\CP_{\tau,1}  \stackrel{\sim}{\longrightarrow}  q_0^*\CP_{\tau,1}$$
for each $\gp \in S_p$ and $\tau \in \Sigma_{\gp,0}$. Furthermore since the composite 
$$p_0^*\CL_{\tau,1}  \hookrightarrow p_0^*\CP_{\tau,1} \stackrel{\sim}{\longrightarrow}  q_0^*\CP_{\tau,1}
      \onto q_0^*\CM_{\tau,1}$$
has trivial pull-back to $S \cong \Delta$, it factors through a morphism
$$\delta_{0,*} \CL_{\tau,1}  = p_0^*\CL_{\tau,1} \otimes_{\CO_{Z_0}} (\CO_{Z_0}/\CI)
  \longrightarrow q_0^*\CM_{\tau,1} \otimes_{\CO_{Z_0}} \CI = \delta_{0,*} \CM_{\tau,1}\otimes_{\CO_{Z_0}} \CI,$$
and hence induces a morphism
$$\beta_{\tau,1}:  \delta_{0,*} (\CL_{\tau,1}\otimes_{\CO_S} \CM_{\tau,1}^{-1})  \longrightarrow  \CI.$$

We then define the sheaf of ideals $\CI_{\tau,1}$ on $Z_0$ to be the image of $\beta_{\tau,1}$,
and we let $Z_{\tau,1}$ denote the subscheme of $Z_0$ defined by $\CI_{\tau,1}$, and
$p_{\tau,1}$ and $q_{\tau,1}$ the resulting projection maps $Z_{\tau,1} \to S$.
By construction the pull-back of $\beta_{\tau,1}$ to $Z_{\tau,1}$ is trivial, and hence so is that of the
morphism $p_0^*\CL_{\tau,1} \to q_0^*\CM_{\tau,1}$, which implies that the pull-back of $\alpha$
maps $p_{\tau,1}^*\CL_{\tau,1} = p_{\tau,1}^*\CF_\tau^{(1)}$ isomorphically to 
$q_{\tau,1}^*\CL_{\tau,1} = q_{\tau,1}^*\CF_\tau^{(1)}$.  It follows from $\CO_F$-linearity that $\alpha$
induces an isomorphism $p_{\tau,1}^*\CG_\tau^{(2)} \stackrel{\sim}{\longrightarrow} q_{\tau,1}^*\CG_\tau^{(2)}$
(if $e_\gp > 1$), and hence an isomorphism
 $$\alpha_{\tau,2}:  p_{\tau,1}^*\CP_{\tau,2}  \stackrel{\sim}{\longrightarrow}  q_{\tau,1}^*\CP_{\tau,2}.$$
The same argument as above now yields a morphism
$$\beta_{\tau,2}:  \delta_{0,*} (\CL_{\tau,2}\otimes_{\CO_S} \CM_{\tau,2}^{-1})  \longrightarrow  \CI/\CI_{\tau,1},$$
whose image is that of a sheaf of ideals on $Z_0$ we denote by $\CI_{\tau,2}$.

Iterating the above construction thus yields, for each $\tau \in \Sigma_{\gp,0}$, a chain of sheaves of ideals
$$0 = \CI_{\tau,0} \subset \CI_{\tau,1} \subset \cdots \subset \CI_{\tau,e_\gp}$$
on $Z_0$ such that $\alpha$ induces 
\begin{itemize}
\item isomorphisms $p_{\tau,j}^*\CF_\tau^{(j)} \stackrel{\sim}{\longrightarrow} q_{\tau,j}^* \CF_\tau^{(j)}$,
\item and surjections $\delta_{0,*} (\CL_{\tau,j}\otimes_{\CO_S} \CM_{\tau,j}^{-1})  \onto  \CI_{\tau,j}/\CI_{\tau,j-1}$,
\end{itemize}
for $j=1,\ldots,e_\gp$, where $Z_{\tau,j}$ denotes the closed subscheme of $Z_0$ defined by $\CI_{\tau,j}$
and $p_{\tau,j}$, $q_{\tau,j}$ are the projections $Z_{\tau,j} \to S$.

Furthermore we claim that the map 
$$\bigoplus_{\gp \in S_p}  \bigoplus_{\tau \in \Sigma_{\gp,0}} \CI_{\tau,e_\gp}  \to \CI$$
is surjective.  Indeed let $\CJ$ denote the image and let $T$ denote the corresponding closed
subscheme of $Z_0$, so $T$ is the scheme-theoretic intersection of the $Z_{\tau,e_\gp}$, and
let\footnote{With apologies for the temporary dual use of $p$.} 
$p$, $q: T \to S$ denote the projection maps.  By construction $\alpha$ pulls back to an isomorphism 
$p^*\CH^1_{\dr}(A/S) \stackrel{\sim}{\longrightarrow} q^*\CH^1_\dr(A/S)$ under which 
$p^*\CF_\tau^{(j)} \stackrel{\sim}{\longrightarrow} q^*\CF_\tau^{(j)}$ for all $\tau$ and $j$.
In particular $t_*\Omega^1_{p^*A/T} = p^*(s_*\Omega^1_{A/S} ) \stackrel{\sim}{\longrightarrow} 
q^*(s_*\Omega^1_{A/S}) = u_*\Omega^1_{q^*A/T}$ (where $t:p^*A \to T$ and $u:q^*A \to T$ are
the structure morphisms), which the Grothendieck--Messing Theorem
implies is induced by an isomorphism $p^*A \cong q^*A$ of abelian schemes lifting the
identity over $S$.  Since the isomorphism respects the filtrations $\CF^\bullet$, and the lifts of
the universal auxiliary structures $\iota$, $\lambda$ and $\eta$ over $T$ are unique, it follows that
$p^*\underline{A} \cong q^*\underline{A}$, which means that $p=q \in S(T)$, so $T = \Delta$.

Now defining $\Omega^1_{S/\CO,\gp,i}  = \delta_0^*\CI_{\tau_{\gp,i},e_\gp}$ and
$\Fil^j( \Omega^1_{S/\CO,\gp,i} )= \delta_0^*\CI_{\tau_{\gp,i},j}$ for $\gp \in S_p$,
$1 \le i \le f_\gp$, $1 \le j \le e_\gp$, we obtain surjective morphisms
$$\CL_{\tau_{\gp,i},j}\otimes_{\CO_S} \CM_{\tau_{\gp,i},j}^{-1}  \onto  \gr^j (\Omega^1_{S/\CO,\gp,i} )
\qquad \mbox{and}\qquad \bigoplus_{\gp,i} \Omega^1_{S/\CO,\gp,i} \onto \Omega^1_{S/\CO}.$$
Since the $\CL_{\tau_{\gp,i},j}\otimes_{\CO_S} \CM_{\tau_{\gp,i},j}^{-1}$ are line bundles and
$\Omega^1_{S/\CO}$ is locally free of rank $d$, it follows that all the maps are isomorphisms.

Finally the constructions above are independent of the polarization $\lambda$, hence are
compatible with the action of $\nu \in \CO_{F,(p),+}^\times$ on $S = \widetilde{Y}_U$.  More precisely,
the pull-back of $\alpha$ via the diagonal map $(\nu,\nu)$ is compatible with the canonical isomorphism
$\nu^*\CH^1_\dr(A/S) \cong \CH^1_\dr(A/S)$ induced by the identification of $\nu^*A$ with $A$,
from which it follows easily that the morphisms in the construction of the filtration are invariant
under the action of $\CO_{F,(p),+}^\times$, hence descend to give the decomposition, filtrations
and isomorphisms in the statement of the theorem.
\epf

Let us also note the interpretation of the Kodaira--Spencer filtration in terms of tangent spaces.
For a closed point $y$ of $S$ corresponding to the data $\underline{A}_0 = (A_0,\iota_0,\lambda_0,\eta_0,\CF_0^\bullet)$
over a finite extension $k$ of the residue field of $\CO$, the fibre $T_y(S)$ of $\Shom_{\CO_S}(\Omega^1_{S/\CO},\CO_S)$
is canonically identified with the set of isomorphism classes of data $\underline{A}_1$ over $k[\epsilon]$ lifting $\underline{A}_0$,
and the decomposition and filtrations of the theorem yield dual decompositions of $T_y(S)$ into components $T_y(S)_{\tau}$
with decreasing filtrations $\Fil^j( T_y(S)_{\tau})$.  From the proof of the theorem one sees immediately that
$\bigoplus_{\tau} \Fil^{j_\tau}( T_y(S)_\tau)$ corresponds to the set of $(A_1,\iota_1,\lambda_1,\eta_1,\CF_1^\bullet)$ such that
$\CF_{1,\tau}^{(j)}$ is the image of $\CF_{0,\tau}^{(j)} \otimes_k k[\epsilon]$ for all $\tau$ and $j \le j_\tau$ under the
canonical isomorphism
$$H^1_\dr(A_1/k[\epsilon]) \cong H^1_\cris(A_0/k[\epsilon]) \cong H^1_\dr(A_0/k) \otimes_k k[\epsilon].$$

We note also that the theorem yields a canonical (Kodaira--Spencer) isomorphism
$$\Omega_{Y_U/\CO}^d  \cong \wedge_{\CO_{Y_U}}^d \Omega^1_{Y_U/\CO}  \cong \CA_{{\bf{2}},{\bf{-1}},\CO}$$
(writing $\sum m\be_\theta$ as $\bm$ for $m \in \ZZ$).
Furthermore the decomposition, filtrations and isomorphisms of the theorem (and hence also the
Kodaira--Spencer isomorphism) are Hecke-equivariant in the obvious sense.  More precisely the same
argument as for the compatibility with the $\CO_{F,(p),+}^\times$-action, but using the quasi-isogeny
in the construction of $\rho_g$, shows that if $U_1$, $U_2$ and $g \in \GL_2(\A_{F,\f}^{(p)})$ are such that $g^{-1}U_1g \subset U_2$,
then $\rho_g^* \Fil^j (\Omega^1_{Y_{U_2}/\CO,\gp,i})$ corresponds to $\Fil^j (\Omega^1_{Y_{U_1}/\CO,\gp,i})$ for all $\gp,i,j$
under the canonical isomorphism  $\rho_g^*\Omega^1_{Y_{U_2}/\CO} \stackrel{\sim}{\longrightarrow} \Omega^1_{Y_{U_1}/\CO}$,
and the resulting diagrams
$$\xymatrix{ \rho_g^*\CA_{2\be_\theta,-\be_\theta,\CO,2}  \ar[r]^{\sim}\ar[d]^{\wr} & \CA_{2\be_\theta,-\be_\theta,\CO,1} \ar[d]^{\wr} \\
\rho_g^* \gr^j (\Omega^1_{Y_{U_2}/\CO,\gp,i}) \ar[r]^{\sim} & \gr^j(\Omega^1_{Y_{U_1}/\CO,\gp,i}) }$$
commute (where the top arrow is defined in the discussion preceding (\ref{eqn:inflevel})).

\section{Partial Hasse invariants}  \label{H}
\subsection{Construction of $H_\theta$ and $G_\theta$}  \label{sec:Hasse}
We now recall the definition, due to Reduzzi and Xiao \cite{RX}, of generalized partial Hasse invariants on Pappas--Rapoport models.
These will be, for each $\theta = \theta_{\gp,i,j} \in \Sigma$, a Hilbert modular form $H_\theta$ of weight $(\bh_\theta,\bf{0})$
with coefficients in $\FF = \CO/\gm_\CO$, where $\bh_\theta: = n_\theta \be_{\sigma^{-1}\theta} - \be_\theta$,
with $n_\theta = p$ if $j=1$ and $n_\theta = 1$ if $j > 1$.  We also define below a(n in)variant $G_\theta$ of weight 
$(\bf{0},\bh_\theta)$.

We will now be working in the mod $p$ setting, so until further notice $S$ will denote $\widetilde{Y}_{U,\FF}$, and
$s:A \to S$ the universal abelian scheme over it.  Thus $\CH^1_\dr(A/S)$ is a locally free sheaf of rank two over
$$\CO_F\otimes \CO_S \cong \bigoplus_{\gp \in S_p}  \bigoplus_{\tau \in \Sigma_{\gp,0}} \CO_S [u] / (u^{e_\gp}),$$
where $u$ acts via $\iota(\varpi_\gp)^*$ on the $\gp$-component of 
$$\CH^1_\dr(A/S) = \bigoplus_{\gp \in S_p} \CH_\gp = \bigoplus_{\gp \in S_p}  \bigoplus_{\tau \in \Sigma_{\gp,0}} \CH_\tau.$$

We will also now be working with a fixed $\gp$ and omit the subscript from the notation, so that
$$\CH =  \CH^1_\dr(A/S)_\gp = \bigoplus _{\tau \in \Sigma_{\gp,0}} \CH_\tau = \bigoplus_{i \in \ZZ/f\ZZ} \CH_i $$
with each $\CH_i$ locally free of rank two over $\CO_S[u]/(u^e)$ (where we have also abbreviated the
subscript $\tau_{\gp,i}$ by $i$).  Furthermore for each $i\in \ZZ/f\ZZ$, we have a filtration
$$0 = \CF_i^{(0)} \subset \CF_i^{(1)}\subset \cdots \subset \CF_i^{(e-1)} \subset \CF_i^{(e)} = (s_*\Omega^1_{A/S})_i$$
by sheaves of $\CO_S[u]/(u^e)$-modules 
such that the quotients $\CL_{i,j} = \CF_i^{(j)}/\CF_i^{(j-1)}$ are line bundles annihilated by $u$.

Firstly note that if $j > 1$, then $u: \CF_i^{(j)} \to \CF_i^{(j-1)}$ induces a morphism $\CL_{i,j} \to \CL_{i,j-1}$.
On the other hand if $j=1$, then the Verschiebung morphism $\phi_S^*A \to A$ over $S$ induces $\CO_S[u]/(u^e)$-linear morphisms
$$\Ver_i^*:  \CH_i = \CH_{\tau_i} \longrightarrow \CH^1_\dr(\phi_S^*(A)/S)_{\tau_i} \cong  \phi_S^*(\CH_{\tau_{i-1}}) = \phi_S^*(\CH_{i-1})$$
with image $\phi_S^*(\CF_{i-1}^{(e)})$ (where $\phi_S$ denotes the absolute Frobenius on $S$).
Note that $\CL_{i,1} = \CF_i^{(1)} \subset u^{e-1} \CH_i$, so that $u^{e-1}$ defines an isomorphism
$$u^{1-e}\CL_{i,1}/u \CH_i \stackrel{\sim}{\longrightarrow} \CF_i^{(1)},$$
and that $\Ver_i^*(u\CH_i) = u\phi_S^*(\CF_{i-1}^{(e)}) \subset \phi_S^*(\CF_{i-1}^{(e-1)})$,
so we obtain a well-defined $\CO_S$-linear morphism ``$\Ver_i^*\circ u^{1-e}$''
$$\CL_{i,1} \stackrel{\sim}{\longleftarrow} u^{1-e}\CL_{i,1} /u \CH_i  \longrightarrow 
  \phi_S^*( \CL_{i-1,e}) \cong \CL_{i-1,e}^{\otimes p}.$$
We have now defined a morphism $\CL_\theta \longrightarrow \CL_{\sigma^{-1}\theta}^{\otimes n_\theta}$
for all $\theta$, and hence a section of $\widetilde{\CA}_{\bh_\theta,\bf{0},\FF}$ over $S$.  Furthermore it
is straightforward to check that the section is invariant under the action of $\CO_{F,(p),+}^\times$ and
therefore descends to an element
$$ H_\theta  \in M_{\bh_\theta,\bf{0}} (U; \FF) = H^0(Y_{U,\FF} , \CA_{\bh_\theta,\bf{0},\FF}),$$
which we call the {\em partial Hasse invariant} (indexed by $\theta$).  Furthermore the partial Hasse
invariants are stable under the Hecke action, in the sense that if
$U_1$, $U_2$ and $g \in \GL_2(\A_{F,\f}^{(p)})$ are such that $g^{-1}U_1g \subset U_2$,
then $[g]H_{\theta ,2} = H_{\theta,1}$. (Note also
that the partial Hasse invariants are dependent on the choice of uniformizer $\varpi = \varpi_\gp$
only up to a scalar in $\FF^\times$: if $\varpi$ is replaced by $a \varpi$ for some $a \in \CO_{F,\gp}^\times$,
then $H_\theta$ is replaced by $\tau(a) H_\theta$ if $j >1$ and by $\tau(a)^{1-e} H_\theta$ if $j=1$.)

We remark also that the line bundles $\CA_{\bf{0},\bh_\theta,\FF}$ have canonical trivializations.
Indeed for each $i \in  \ZZ/f\ZZ$ and $j=2,\ldots,e$, we have the exact sequence
$$0 \to \CG_i^{(j-1)}/\CF_i^{(j-1)}  \longrightarrow \CG_i^{(j)}/\CF_i^{(j-1)} \stackrel{u}{\longrightarrow} \CF_i^{(j-1)}/\CF_i^{(j-2)} \to 0$$
over $S$, i.e., $0 \to \CM_{i,j-1} \to \CP_{i,j} \to \CN_{i,j-1} \to 0$, inducing an isomorphism
$$ \CN_{i,j-1} = \CL_{i,j-1}\otimes_{\CO_S} \CM_{i,j-1} \cong \wedge^2_{\CO_S} \CP_{i,j} = \CN_{i,j}$$
and hence $\CO_S \simeq \CN_{i,j-1}^{-1}\otimes_{\CO_S}\CN_{i,j} = \widetilde{\CA}_{\bf{0},\bh_\theta,\FF}$ 
for $\theta = \theta_{\gp,i,j}$, which it is straightforward to check descends to $Y_{U,\FF}$.   Similarly we have the exact sequence
$$0 \to\phi_S^* (\CG_{i-1}^{(e)}/\CF_{i-1}^{(e)})  \stackrel{\Frob_A^*}{\longrightarrow}  \CG_{i}^{(1)} 
 \stackrel{\Ver_A^*u^{1-e}}{\longrightarrow}  \phi_S^*(\CF_{i-1}^{(e)}/\CF_{i-1}^{(e-1)}) \to 0$$
inducing an isomorphism $\phi_S^*(\CN_{i-1,e})  \cong \CN_{i,1}$ and hence $\CO_S \simeq \widetilde{\CA}_{\bf{0},\bh_\theta,\FF}$ 
for $\theta = \theta_{\gp,i,1}$ descending to $Y_{U,\FF}$.  Furthermore these isomorphisms
are Hecke-equivariant in the usual sense, but note that they depend via $u$ on the choice of $\varpi_\gp$.
For each $\theta$, we let $G_\theta \in M_{\bf{0},\bh_\theta}(U;\FF)$ denote the canonical trivializing section.

\subsection{Stratification}  \label{sec:stratification}  We also recall how the partial Hasse invariants define a
stratification of the Hilbert modular variety in characteristic $p$.  For any $\theta \in \Theta$, we define
$\widetilde{Z}_\theta$ (resp.~$Z_\theta$) to be the closed subscheme of $S = \widetilde{Y}_{U,\FF}$
(resp.~$Y_{U,\FF}$) defined by the vanishing of $H_\theta$, and for any subset $T \subset \Sigma$, we let 
$$\widetilde{Z}_T = \bigcap_{\theta \in T} \widetilde{Z}_\theta\quad \mbox{and} \quad Z_T = \bigcap_{\theta \in T} Z_\theta.$$
Note that the schemes $Z_T$ are stable under the Hecke action, in the strong sense that $Z_{T,1}$ is the pull-back
of $Z_{T,2}$ under $\rho_g:Y_{U_1} \to Y_{U_2}$.

We then have the following consequence (\cite[Thm.~3.10]{RX}) of the description of the Kodaira--Spencer filtration on
tangent spaces at closed points; we give a proof here as some of the details are relevant to the construction of 
$\Theta$-operators in \S\ref{sec:theta}.
\begin{proposition} \label{prop:smooth}  The schemes $\widetilde{Z}_T$ and $Z_T$ are smooth over $\FF$ of dimension $|\Sigma - T|$.
\end{proposition}
\begpf  We prove the result for $\widetilde{Z}_T$, from which the result for $Z_T$ is immediate.

Let $y$ be a closed point of $S$ with local ring $R = \CO_{S,y}$, maximal ideal $\gm$
and residue field $k = R/\gm$.  For each $\theta \in \Sigma$, choose a basis $b_\theta$ for ${\CL}_{\theta,y}$
over $R$ and write $H_{\theta,y} b_\theta = x_\theta b_{\sigma^{-1}\theta}^{n_\theta}$.
Thus if $y \in \widetilde{Z}_\theta$, then $x_\theta \in \gm$, and we let $\overline{x}_\theta$ denote its
image in $\gm/\gm^2$.

Identifying $\gm/\gm^2$ with the fibre of $\Omega^1_{\widetilde{Y}_U/\CO}$ at $y$ and writing 
$\Fil^j(\gm/\gm^2)_\tau$ for the subspaces obtained from the Kodaira--Spencer filtration, we claim
that if $y \in Z_\theta$, then
\begin{equation} \label{eqn:derivs}
 \Fil^j(\gm/\gm^2)_\tau = k\overline{x}_\theta +  \Fil^{j-1}(\gm/\gm^2)_\tau, \end{equation}
where $\tau = \tau_{\gp,i}$ and $\theta = \theta_{\gp,i,j}$.
Comparing dimensions, we see it suffices to prove the inclusion of the left-hand side
in the right, or equivalently that if
$$ v \in T_y(S) = \bigoplus_{\tau' \in \Sigma_0} T_y(S)_{\tau'}$$
is such that its $\tau$-component $v_\tau$ lies in $\Fil^{j-1} T_y(S)_\tau$ and
$v$ is orthogonal to $\overline{x}_\theta$, then in fact $v_\tau \in \Fil^j(T_y(S)_\tau)$
(using the notation of the discussion following the proof of
Theorem~\ref{thm:KSfil}).  

Let $\underline{A}_0 = (A_0,\iota_0,\lambda_0,\eta_0,\CF_0^\bullet)$ denote the data corresponding to the
point $y \in S(k)$ and $\underline{A}_1 = (A_1,\iota_1,\lambda_1,\eta_1,\CF_1^\bullet)$ that of its lift $v \in S(k[\epsilon])$.
With $\tau = \tau_{\gp,i}$ fixed for now, we will suppress $\gp$ from the notation and replace the subscript $\tau_{p,i'}$
by $i'$ (for $i' = i,i-1$).  Recall the assumption that $v_i \in \Fil^{j-1}( T_y(S)_i)$ means that 
$\CF_{1,i}^{(j')}$ corresponds to $\CF_{0,i}^{(j')} \otimes_k k[\epsilon]$ for $j' = 1,\ldots,j-1$ under the
canonical isomorphism
\begin{equation}
\label{eqn:cris}
H^1_\dr(A_1/k[\epsilon]) \cong H^1_\cris(A_0/k[\epsilon]) \cong H^1_\dr(A_0/k) \otimes_k k[\epsilon].
\end{equation}
For $v_i$ to be orthogonal to $\overline{x}_\theta$ means that the morphism
\begin{equation}
\label{eqn:hasse}
 \CL_{1,\theta} \longrightarrow \CL_{1,\sigma^{-1}\theta}^{\otimes n_\theta}
 \end{equation}
induced by $H_\theta$ vanishes, and we need to show this implies that
$\CF_{1,i}^{(j)}$ is the image of $\CF_{0,i}^{(j)} \otimes_k k[\epsilon]$.

Suppose first that $j > 1$.  Then (\ref{eqn:hasse}) is simply
$$u: \CF_{1,i}^{(j)}/\CF_{1,i}^{(j-1)} \longrightarrow \CF_{1,i}^{(j-1)}/\CF_{1,i}^{(j-2)},$$
whose vanishing means $\CF_{1,i}^{(j)} = u^{-1}\CF_{1,i}^{(j-2)}$.
Since (\ref{eqn:cris}) sends $\CF_{1,i}^{(j-2)}$ to  $\CF_{0,i}^{(j-2)} \otimes_k k[\epsilon]$
and is compatible with $u$, it follows that it also sends 
$\CF_{1,i}^{(j)}$ to  $\CF_{0,i}^{(j)} \otimes_k k[\epsilon]$.

On the other hand if $j=1$, then the vanishing of (\ref{eqn:hasse}) means that
$u^{1-e}\CF_{1,i}^{(1)}$ is the preimage of 
 $\phi_1^*(\CF_{1,i-1}^{(e-1)})$ under $\Ver_i^*$
(where $\phi_1$ is the absolute Frobenius on $k[\epsilon]$, and $\phi_0$ will
denote the absolute Frobenius on $k$).  
Since the diagram
$$\xymatrix{H^1_\dr(A_1/k[\epsilon]) \ar[r] \ar[d] &  H^1_\dr(A_0/k) \otimes_k k[\epsilon] \ar[d]\\
\phi_1^*H^1_\dr(A_1/k[\epsilon]) \ar[r] & \phi_0^*H^1_\dr(A_0/k) \otimes_k k[\epsilon]}$$
commutes, where the vertical maps are induced by Verschiebung, the top arrow
is (\ref{eqn:cris}) and the bottom one is given by the identification of $\phi_1^*A_1$
with $\phi_0^*A_0 \otimes_k k[\epsilon]$, so in particular identifies 
$\phi_1^*(\CF_{1,i-1}^{(e-1)})$ with $\phi_0^*(\CF_{0,i-1}^{(e-1)}) \otimes_k k[\epsilon]$,
it follows that the top arrow sends  
$u^{1-e}\CF_{1,i}^{(1)}$ to  $u^{1-e}\CF_{0,i}^{(1)} \otimes_k k[\epsilon]$,
and hence $\CF_{1,i}^{(1)}$ to  $\CF_{0,i}^{(1)} \otimes_k k[\epsilon]$.
This completes the proof of the claim.

Now note that if $y \in \widetilde{Z}_T$, then (\ref{eqn:derivs}) implies that the elements $\overline{x}_\theta$ for $\theta \in T$ can
be extended to a basis for $\gm/\gm^2$ over $k$, hence are linearly independent.  Since $R$ is regular of dimension
$d = |\Sigma|$, it follows that $\CO_{\widetilde{Z}_T,y} = R/\langle x_\theta \rangle_{\theta \in T}$ is regular of dimension
of $d - |T|$, and hence that $\widetilde{Z}_T$ is smooth over $\FF$ of dimension $d - |T|$.
\epf

Finally we recall the definition of the minimal weight of a non-zero mod $p$ Hilbert modular form.
If $f \in M_{\bk,\bl}(U;\FF)$, then $\bk_{\min}(f)$ is defined to be $\bk - \sum_{\theta} m_\theta \bh_\theta$
where $\sum_\theta m_\theta \be_\theta$ is the unique maximal element of the set
$$\left\{\,\left. \sum_\theta m_\theta \be_\theta \in \ZZ_{\ge 0}^\Sigma\,\right|\, 
        f = g\prod_{\theta \in \Sigma} H_\theta^{m_\theta} \mbox{\,for some\,}
        g \in M_{\bk - \sum_\theta m_\theta \bh_\theta,\bl}(U;\FF)\,\right\}.$$
By the main result of \cite{DK2}, the minimal weight of $f$ always lies in the cone:
\begin{equation}\label{eqn:Ximin}
 \Xi^{\min} := \left\{\, \left.\sum_\theta k_\theta \be_\theta\,\right|\, \mbox{$n_\theta k_\theta \ge k_{\sigma^{-1}\theta}$ for all $\theta \in \Sigma$}\,\right\}.
 \end{equation}
Note that the result stated in \cite{DK2} applies to forms on a finite \'etale cover of $Y_{U,\FF}$, from which the
analogous result for forms on $Y_{U,\FF}$ is immediate.

\section{Partial Theta operators}    \label{Theta}
 \subsection{Fundamental Hasse invariants}   \label{sec:Igusa}
In order to define the partial $\Theta$-operators (in \S\ref{sec:theta} below), we first define
a canonical factorization of the partial Hasse invariants over a finite flat (Igusa) cover of the
Hilbert modular variety over $\FF$.

We fix a sufficiently small $U$ that the line bundles $\CL_\theta$, $\CM_\theta$, $\CN_\theta$ (and hence $\widetilde{\CA}_{\bk,\bl,\FF}$)
on $\widetilde{Y}_{U,\FF}$ descend to $Y_{U,\FF}$ for all $\theta \in \Sigma$ (and all $\bk,\bl\in \ZZ^\Sigma$), and we write simply
$\ol{Y}$ for $Y_{U,\FF}$, and $\ol{\CL}_{\tau,j}$, $\ol{\CM}_{\tau,j}$ and $\ol{\CN}_{\tau,j}$ for the line bundles on $\ol{Y}$.
For each $\gp \in S_p$ and $\tau \in \Sigma_{\gp,0}$, we let 
$$H_\tau = \prod_{j=1}^{e_\gp} H_{\tau,j} \in H^0(\ol{Y}, \ol{\CL}_{\tau,e_\gp}^{-1}\otimes_{\CO_{\ol{Y}}} \ol{\CL}_{\phi^{-1}\circ\tau,e_\gp}^{\otimes p}).$$
Viewing each $H_\tau$ as a morphism $(\ol{\CL}_{\phi^{-1}\circ\tau,e_\gp}^{-1})^{\otimes p} \to \ol{\CL}_{\tau,e_\gp}^{-1}$ and
$H_\gp  := \prod_{\tau \in \Sigma_{\gp,0}}  H_\tau$ as a morphism 
$\otimes_{\tau \in \Sigma_{\gp,0}} (\ol{\CL}_{\tau,e_\gp}^{-1})^{\otimes (p-1)}  \to \CO_{\ol{Y}}$,
we define the Igusa cover\footnote{The cover has a natural moduli-theoretic interpretation in terms of
$A[\gq]$, but we will not need this here.}  of $\ol{Y}$ (of level $\gq = \prod \gp$) to be
$${Y}^{\mathrm{Ig}}  =  \mathbf{Spec}  \left(\mathrm{Sym}_{\mathcal{O}_{\ol{Y}}}
(\bigoplus_{\tau \in \Sigma_0}  \ol{\CL}_{\tau,e_\gp}^{-1})/\mathcal{I}\right),$$
where $\CI$ is the sheaf of ideals generated by the $\CO_{\ol{Y}}$-submodules
$$\mbox{$(H_\tau - 1) \ol{\CL}_{\tau,e_\gp}^{-1}$ for $\tau \in \Sigma_{0}$, and
 $(H_\gp -1)\left(\bigotimes_{\tau \in \Sigma_{\gp,0}} (\ol{\CL}_{\tau,e_\gp}^{-1})^{\otimes (p-1)}\right)$ for $\gp \in S_p$}$$
 (where all tensor products are over $\CO_{\ol{Y}}$).
 We then define an action of $(\CO_F/\gq)^\times$ on ${Y}^{\Ig}$ over $\ol{Y}$ by letting
 $\alpha \in (\CO_F/\gq)^\times$ act on the structure sheaf by the automorphism of sheaves
 of $\CO_{\ol{Y}}$-algebras induced by multiplication by $\tau({\alpha})^{-1}$ on $\ol{\CL}_{\tau,e_\gp}^{-1}$
 for each $\tau$.  We then see, exactly as in the proof of parts (1) and (2) of \cite[Prop.~8.1.1]{DS},
 that the projection $\pi:{Y}^{\Ig} \to \ol{Y}$ is finite flat, generically \'etale, and identifies $\ol{Y}$
 with the quotient of ${Y}^{\Ig}$ by the action of $(\CO_F/\gq)^\times$.  
 
For each $\tau \in \Sigma_{\gp,0}$, we let $h_{\tau,e_\gp}$ denote the tautological section of $\pi^*\ol{\CL}_{\tau,e_\gp}$
induced by the inclusion $\ol{\CL}_{\tau,e_\gp}^{-1} \hookrightarrow \pi_*\CO_{{Y}^\Ig}$.  We also define the section
$$h_{\tau,j} = \pi^*(H_{\tau,j+1}\cdots H_{\tau,e_\gp})h_{\tau,e_\gp}$$
of $\pi^*\ol{\CL}_{\tau,j}$ for $j=1,\ldots,e_{\gp}-1$.
Note that since ${Y}^\Ig$ is reduced (or since $\prod_{\tau \in \Sigma_{\gp,0}} h_{\tau,e_\gp}^{p-1} = \pi^*H_\gp$
by construction), the sections $h_{\tau,e_\gp}$ are injective, and hence so are the $h_{\tau,j}$ for all $\tau$ and $j$.
We write simply $h_\theta$ for the section $h_{\tau,j}$ of $\pi^*\ol{\CL}_\theta = \pi^*\ol{\CL}_{\tau,j}$, and we call
the $h_\theta$ the {\em fundamental Hasse invariant} (indexed by $\theta$).  

\subsection{Construction of $\Theta_\tau$}   \label{sec:theta}
We now explain how the construction of $\Theta$-operators in \cite{DS} directly generalizes to the case
where $p$ is ramified in $F$, yielding an operator that shifts the weight $\bk$ by $(1,1)$ in the final two components
corresponding to embeddings with the same reduction, i.e., $\theta_{\gp,i,e_{\gp}-1}$, $\theta_{\gp,i,e_\gp}$
(and hence, by composing with multiplication by partial Hasse invariants, one can shift weights by $+1$ for
 any pair of embeddings with the same reduction).

Indeed for each $\tau \in \Sigma_0$, we define the operator $\Theta_{\tau}$  
{\em exactly} as in \cite[\S8]{DS}, but using the morphism
$$KS_{\tau}: \Omega^1_{\overline{Y}/\F} \longrightarrow 
\gr^{e_{\gp}}(\Omega^1_{\overline{Y}/\F})_{\tau} \stackrel{\sim}{\longrightarrow}
\ol{\CL}_{\tau,e_\gp}\otimes_{\CO_{\ol{Y}}} \ol{\CM}^{-1}_{\tau,e_\gp}$$
provided by Theorem~\ref{thm:KSfil} via projection to the top graded piece
of the filtration of the $\tau$-component of $\Omega^1_{Y/\F}$.  More precisely,
fix $\gp_0 \in S_p$ and $\tau_0 = \tau_{\gp_0,i}$, let $\theta_0 = \theta_{\gp_0,i,e_{\gp_0}}$,
and consider the morphism
$$KS_{\tau_0}^{\Ig}:   \Omega^1_{{Y}^{\mathrm{Ig}}/\F}\otimes_{\mathcal{O}_{{Y}^{\mathrm{Ig}}}}\mathcal{F}^{\mathrm{Ig}} 
\cong  \pi^*\Omega^1_{\overline{Y}/\F}\otimes_{\mathcal{O}_{{Y}^{\mathrm{Ig}}}} 
\mathcal{F}^{\mathrm{Ig}} \longrightarrow \pi^* (\CA_{2\be_{\theta_0},-\be_{\theta_0},\FF})
\otimes_{\mathcal{O}_{Y^{\mathrm{Ig}}}} \mathcal{F}^{\mathrm{Ig}}$$
induced by $KS_{\tau_0}$, where $\CF^\Ig$ is the sheaf of total fractions on ${Y}^\Ig$.  
Suppose now that $f \in M_{\bk,\bl}(U;\F)$, and write $h^{\bk} = \prod_{\theta \in \Sigma} h_\theta^{k_\theta}$ and
$g^{\bl} = \prod_{\theta \in \Sigma} g_\theta^{l_\theta}$ for any choice of trivializations $g_\theta$ of the line bundles $\ol{\CN}_\theta$
on $\ol{Y}$.  We then define the section
$$\Theta_{\tau_0}^\Ig(f) := h^{\bk} \pi^*(g^{\bl} H_{\theta_0}) KS_{\tau_0}^\Ig(d(h^{-\bk} \pi^*(g^{-\bl}f))),$$
where
\begin{equation} \label{eqn:thetawt} \bk' = \bk + n_{\theta_0} \be_{\sigma^{-1}\theta_0} + \be_{\theta_0}\quad\mbox{and}\quad \bl' = \bl + \be_{\theta_0}.\end{equation}
Furthermore, the section is independent of the choices of $g_\theta$ and invariant under the action of $(\CO_F/\gq)^\times$,
hence descends to a section of $\CA_{\bk',\bl',\F} \otimes_{\CO_{\ol{Y}}} \CF$, where $\CF$ is the sheaf of total fractions on $\ol{Y}$.
Denoting the section $\Theta_{\tau_0}(f)$, we have the following generalization of \cite[Thm.~8.2.2]{DS}:
\begin{theorem} \label{thm:theta} If $f \in M_{\bk,\bl}(U;\F)$, then $\Theta_{\tau_0}(f) \in M_{\bk',\bl'}(U;\F)$.  Moreover $\Theta_{\tau_0}(f)$ is divisible
by $H_{\theta_0}$ if and only if either $f$ is divisible by $H_{\theta_0}$ or $p|k_{\theta_0}$.
\end{theorem}
\begpf  We see exactly as in \cite{DS} that $\Theta_{\tau_0}(f)$ is regular on the ordinary locus of $\ol{Y}$, i.e., the complement
of the divisor $\cup_{\theta \in \Sigma} Z_\theta$ (where $Z_\theta$ was defined in \S\ref{sec:stratification}), so the theorem
reduces to proving that if $z$ is the generic point of an
irreducible component of $Z_{\theta_1}$ for some $\theta_1 \in \Sigma$, then 
\begin{itemize}
\item $\ord_z(\Theta_{\tau_0}(f)) \ge 0$,
\item if $\theta_1 = \theta_0$, then $\ord_z(\Theta_{\tau_0}(f)) > 0$ if and only if $p|k_{\theta_0}$ or $\ord_z(f) > 0$.
\end{itemize}
Let $R$ denote the discrete valuation ring $\CO_{\ol{Y},z}$, and for each $\tau \in \Sigma_{\gp,0}$ and $\theta \in \Sigma_\tau$,
let $y_\theta = y_{\tau,j}$ be a basis for the stalk $\ol{\CL}_{\theta,z} = \ol{\CL}_{\tau,j,z}$ over $R$ (for $j = 1,\ldots,e_\gp$).
For each $\theta \in \Sigma$, we may then write $$H_\theta y_\theta = r_\theta y_{\sigma^{-1}\theta}^{n_\theta}$$
for some $r_\theta = r_{\tau,j}  \in R$, and we let $r_\tau = \prod_{j=1}^{e_\gp} r_{\tau,j}$. By construction, we have
$T: = (\pi_*\CO_{{Y}^{\Ig}})_z = R[x_\tau]_{\tau \in \Sigma_0} / I$, where $I$ is the ideal generated by
$$x_{\phi^{-1}\circ\tau}^p  -  r_\tau x_\tau \mbox{\,\,for $\tau \in \Sigma_0$,}\quad\mbox{and}\quad
    \prod_{\tau \in \Sigma_{\gp,0}} x_\tau^{p-1} - \prod_{\tau \in \Sigma_{\gp,0}}r_{\tau}^{p-1} \mbox{\\,,for $\gp \in S_p$,}$$
where each $x_\tau$ is the dual basis of $y_{\tau,e_\gp}$.  We then have $h_{\tau,e_\gp} = x_\tau y_{\tau,e_\gp}$
(in $(\pi_*\pi^*\CL_{\tau,e_\gp})_z$), from which it follows that
$$h_{\tau,j} = r_{\tau,j+1} r_{\tau,j+2} \cdots r_{\tau,e_\gp} x_\tau y_{\tau,j}$$
for $j=1,\ldots,e_{\gp} - 1$, and hence that $h^{\bk} =   \varphi_{\bk} y^{\bk}$, where 
$y^{\bk} = \prod_{\theta\in \Sigma} y_\theta^{k_\theta}$ and
$$ \varphi_{\bk}  = \prod_{\gp \in S_p} \prod_{\tau \in \Sigma_{\gp,0}} \left( (r_\tau x_\tau)^{\sum_{\theta \in \Sigma_\tau} k_\theta}
   \prod_{j=1}^{e_\gp}  r_{\tau,j}^{-\sum_{j' = j}^{e_\gp} k_{\tau,j'}} \right)$$
(writing $k_{\tau_{\gp,i},j}$ for $k_{\theta_{\gp,i,j}}$ as usual,
and working over the field of fractions of $T$).

Writing $f = \varphi_f y^{\bk} g^{\bl}$, we see that
$$\Theta_{\tau_0}^\Ig(f) = KS_{\tau_0}^\Ig(  r_{\theta_0}  \varphi_{\bk} d (\varphi_f\varphi_{\bk}^{-1} ) )  y_{\theta_0}^{-1} y_{\sigma^{-1} \theta_0}^{n_{\theta_0}} y^{\bk} g^{\bl}.$$
Since $r_\tau  x_\tau = x_{\phi^{-1}\circ\tau}^p $, we have $d(r_\tau x_\tau) = 0$ and 
\begin{equation}
\label{eqn:altdeftheta}
\Theta_{\tau_0}(f) = KS_{\tau_0} \left(r_{\theta_0} d\varphi_f  + r_{\theta_0} \varphi_f \sum_{\theta \in \Sigma} k_\theta' \frac{dr_\theta}{r_\theta}\right) 
y_{\theta_0}^{-1} y_{\sigma^{-1} \theta_0}^{n_{\theta_0}} y^{\bk} g^{\bl},
\end{equation}
where $k'_\theta = k_{\tau,j} + k_{\tau,j+1} + \cdots k_{\tau,e_\gp}$ if $\tau = \tau_{\gp,i}$ and $\theta = \theta_{\gp,i,j}$.
We are therefore reduced to showing that
$\ord_z KS_{\tau_0} (dr_{\theta_1}) > 0$ if and only if $\theta_1 = \theta_0$.
However the proof of Proposition~\ref{prop:smooth} shows that if $y$ is a closed point of $Z_{\theta_1}$,
then $KS_{\tau_0}(dr_{\theta_1})$ vanishes at $y$ if and only if $\theta_1 = \theta_0$.
\epf

\begin{remark}  
The Kodaira--Spencer isomorphism is defined in \cite{DS} using the Gauss--Manin connection.
Much of the work in \cite[\S8.2]{DS} amounts to an explicit translation of this to the context of
deformation theory.  Here however we defined the morphism $KS_{\tau_0}$ more directly using
deformation theory, so the analogue of \cite[Lem.~8.2.3]{DS} was not needed here.
\end{remark}

\begin{remark}
It is straightforward to check directly that the right-hand side of (\ref{eqn:altdeftheta}) is independent of the choice
of local trivializations $y_\tau$ and $g_\tau$, and can therefore be used to {\em define} the partial $\Theta$-operator
without reference to the Igusa cover ${Y}^\Ig$.
\end{remark}

We call $\Theta_{\tau_0}$ the {\em partial $\Theta$-operator} (indexed by $\tau_0$).
It is immediate from its definition that the resulting map on $\FF$-algebras
$$\bigoplus_{\bk,\bl \in \ZZ^\Sigma} M_{\bk,\bl}(U;\F)  \longrightarrow \bigoplus_{\bk,\bl \in \ZZ^\Sigma} M_{\bk,\bl}(U;\F),$$
given by the direct sum over all weights of the operators $\Theta_{\tau_0}$, is an $\FF$-linear derivation, i.e. that
$$\Theta_{\tau_0}(f_1f_2) =  f_1 \Theta_{\tau_0}(f_2) + \Theta_{\tau_0}(f_1)f_2$$
for all $f_1,f_2$ in $\oplus M_{\bk,\bl}(U;\FF)$.  It is also clear that $\Theta_{\tau_0}(H_\theta) = 0$
for all $\theta \in \Sigma$, and hence that $\Theta_{\tau_0}$ commutes with multiplication by
partial Hasse invariants.

It is also straightforward to check that the operator $\Theta_{\tau_0}$ commutes with the Hecke action in the 
obvious sense, and hence induces a $\GL_2(\A_{F,\f}^{(p)})$-equivariant map
$$\Theta_\tau:  M_{\bk,\bl}(\F)  \longrightarrow M_{\bk',\bl'}(\F)$$
where $M_{\bk,\bl}(\F)$ (and $M_{\bk',\bl'}(\F)$, with their $\GL_2(\A_{F,\f}^{(p)})$ actions) are 
defined in (\ref{eqn:inflevel}) as direct limits over suitable open compact $U$.

Let us also make the effect of $\Theta_{\tau_0}$ on the weight $\bk$ more explicit.  Note that
if $\bk = \sum_\theta k_\theta\be_\theta$ and $\tau_0 = \tau_{\gp_0,i_0}$, then $\bk' = \sum_\theta k'_\theta\be_\theta$, where
\begin{itemize}
\item if $e_{\gp_0} = f_{\gp_0} = 1$, then 
$k'_\theta = \left\{\begin{array}{ll}  
k_\theta + p + 1, & \mbox{if $\theta = \theta_0 = \theta_{\gp_0,1,1}$,}\\
k_\theta,&{otherwise;}\end{array}\right.$
\item if $e_{\gp_0} = 1$ and $f_{\gp_0} > 1$, then
$k'_\theta = \left\{\begin{array}{ll}  
k_\theta + 1, & \mbox{if $\theta = \theta_0 = \theta_{\gp_0,i_0,1},$}\\
k_\theta + p, & \mbox{if $\theta = \sigma^{-1}\theta_0 = \theta_{\gp_0,i_0-1,1}$,}\\
k_\theta,&{otherwise;}\end{array}\right.$
\item if $e_{\gp_0} > 1$, then
$k'_\theta = \left\{\begin{array}{ll}  
k_\theta + 1, & \mbox{if $\theta = \theta_0 = \theta_{\gp_0,i_0,e_{\gp_0}}$ or $\theta = \sigma^{-1}\theta_0 = \theta_{\gp_0,i_0,e_{\gp_0}-1}$,}\\
k_\theta,&{otherwise.}\end{array}\right.$
\end{itemize}

\begin{remark} Considerations from the theory of Serre weights from the point of \cite{DS} suggest that the above weight shifts
are in a certain sense optimal.  One can also define cruder partial $\Theta$-operators by composing
with multiplication by (products of) partial Hasse invariants.  For example, the operator
$H_{\tau_0,1}H_{\tau_0,2}\cdots H_{\tau_0,e_{\gp_0} - 1}\Theta_{\tau_0}$ is the one constructed in \cite{DDW},
and for any $j=1,\ldots,e_{\gp_0}$, the operator
$$H_{\tau_0,j}H^2_{\tau_0,j+1}\cdots H^2_{\tau_0,e_{\gp_0}-1}H_{\tau_0,e_{\gp_0}}\Theta_{\tau_0}$$
shifts the weight $\bk$ by $\be_\theta + n_{\sigma^{-1}\theta}\be_{\sigma^{-1}\theta}$, where $\theta = \theta_{\gp_0,i_0,j}$.
\end{remark}

\section{Partial Frobenius operators}  \label{Frob}
\subsection{Partial Frobenius endomorphisms}  \label{sec:Frobenius}
In order to define partial Frobenius operators on Hilbert modular forms (in \S\ref{sec:V} below), we first need to define
partial Frobenius endomorphisms of Hilbert modular varieties over $\FF$.

Fix a prime $\gp$ dividing $p$, and a level $U$, assumed as usual to be sufficiently small and prime to $p$. 
We will draw on ideas from \cite[\S7.1]{DKS} to construct an isogeny on the universal abelian variety $s: A \to S$,
where $S = \widetilde{Y}_{U,\FF}$.

We begin by associating Raynaud data to the line bundles ${\CL}_{\gp,i,e_\gp}$ over $S$, which we write simply
as ${\CL}_i$ for $i \in \ZZ/f\ZZ =  \ZZ/f_\gp\ZZ $ (omitting the subscripts for the fixed $\gp$ and $j = e = e_\gp$).
We define $f_i:  {\CL}^{\otimes p}_i \to {\CL}_{i+1}$ to
be zero, and we define $v_i:{\CL}_{i+1}\to{\CL}^{\otimes p}_i$ to be the morphism induced by
the restriction of
$$\Ver^*_A:  \CH  = \CH^1_\dr(A/S)_\gp \to \CH^1_\dr((\phi_S^*A)/S)_\gp = \phi_S^*\CH^1_\dr(A/S)_\gp = \phi_S^*\CH$$
to $\CF_{i}^{(e)} = (s_*\Omega^1_{A/S})_{i}$ (abbreviating subscripts $\tau_{\gp,i}$ by $i$).  Note that since the image of 
$\CH_{i+1}$ under $\Ver_A^*$ is $\phi_S^*(\CF_{i}^{(e)})$, the inclusions $\CF_{i+1}^{(e-1)} \subset u \CH_{i+1}$
and $u\CF_i^{(e)} \subset \CF_i^{(e-1)}$ ensure that $\Ver_A^*(\CF_{i+1}^{(e-1)}) \subset \phi_S^*(\CF_i^{(e-1)})$, so
the morphism $v_i$ is well-defined.  We then let $H$ denote the
finite flat $(\CO_F/\gp)$-vector space scheme over $S$ associated to the Raynaud data 
$({\CL}_i, f_i, v_i)_{i \in\ZZ/f\ZZ}$.  

Recall that the Dieudonn\'e crystal of $\ker(\Frob_A)$ is canonically isomorphic to $\Phi^*(s_*\Omega^1_{A/S})$,
with $F = 0$ and $V$ induced by $\Phi^*(\Ver_A^*)$ (in the notation of \cite[\S4.4.3]{BBM}).  On the other hand the Dieudonn\'e crystal
of $H$ is identified with $\Phi^*(\oplus_i {\CL}_i)$ with $F = \Phi^*(\oplus_i f_i) = 0$ and $V = \Phi^*(\oplus_i v_i)$
(as a simple special case of \cite[Prop.~7.1.3]{DKS}).  Therefore the canonical projection 
$s_*\Omega^1_{A/S} \to \oplus_{i} {\CL}_i$ induces a surjective morphism
of  Dieudonn\'e crystals $\D(\ker(\Frob_A)) \to \D(H)$.  As the base $S$ is smooth over $\FF$,
the exact contravariant functor $\D$ is fully faithful on finite flat $p$-group schemes over $S$ (\cite[Thm.~4.1.1]{BM}),
so the surjection arises from a closed immersion $H \hookrightarrow \ker(\Frob_A)$, and we let
$$\alpha:  A \longrightarrow  A' :=  A/H$$
denote the resulting isogeny of abelian varieties over $S$.  Note that $A'$ naturally
inherits an $\CO_F$-action $\iota'$ from the action $\iota$ on $A$.

Let $\CI$ denote the image of the morphism $\alpha^*: \CH^1_\dr(A'/\widetilde{Y}_{U,\FF})_\gp \to \CH$.
By construction, we have the exact sequence
$$\begin{array}{ccccccc} \D(A'[p])_{S}  & \longrightarrow &  \D(A[p])_{S}  
& \longrightarrow & \D(H)_{S} & \longrightarrow & 0\\
\parallel\wr&&\parallel\wr && \parallel\wr &&\\
\CH^1_\dr(A'/S) &\stackrel{\alpha^*}{\longrightarrow} &\CH^1_\dr(A/S) & \stackrel{\Ver_A^*}{\longrightarrow}
& \bigoplus_{i} \phi_S^*(\CF_\tau^{(e)}/\CF_\tau^{(e-1)})  & \longrightarrow & 0,
\end{array}$$
showing that $\CI = \oplus_{i} \CI_i$, where $\CI_i$ is the preimage of $\phi_S^*(\CF_{i-1}^{(e)})$ under
$\Ver_{A,i}^*:\CH_i \to \phi_S^*(\CH_{i-1})$.  Note in particular that $u\CH_i \subset \CI_i$ for all $i$, so
that $H \subset A[\gp]$ and there is a unique isogeny $\beta: \gp \otimes_{\CO_F} A' \longrightarrow  A$
such that $\alpha\circ \beta$ is the canonical isogeny $\gp\otimes_{\CO_F} A'  \to A'$.

We now equip $A'$ with auxiliary data corresponding to an element of $\widetilde{Y}_{U,\FF}(S)$.

Since $\alpha$ induces isomorphisms $T^{(p)}(A_{\overline{s}}) \stackrel{\sim}{\longrightarrow} T^{(p)}(A'_{\overline{s}})$
for all geometric points $\overline{s}$ of $S$, we immediately have a level $U^p$ structure $\eta'$ on $A'$
inherited from $A$.

Next we claim that the quasi-polarization $\lambda$ on $A$ induces an 
isomorphism\footnote{Here $\gc$ depends on the connected component of $S$.}
$\gp\gc\gd \otimes_{\CO_F}  A' \to (A')^\vee$,
or equivalently $A' \to \gp^{-1}\gc^{-1}\gd^{-1} \otimes_{\CO_F} (A')^\vee$,
which amounts to the claim that $H$ corresponds to the kernel of
$$\gc^{-1}\gd^{-1} \otimes \beta^\vee:  \gc^{-1}\gd^{-1} \otimes_{\CO_F} A^\vee \longrightarrow \gp^{-1}\gc^{-1}\gd^{-1} \otimes_{\CO_F}( A')^\vee$$
under the isomorphism induced by $\lambda$.  Denoting this kernel by $I$, we have that $H$ and $I$ are finite flat group schemes
over $S$ of the same rank, so it suffices to prove that the composite
$$I \longrightarrow \gc^{-1}\gd^{-1} \otimes_{\CO_F} A^\vee[p] \stackrel{\sim}{\longrightarrow}  A[p] \longrightarrow A'[p]$$
is trivial.  Taking Dieudonn\'e modules, this in turn amounts to the vanishing of the composite
$$\D(A'[p])_{S}  \longrightarrow \D(A[p])_{S} \longrightarrow
  \D(\gc^{-1}\gd^{-1} \otimes_{\CO_F} A^\vee[p])_{S} \longrightarrow \D(I)_{S}.$$
We have already noted that the image of the first map has $\gp$-component
 $\oplus_{i} \CI_i$; on the other hand the kernel of the last map is the image of the map
 $$         \D(\gp^{-1}\gc^{-1}\gd^{-1} \otimes_{\CO_F} (A')^\vee[p])_{S} \longrightarrow
\D(\gc^{-1}\gd^{-1} \otimes_{\CO_F} A^\vee[p])_{S} $$
corresponding to the adjoint of 
$$\beta^*:  \CH^1_\dr(A/S)  \to \CH^1_\dr((\gp\otimes_{\CO_F} A')/S) \cong
     \gp^{-1} \otimes_{\CO_F} \CH^1_\dr(A'/S)$$
under the canonical isomorphisms
$$\begin{array}{rl} \D(\gc^{-1}\gd^{-1} \otimes_{\CO_F} A^\vee[p])_S
    \cong \CH^1_\dr((\gc^{-1}\gd^{-1} \otimes_{\CO_F} A^\vee / S) &\\
     \cong \Shom_{\CO_S}(\gd^{-1}\otimes_{\CO_F} \CH^1_\dr(A/ S), \CO_{S})&
      \cong \Shom_{\CO_F\otimes \CO_{S}}(\CH^1_\dr(A/ S), \CO_F \otimes \CO_{S})
      \end{array}$$
 and similarly
$$  \D(\gp^{-1}\gc^{-1}\gd^{-1} \otimes_{\CO_F} (A')^\vee[p])_{S }
    \cong  \Shom_{\CO_F\otimes \CO_{S}}(\gp^{-1} \otimes_{\CO_F}\CH^1_\dr(A'/S), \CO_F \otimes \CO_S)$$
 obtained from duality. 
 We are therefore reduced to proving that $\CI_i$ is orthogonal to the kernel of 
 $\beta_i^*$ for each $i \in  \ZZ/f\ZZ $ under the pairing $\langle\cdot,\cdot \rangle_i$
 defined by (\ref{eqn:pairing2}).  Note however that the kernel of $\beta_i^*$ is $u^{e-1} \CI_i$, as can be seen
 for example from the commutative diagram
 $$\xymatrix{ H^1_\crys(A_{\os}/W(\Fpbar))_i  \ar@{^{(}->}[r]\ar@{->>}[d] & (\gp^{-1} \otimes_{\CO_F} H^1_\crys(A'_{\os}/W(\Fpbar)))_i   \ar@{->>}[d] \\
  H^1_\dr(A_{\os}/\Fpbar))_i  \ar[r] & (\gp^{-1} \otimes_{\CO_F} H^1_\dr(A'_{\os}/\Fpbar))_i }$$
 of $W(\Fpbar)[u]$-modules for $\os \in S(\Fpbar)$.
Finally the orthogonality of $\CI_i$ and $u^{e-1}\CI_i$ is immediate from that of
$\CF_{i-1}^{(e-1)}$ and $u^{-1}\CF_{i-1}^{(e-1)}$ provided by Lemma~\ref{lem:ortho},
completing the proof of the claim.  We may then define the quasi-polarization on $A'$ by $\alpha^*(\lambda') = \delta \lambda$
for any totally positive generator $\delta = \delta_\gp$ of $\gp\CO_{F,(p)}$, so that $\lambda'$ induces an isomorphism
$\gc'\gd\otimes_{\CO_F} A' \stackrel{\sim}{\longrightarrow} (A')^\vee$ where $\gc' = \delta^{-1}\gp\gc$.

Finally we define a Pappas--Rapoport filtration on $\CF'_\tau := (s'_*\Omega^1_{A'/S})_\tau$ for all $\tau \in \Sigma_0$.

First note that if  $\tau \not\in \Sigma_{\gp,0}$, then $\alpha^*$ induces an isomorphism $\CF'_\tau \simeq \CF_\tau^{(e)}$,
and we define $(\CF_\tau')^{(j)}$ as the pre-image of $\CF_\tau^{(j)}$.

Suppose now that $\tau = \tau_{\gp,i}$.  Recall from the construction of $A' = A/H$ that 
$\Ver_A^*(\CF_{i}^{(e-1)}) \subset \phi_S^*(\CF_{i-1}^{(e-1)})$, so we have
$\CF_i^{(e-1)} \subset \CI_i$.  It follows that $(\alpha_i^*)^{-1}(\CF_i^{(e-1)})$
is a subbundle of $\CH'_i := \CH^1_\dr(A')_i$ of the same rank as $\CF'_i$, namely $e$.
Furthermore we have
$$\begin{array}{rl} \phi_S^*(\alpha^*(\CF_i')) =  (\phi_S^*(\alpha))^* (\CF_i') =  &(\phi_S^*(\alpha))^*(\Ver_{A'}^*(\CH_{i+1}'))\\
= &\Ver_A^*(\alpha^*(\CH_{i+1}')) = \Ver_A^*(\CI_{i+1}) \subset \phi_S^*(\CF_i^{(e+1)}),\end{array}$$
so in fact $\CF_i' \subset (\alpha_i^*)^{-1}(\CF_i^{(e-1)})$, and hence equality holds.
We thus obtain an exact sequence
$$0 \to \ker(\alpha_i^*) \longrightarrow  \CF'_i  \stackrel{\alpha_i^*}{\longrightarrow} \CF_i^{(e-1)} \to 0.$$
We may thus define a Pappas--Rapoport filtration on $\CF_i'$ by setting
$$(\CF_i')^{(j)} = (\alpha_i^*)^{-1}(\CF_i^{(j-1)})$$
for $j=1,\ldots,e$, so in particular $(\CF_i')^{(1)} = \ker(\alpha_i^*)$.

We now define $\widetilde{\Phi}_\gp:\widetilde{Y}_{U,\FF} = S \to \widetilde{Y}_{U,\FF}$ to be the endomorphism corresponding
to the data $(A',\iota',\lambda',\eta',(\CF')^\bullet)$.  Note that $\widetilde{\Phi}_\gp$ depends on the choice of $\delta$
in the definition of $\lambda'$; however it is straightforward to check that $\widetilde{\Phi}_\gp$ is compatible with the
$\CO_{F,(p),+}^\times$-action on $S$ and descends to an endomorphism $\Phi_\gp$ of $\ol{Y}_U$
which is independent of this choice.  We call $\widetilde{\Phi}_\gp$ (resp.~$\Phi_\gp$) the {\em partial Frobenius
endomorphism} (indexed by $\gp$) of $\widetilde{Y}_{U,\FF}$ (resp.~$\ol{Y}_U$); the terminology is justified
by the next proposition.

For the statement of the proposition, we also define the endomorphism $\widetilde{\Phi}$ of $S = \widetilde{Y}_{U,\FF}$
corresponding to the data $\phi_S^*(\uA) = (\phi_S^*A, \phi_S^*\iota,\phi_S^*\lambda,\phi_S^*\eta,(\phi_S^*\CF)^\bullet)$, 
where $(\phi_S^*\CF)^\bullet$ is the collection of filtrations on the vector bundles 
$$((\phi_S^*s)_*(\Omega^1_{(\phi_S^*A)/S}))_\tau  = (\phi_S^*(s_*\Omega^1_{A/S}))_\tau = \phi_S^*((s_*\Omega^1_{A/S})_{\phi^{-1}\circ\tau})$$
given by $(\phi_S^*\CF)_\tau^{(j)} = \phi_S^*(\CF_{\phi^{-1}\circ\tau}^{(j)})$.
Note that $\widetilde{\Phi}$ is not the absolute Frobenius $\phi_S$  on $S$ (unless $\FF = \FF_p$),
but we may write $\phi_S = \widetilde{\Phi}\circ\widetilde{\epsilon}$ where $\widetilde{\epsilon}$ is the isomorphism
defined by the commutative diagram
$$\xymatrix{S \ar[drr]^{\widetilde{\epsilon}}  \ar[dr]^{\varepsilon} \ar[ddr] &&\\
& \phi^*S \ar[r]^{\sim} \ar[d] & S \ar[d] \\
& \Spec \FF \ar[r]^{\sim}_{\phi} & \Spec \FF},$$
where the square is Cartesian and $\varepsilon$ is the 
inverse of the isomorphism associated to $\phi^*A = A \times_{\FF,\phi} \FF$ with the evident auxiliary data. 
We thus have an isomorphism $\widetilde{\epsilon}^*A \cong A$ compatible with $\iota$, $\lambda$ and $\eta$, and
 inducing $\widetilde{\epsilon}^*\CF_\tau^{(j)} \cong \CF_{\phi\circ\tau}^{(j)}$ for all $\tau$ and $j$.  (Note also that $\widetilde{\Phi}$
 may be viewed as the base-change of the absolute Frobenius on
the descent of $S$ to $\FF_p$ defined by the diagram.)  

The endomorphism $\widetilde{\Phi}$ is compatible  with the $\CO_{F,(p),+}$-action on $S = \widetilde{Y}_{U,\FF}$, and we let $\Phi$ denote
the resulting endomorphism of $\ol{Y}_U$.  Similarly $\widetilde{\epsilon}$ descends to a $\phi$-linear 
automorphism $\epsilon$ of $\ol{Y}_U$ such that the absolute Frobenius on $\ol{Y}_U$ is $\Phi\circ\epsilon$.

\begin{proposition}  The morphisms $\Phi_\gp$ are finite flat of (constant) degree $\Nm_{F/\QQ}(\gp)$, commute with each other,
and satisfy the formula
$$\prod_{\gp \in S_p} \Phi_\gp^{e_\gp} = \Phi.$$
\end{proposition}
\begpf  We first prove the commutativity and analogous formula for the maps $\widetilde{\Phi}_\gp$ on
$S = \widetilde{Y}_{U,\FF}$, from which the corresponding assertions for $\Phi_\gp$ follow.  To that end it
suffices to consider the maps on geometric closed points $\overline{s}\in S(\Fpbar)$, which we will do in order
to facilitate computations on Dieudonn\'e modules.

Let $\uA_0$ denote the data corresponding to $\overline{s} \in S(\Fpbar)$ and $\uA_0' = \uA_{0,\gp}'$ denote the data corresponding
to $\widetilde{\Phi}_\gp(\overline{s})$.   Let $D = H^1_\crys(A_0/W(\Fpbar))$ and $D' = H^1_\crys(A_0'/W(\Fpbar))$, so we may decompose
the $\CO_F\otimes W(\Fpbar)$-modules $D = \oplus_{\tau\in \Sigma} D_\tau$ and $D' = \oplus_{\tau \in \Sigma'} D_\tau'$
where $D_\tau$ and $D_\tau'$ are free $W(\Fpbar)[u]/(f^\tau(u))$-modules of rank two.  Furthermore
the canonical isogeny $\alpha: A_0 \to A_0'$ induces an injective $W(\Fpbar)[u]/(f^\tau(u))$-linear map
$\alpha_\tau^*: D_\tau' \to D_\tau$ for each $\tau$, compatible in the obvious sense with the maps
$$\Frob_{A_0,\tau}^*: \phi^*(D_{\phi^{-1}\circ\tau})  \to D_\tau\quad\mbox{and}\quad\Frob_{A_0',\tau}^*:\phi^* (D_{\phi^{-1}\circ\tau}') \to D'_\tau,$$
as well as $\Ver_{A_0,\tau}^* = p(\Frob_{A_0,\tau}^*)^{-1}$ and $\Ver_{A_0',\tau}^* = p(\Frob_{A_0',\tau}^*)^{-1}$.
Letting $F_\tau^{(j)} \subset D_\tau$ denote the pre-image of $\CF_\tau^{(j)}$ under
the canonical surjection
$$D_\tau \longrightarrow (D/pD)_\tau \cong H^1_\dr(A_0/\Fpbar)_\tau,$$
we have  $\alpha_\tau^*(D_\tau') = D_\tau$ if $\tau \not\in \Sigma_{\gp,0}$ and
$$\alpha_\tau^*(D_\tau') = (\Ver_{A_0,\tau}^*)^{-1}(\phi^*(F_{\phi^{-1}\circ\tau}^{(e_\gp-1)}))$$
if $\tau \in \Sigma_{\gp,0}$ (by the construction of $\Phi_\gp$).  Furthermore writing $F'^{(j)}_\tau$ for the
submodules of $D_\tau'$ similarly defined by the Pappas--Rapoport filtration on $H^0(A_0',\Omega^1_{A_0'/\Fpbar})_\tau$,
we have 
\begin{itemize}
\item $\alpha_\tau^*(F'^{(j)}_\tau) = F^{(j)}_\tau$ if $\tau \not\in \Sigma_{\gp,0}$, 
\item $\alpha_\tau^*(F'^{(j)}_\tau) = F^{(j-1)}_\tau$ for $j=2,\ldots,e_\gp$ if $\tau \in \Sigma_{\gp,0}$,
\item and $\alpha_\tau^*(F'^{(1)}_\tau) = pD_\tau = \Frob_{A_0,\tau}^*(\phi^*(F^{(e_\gp)}_{\phi^{-1}\circ\tau}))$ if $\tau\in \Sigma_{\gp,0}$.
\end{itemize}

Thus if $\gp_1$ and $\gp_2$ are distinct elements of $S_p$, then $\widetilde{\Phi}_{\gp_1}(\widetilde{\Phi}_{\gp_2}(\overline{s}))$ corresponds to the data
$\uA_0''$ for which we have an isogeny $\alpha':A_0 \to A_0''$ such that if $\tau \not\in \Sigma_{\gp_1,0} \cup \Sigma_{\gp_2,0}$
then $(\alpha')_\tau^*(D_\tau'') = D_\tau$ and $(\alpha')_\tau^*(F''^{(j)}_\tau) = F^{(j)}_\tau$
(with the obvious notation), but if $\tau \in \Sigma_{\gp_i,0}$ for $i=1$ or $2$, then
$$(\alpha')_\tau^*(D_\tau'') = (\Ver_{A_0,\tau}^*)^{-1}(\phi^*(F_{\phi^{-1}\circ\tau}^{(e_{\gp_i}-1)})),$$
$(\alpha')_\tau^*(F''^{(1)}_\tau) = pD_\tau$ and $(\alpha')_\tau^*(F''^{(j)}_\tau) = F^{(j-1)}_\tau$ for $j=2,\ldots,e_{\gp_i}$.
Furthermore we have $\eta'' = \alpha'\circ\eta$ and $(\alpha')^*(\lambda'') = \delta_{\gp_2}\delta_{\gp_1}\lambda$,
from which it follows that the isomorphism class of the data $\uA_0''$ also corresponds to $\widetilde{\Phi}_{\gp_2}(\widetilde{\Phi}_{\gp_1}(\overline{s}))$.

Now consider the data $\uA_0^{(r)}$ associated to $\widetilde{\Phi}_\gp^r(\overline{s})$ for $r=1,\ldots,e_\gp$, and write $D_r = \oplus D_{r,\tau}$
for $H^1_\crys(A_0^{(r)}/W(\Fpbar))$, $F_{r,\tau}$ for the submodule of $D_{r,\tau}$ determined as above by the Pappas--Rapoport
filtration, and $\alpha_r$ for the composite isogeny $A_0 \to A_0^{(1)} \to \cdots \to A_0^{(r)}$.  By induction on $r$, we find that
if $\tau \not\in \Sigma_{\gp,0}$, then $\alpha_{r,\tau}^*(D_{r,\tau}) = D_\tau$ and $\alpha_{r,\tau}^*(F_{r,\tau}^{(j)}) = F^{(j)}_\tau$,
but if $\tau\in \Sigma_{\gp,0}$, then 
\begin{itemize}
\item $\alpha_{r,\tau}^*(D_{r,\tau}) = (\Ver_{A_0,\tau}^*)^{-1}(\phi^*(F_{\phi^{-1}\circ\tau}^{(e_{\gp}-r)}))$,
\item $\alpha_{r,\tau}^*(F^{(j)}_{r,\tau}) = F^{(j-r)}_\tau$ for $j=r+1,\ldots,e_\gp$,
\item $\alpha_{r,\tau}^*(F^{(j)}_{r,\tau}) = \Frob_{A_0,\tau}^*(\phi^*(F^{(e_\gp+j-r)}_{\phi^{-1}\circ\tau}))$ for $j=1,\ldots,r$.
\end{itemize}
In particular taking $r=e_\gp$ gives 
$$\alpha_{e_\gp,\tau}^*(D_{e_\gp,\tau}) = (\Ver_{A_0,\tau}^*)^{-1}(p\phi^*(D_{\phi^{-1}\circ\tau})) = \Frob_{A_0,\tau}^*(\phi^*(D_{\phi^{-1}\circ\tau}))$$
and $\alpha_{e_\gp,\tau}^*(F^{(j)}_{e_\gp,\tau}) = \Frob_{A_0,\tau}^*(\phi^*(F^{(j)}_{\phi^{-1}\circ\tau}))$ for $j=1,\ldots,e_\gp$ and
$\tau \in \Sigma_{\gp,0}$.  It then follows that $\prod_{\gp \in S_p}\widetilde{\Phi}_\gp^{e_\gp}(\overline{s})$ corresponds to $\uA_0''$ with $\alpha':A_0\to A_0''$
satisfying $(\alpha')_\tau^*(D_\tau'') = \Frob_{A_0,\tau}^*(\phi^*(D_{\phi^{-1}\circ\tau}))$
and $(\alpha')_\tau^*(F''^{(j)}_{\tau}) = \Frob_{A_0,\tau}^*(\phi^*(F^{(j)}_{\phi^{-1}\circ\tau}))$ for all $\tau$ and $j$.
Furthermore we have $\eta'' = \alpha'\circ\eta$ and $(\alpha')^*(\lambda'') = \prod_{\gp} \delta_{\gp}^{e_\gp}\lambda$,
from which it follows that $\uA_0''$ is isomorphic to $(\phi^*A_0,\phi^*\iota,\nu\phi^*\lambda,(\phi^*\CF)^\bullet)$ with
$\nu = p^{-1}\prod_{\gp} \delta_{\gp}^{e_\gp} \in \CO_{F,(p),+}^\times$.  This proves that
$$\prod_{\gp \in S_p} \widetilde{\Phi}_\gp^{e_\gp}= \nu \cdot \widetilde{\Phi},$$
which in turn implies the desired formula.

Since $\Phi$ is finite (and $\ol{Y}_U$ is separated), it follows that $\Phi_\gp$ is finite, and therefore
also flat since $\ol{Y}_U$ regular.   Note furthermore that $\Phi_\gp$ is therefore bijective on closed points 
and induces isomorphisms on their residue fields, so the degree of $\Phi_\gp$ in a neighborhood of any
closed point $x$ of  $\ol{Y}_U$ is that of the extension of completed regular local rings 
$\Phi_{\gp,x}^*:  \CO_{\ol{Y}_U,y}^\wedge \to \CO_{\ol{Y}_U,x}^\wedge$, where $y = \Phi_\gp(x)$.
Since $\Phi$ factors through $\Phi_\gp$, so does the absolute Frobenius on $\ol{Y}_U$, and hence
the absolute Frobenius on $\CO_{\ol{Y}_U,x}^\wedge$ factors through $\Phi_{\gp,x}^*$.
Therefore it follows from \cite[Cor.~2]{KN} that $\deg(\Phi_{\gp,x}^*) = p^n$ where $n = n_\gp$ is the
dimension of the kernel of the induced map on tangent spaces $T_x(\ol{Y}_U) \to T_y(\ol{Y}_U)$.
Furthermore since $\prod_{\gp \in S_p} \Phi_\gp^{e_\gp} = \Phi$ has degree 
$p^{[F:\QQ]} = \prod_{\gp \in S_p}p^{ e_\gp f_\gp}$, it suffices to prove that
$n_\gp \ge  f_\gp = Nm_{F/\QQ}(\gp)$ for each $\gp$.   Note also that we may replace $\ol{Y}_U$
by $S=\widetilde{Y}_{U,\FF}$, $x$ by any point in its pre-image in $S$ and $\Phi_\gp$ by
$\widetilde{\Phi}_{\gp}$.

Suppose then that $x$ corresponds to the data $(A_0,\iota_0,\lambda_0,\eta_0,\CF_0^\bullet)$
over the residue field $k$, and its image $y =\widetilde{\Phi}_\gp(x)$ corresponds to the data $(A'_0,\iota'_0,\lambda'_0,\eta'_0,\CF'^{\bullet}_0)$.
Recall that the Kodaira--Spencer filtration on the fibre of $\Omega^1_{S/\Fpbar}$ at $x$
is dual to one on $T_x(S)$ which was described using Grothendieck--Messing deformation theory
(see the discussion following the proof of Theorem~\ref{thm:KSfil}).
In particular, we have a decomposition $T_x(S) = \oplus_{\tau \in \Sigma_0} T_x(S)_\tau$ and a decreasing
filtration of length $e_{\gp'}$ on $T_x(S)_\tau$ for each $\tau \in \Sigma_{\gp',0}$ (where $\gp' \in S_p$)
such that
$$\bigoplus_{\tau \in \Sigma_0}  \Fil^{j_\tau} T_x(S)_\tau$$
corresponds to the set of lifts of $\uA_0$ to $\uA_1 = (A_1,\iota_1,\lambda_1,\eta_1,\CF_1^\bullet) \in S(k[\epsilon])$ such that
$\CF_{1,\tau}^{(j)}$ is the image of $\CF_{0,\tau}^{(j)} \otimes_k k[\epsilon]$ for all $\tau$ and $j \le j_\tau$ under the
canonical isomorphism
$$H^1_\dr(A_1/k[\epsilon]) \cong H^1_\cris(A_0/k[\epsilon]) \cong H^1_\dr(A_0/k) \otimes_k k[\epsilon].$$
We claim that the $f_\gp$-dimensional subspace
$\oplus_{\tau \in \Sigma_{\gp,0}}  \Fil^{e_\gp-1} T_x(S)_\tau$ is contained in the kernel of 
$T_x(S) \to T_y(S)$.  Indeed if $\uA_1$ is a lift corresponding to an element of this subspace
and $\uA_1'$ is its image in $T_y(S)$ and $\alpha_i:A_i \to A_i'$ are the specializations of the
universal isogeny $\alpha:A \to A'$, then the commutativity of the diagram
$$\xymatrix{ H^1_\dr(A_1'/k[\epsilon]) \ar[r]^{\alpha_1^*} \ar[d]^{\wr} &  H^1_\dr(A_1/k[\epsilon]) \ar[d]^{\wr} \\
 H^1_\dr(A'_0/k) \otimes_k k[\epsilon] \ar[r]^{\alpha_0^*\otimes 1}   & H^1_\dr(A_0/k) \otimes_k k[\epsilon]}$$
 and the definition of $\widetilde{\Phi}_\gp$ imply that 
$\CF'^{(j)}_{1,\tau}$ corresponds to $\CF'^{(j)}_{0,\tau} \otimes_k k[\epsilon]$ for all $\tau$ and $j$.
(Note in particular that $\CF_{1,\tau}^{(1)} = \ker(\alpha_{1,\tau}^*)$ for all $\tau \in \Sigma_{\gp,0}$,
and that $H^0(A_1',\Omega_{A_1'/k[\epsilon]})_\tau$ corresponds to $H^0(A_0,\Omega_{A_0'/k})_\tau \otimes_k k[\epsilon]$
for all $\tau \in \Sigma_0$.)  It follows that $\uA_1'$ is the trivial deformation of $\uA_0'$, so the kernel of
$T_x(S) \to T_y(S)$ has dimension $n \ge f_\gp$ as required.
\epf

\begin{remark}  \label{rmk:Phidiffs}  The final part of the proof of the proposition shows that the kernel of the pull-back map
$\Phi_\gp^*\Omega^1_{\ol{Y}_U/\FF} \to \Omega^1_{\ol{Y}_U/\FF}$ is precisely 
$$\Phi_\gp^*\left(\oplus_{\tau \in \Sigma_{\gp,0}}  \Fil^1(\Omega^1_{\ol{Y}_U/\FF})_\tau\right).$$
Furthermore a similar argument shows that the map preserves the Kodaira--Spencer decomposition
and filtration, in the obvious sense, and induces isomorphisms
$$\begin{array}{cc}
 \Phi_\gp^*\left( \Fil^j(\Omega^1_{\ol{Y}_U/\FF})_\tau/\Fil^1(\Omega^!_{\ol{Y}_U/\FF})_\tau \right)
\stackrel{\sim}{\longrightarrow}\Fil^{j-1}(\Omega^1_{\ol{Y}_U/\FF})_\tau, &
\mbox{if $\tau \in \Sigma_{\gp,0}, j = 1,\ldots, e_\gp$,}\\
\Phi_\gp^*\left( \Fil^j(\Omega^1_{\ol{Y}_U/\FF})_\tau \right)
\stackrel{\sim}{\longrightarrow}\Fil^j(\Omega^1_{\ol{Y}_U/\FF})_\tau , &
\mbox{if $\tau \not\in \Sigma_{\gp,0}$.}\end{array}$$
\end{remark}

\subsection{Construction of $V_\gp$}  \label{sec:V}
In this section we generalize the construction\footnote{The operators defined here differ slightly
from the ones defined in \cite{DS} in the unramified case.  The construction there is tailored to
be compatible with the classical case and to be simply interpreted on $q$-expansions at cusps at $\infty$.
Doing this in the general ramified case would introduce complications that make it seem not worthwhile.}
of \cite[\S9.8]{DS} to define partial Frobenius operators, similar to the $V_p$-operator on classical
modular forms.  

We maintain the notation of \S\ref{sec:Frobenius}, so that $\widetilde{\Phi}_\gp$ is an endomorphism of
$S = \widetilde{Y}_{U,\FF}$ corresponding to the data $(A',\iota',\lambda',\eta',(\CF')^\bullet)$, where
$A' = A/H$ for a certain finite flat subgroup scheme $H \subset A[\gp]$, and $\alpha$ is the projection
$A \to A'$.

It is immediate from the definition of $\CF'^{(j)}_\tau$ that $\alpha_\tau^*$ induces an isomorphism
$\CL_{\tau,j}'  \stackrel{\sim}{\longrightarrow} \CL_{\tau,j}$ for all $j$ if $\tau\not\in \Sigma_{\gp,0}$, as well as
$$\CL_{\tau,j}' =  \CF'^{(j)}_\tau/\CF'^{(j-1)}_\tau \stackrel{\sim}{\longrightarrow} 
 \CF^{(j-1)}_\tau/\CF^{(j-2)}_\tau   = \CL_{\tau,j-1}$$
 for $j=2,\ldots,e_\gp$ if $\tau\in \Sigma_{\gp,0}$.  To describe $\CL_{\tau,1}'$ for $\tau\in \Sigma_{\gp,0}$,
 note that since $H \subset \ker(\Frob_A)$, there is an isogeny $\gamma:A' \to \phi_S^*A$ such that
 $\gamma\circ\alpha = \Frob_A$.  We thus obtain a surjection 
 $$\gamma_\tau^*:  \phi_S^*(\CF_{\phi^{-1}\circ\tau}^{(e_\gp)})  = \ker(\Frob_A)_\tau^*
     \longrightarrow \ker(\alpha_\tau^*) = \CL'_{\tau,1}.$$
 Furthermore $\phi_S^*(\CF_{\phi^{-1}\circ\tau}^{(e_\gp-1)}) = \ker(\gamma_\tau^*)$, as can be seen on
 closed points, so we obtain a surjection, hence isomorphism, of line bundles
 $$\CL_{\phi^{-1}\circ\tau,e_\gp}^{\otimes p} \cong \phi_S^*(\CL_{\phi^{-1}\circ\tau,e_\gp}) = 
 \phi_S^*(\CF_{\phi^{-1}\circ\tau}^{(e_\gp)} )/ \phi_S^*(\CF_{\phi^{-1}\circ\tau}^{(e_\gp-1)})
 \stackrel{\gamma_\tau^*}{\longrightarrow} \CL'_{\tau,1}.$$
 By construction, we have $\widetilde{\Phi}_\gp^*\CL_{\tau,j}  = \CL'_{\tau,j}$, so we have now defined
 isomorphisms
 $$\widetilde{\Phi}_\gp^*\CL_{\theta}  \cong \left\{  \begin{array}{ll}  \CL_{\sigma^{-1}\theta}^{\otimes n_\theta},&\mbox{if $\theta\in\Sigma_{\gp}$;}\\
 \CL_{\theta},&\mbox{if $\theta \not\in\Sigma_{\gp}$.}\end{array}\right.$$
 Similarly we find that if $\tau\in \Sigma_{\gp,0}$, then
$\gamma^*$ induces $\phi_S^*(\CP_{\phi^{-1}\circ \tau,e_\gp}) \cong \CP'_{\tau,1}$
and $\alpha^*$ induces $\CP'_{\tau,j} \cong \CP_{\tau,j-1}$ for $j=2,\ldots,e_\gp$, so 
that $\widetilde{\Phi}_\gp^* \CN_\theta \cong \CN_{\sigma^{-1}\theta}^{\otimes n_\theta}$ if $\theta\in \Sigma_\gp$,
while $\widetilde{\Phi}_\gp^*\CN_\theta \cong \CN_\theta$ if $\theta\not\in \Sigma_\gp$.
Taking tensor products, we thus obtain isomorphisms
$$\Phi_\gp^*\widetilde{\CA}_{\bk,\bl,\FF} \cong \widetilde{\CA}_{\bk'',\bl'',\FF}$$
for all $\bk,\bl\in \ZZ^\Sigma$, where 
\begin{itemize}
\item $k''_\theta = k_\theta$ and $l_\theta'' = l_\theta$ if $\theta\not\in\Sigma_\gp$, and
\item $k''_\theta = n_{\sigma\theta} k_{\sigma\theta}$ and $l''_\theta = n_{\sigma\theta} l_{\sigma\theta}$ if $\theta\in \Sigma_\gp$.
\end{itemize}
Furthermore it is straightforward to check that the isomorphisms are compatible with the descent data 
relative to $S=\widetilde{Y}_{U,\FF} \to Y_{U,\FF}$, so we obtain isomorphisms
\begin{equation}  \label{eqn:Phibundles} \Phi_\gp^*\CA_{\bk,\bl,\FF} \cong \CA_{\bk'',\bl'',\FF} \end{equation}
for all $\bk,\bl \in \ZZ^\Sigma$ and sufficiently small $U$.   Note that
\begin{equation} \label{eqn:Phiweight} \bk'' = \bk + \sum_{\theta\in \Sigma_\gp} 
k_\theta \bh_\theta\quad\mbox{and}\quad\bl'' = \bl + \sum_{\theta\in \Sigma_\gp} l_\theta \bh_\theta.\end{equation}

\begin{remark}  One can check that the resulting isomorphisms
$\Phi_\gp^*\CA_{2\be_\theta,-\be_\theta,\FF} \cong \CA_{2\be_\theta,-\be_\theta,\FF}$
(for $\theta \not\in \Sigma_\gp$), and 
$$\Phi_\gp^*\CA_{2\be_\theta,-\be_\tau,\FF} \cong \CA_{2\be_{\sigma^{-1}\theta},-\be_{\sigma^{-1}\theta},\FF}$$
(for $\theta = \theta_{\gp,i,j}$, $j=2,\ldots,e_\gp$) are compatible via the Kodaira--Spencer isomorphisms of 
Theorem~\ref{thm:KSfil} with the corresponding isomorphisms 
$$\Phi_\gp^*(\gr^j(\Omega^1_{\ol{Y}_U/\FF})_\tau) \cong \gr^{j'}(\Omega^1_{\ol{Y}_U/\FF})_\tau$$
given by Remark~\ref{rmk:Phidiffs}, where $j' = j$ if $\tau \not\in \Sigma_{\gp,0}$ and $j' = j-1$ if $\tau \in  \Sigma_{\gp,0}$.
\end{remark}

We are now ready to define the {\em partial Frobenius operator} (indexed by $\gp$)
$$V_\gp:  M_{\bk,\bl}(U;\FF)  \longrightarrow M_{\bk'',\bl''}(U,\FF)$$
as the composite
$$H^0(Y_{U,\FF},\CA_{\bk,\bl,\FF}) \stackrel{\Phi_\gp^*}{\longrightarrow} H^0(Y_{U,\FF},\Phi_\gp^*\CA_{\bk,\bl,\FF})
   \stackrel{\sim}{\longrightarrow}H^0(Y_{U,\FF},\CA_{\bk'',\bl'',\FF}),$$
where the second map is the isomorphism (\ref{eqn:Phibundles}).
It is immediate from the definition that $V_\gp$ is injective, and that taking the direct sum 
over all weights yields an $\FF$-algebra homomorphism 
$$\bigoplus_{\bk,\bl \in \ZZ^\Sigma} M_{\bk,\bl}(U;\FF) \longrightarrow \bigoplus_{\bk,\bl \in \ZZ^\Sigma} M_{\bk,\bl}(U;\FF)$$
for all sufficiently small $U$ containing $\GL_2(\CO_{F,p})$.
It is also straightforward to check that $V_\gp$ is compatible with the Hecke action in the usual sense, and hence defines
a $\GL_2(\AA_{F,\f}^{(p)})$-equivariant map 
$$M_{\bk,\bl}(\FF)  \longrightarrow M_{\bk'',\bl''}(\FF),$$
where the spaces are defined in (\ref{eqn:inflevel}) as direct limits over sufficiently small $U$ containing $\GL_2(\CO_{F,p})$.

It will also be convenient at times to consider instead the operator
\begin{equation} \label{eqn:Vp0} V^0_\gp:  M_{\bk,\bl}(U;\FF) \longrightarrow M_{\bk'',\bl}(U,\FF)\end{equation}
defined by $V^0_\gp(f) = V_\gp(f) \prod_{\theta\in \Sigma_\gp} G_\theta^{-l_\theta}$,
where $G_\theta$ is the trivialization of $\CA_{\bf{0},\bh_\theta,\FF}$ defined
at the end of \S\ref{sec:Hasse}.  Thus $V_\gp^0$ is also Hecke-equivariant,
but depends on the choice of uniformizer $\varpi_\gp$.

We also record the relation between the partial Frobenius operators and the $p$-power map.
First note that the identification $\widetilde{\Phi}^*(\CF_\tau^{(j)}) = \phi_S^*(\CF_{\phi^{-1}\circ\tau}^{(j)})$ arising from
the definition of $\widetilde{\Phi}$ yields isomorphisms 
$$\widetilde{\Phi}^*\CL_{\tau,j}  \cong \phi_S^*(\CL_{\phi^{-1}\circ\tau,j} ) \cong \CL_{\phi^{-1}\circ\tau,j}^{\otimes p}$$
for all $\tau$ and $j$.   We similarly have $\widetilde{\Phi}^*\CN_{\tau,j}    \cong \CN_{\phi^{-1}\circ\tau,j}^{\otimes p}$,
and taking tensor products and descending to $\ol{Y}$ yields isomorphisms
$\Phi^*\CA_{\bk,\bl,\FF} \cong \CA_{p\bk^\phi,p\bl^\phi,\FF}$ for all $\bk,\bl$, where $k^\phi_{\theta_{\gp,i,j}} = k_{\theta_{\gp,i+1,j}}$,
and hence an operator $V_p: M_{\bk,\bl}(U;\FF) \to M_{p\bk^\phi,p\bl^\phi}(U;\FF)$.
Similarly the isomorphisms $\widetilde{\epsilon}^*\CF_{\phi^{-1}\circ\tau}^{(j)} \cong \CF_{\tau}^{(j)}$ yield 
$\widetilde{\epsilon}^*\CL_{\phi^{-1}\circ\tau,j}  \cong \CL_{\tau,j}$ and $\widetilde{\epsilon}^*\CN_{\phi^{-1}\circ\tau,j}  \cong \CN_{\tau,j}$
whose tensor products descend to isomorphisms 
$\epsilon^*\CA_{\bk^\phi,\bl^\phi,\FF} \stackrel{\sim}{\longrightarrow}\CA_{\bk,\bl,\FF}$, yielding a $\phi$-linear
isomorphism $M_{\bk^\phi,\bl^\phi}(U;\FF) \to M_{\bk,\bl}(U;\FF)$ which we denote $\epsilon_{\bk,\bl}$.
Furthermore the above isomorphisms of line bundles on $S$ are compatible in the sense that the resulting diagram
 $$\xymatrix{\widetilde{\epsilon}^*\widetilde{\Phi}^*\CL_{\tau,j}   \ar[r]^{\sim}\ar[d]_{\wr} &
  \widetilde{\epsilon}^*\CL_{\phi^{-1}\circ\tau,j}^{\otimes p} \ar[d]^{\wr} \\
  \phi_S^*\CL_{\tau,j} \ar[r]^{\sim} & \CL_{\tau,j}^{\otimes p}  }$$
commutes, as does its analogue for the $\CN_\theta$, from which it follows that the composite
$$M_{\bk,\bl}(U;\FF) \stackrel{V_p}{\longrightarrow} M_{p\bk^\phi,p\bl^\phi}(U;\FF)
\stackrel{\epsilon_{p\bk,p\bl}}{\longrightarrow} M_{p\bk,p\bl}(U;\FF)$$
is the $p$-power map.

Returning to the partial Frobenius operators, 
the isomorphisms between $\widetilde{\Phi}_\gp^* \CL_\theta$ and $\CL_{\sigma^{-1}\theta}^{\otimes n_\theta}$ (resp.~$\CL_\theta$)
for $\theta \in \Sigma_\gp$ (resp.~$\theta\not\in\Sigma_\gp$) for different $\gp \in S_p$ are compatible with each other
in the obvious sense, and taken together with the formula $\prod_{\gp \in S_p} \widetilde{\Phi}_\gp^{e_\gp}= \nu \cdot \widetilde{\Phi}$
and the canonical isomorphism $\nu^*\CL_\theta \cong \CL_\theta$ yield the isomorphisms
$\widetilde{\Phi}^*\CL_{\tau,j} \cong \CL_{\phi^{-1}\circ\tau,j}^{\otimes p}$ defined above.
A similar assertion holds for the line bundles $\CN_\theta$, and
it follows that the operators $V_\gp$ for $\gp \in S_p$ commute with each other
and that $\prod_\gp V_\gp^{e_\gp} = V_p$, so that
\begin{equation}  \label{eqn:ppower} \left(\epsilon_{p\bk,p\bl} \prod_{\gp\in S_p}  V_\gp^{e_\gp} \right)   (f)  = f^p. \end{equation}

\section{Compactifications and $q$-expansions}  \label{q1}
 \subsection{Toroidal compactifications}  \label{sec:tor}
We next recall how $q$-expansions of Hilbert modular forms are obtained using
compactifications of Hilbert modular varieties.  In this section we review properties of the toroidal compactification constructed
by Rapoport~\cite{rap} (see also~\cite{chai} and~\cite{dim}).  We will consider toroidal compactifications only in the case $U = U(N)$,
but we first describe the set of cusps adelically for any $U$ of level prime to $p$.

For an arbitrary open compact subgroup $U$ of $\GL_2(\A_{F,\f})$ containing $\GL_2(\CO_{F,p})$,
we define the set of {\em cusps} of $Y_U$ to be
$$Y_U^\infty = B(\CO_{F,(p)})_+ \backslash \GL_2(\A_{F,\f}^{(p)}) /  U^p =  B(F)_+\backslash \GL_2(\A_{F,\f}) / U,$$
where $B$ denotes the subgroup of $\GL_2$ consisting of upper-triangular matrices.  Similarly we
define the set of cusps of $\widetilde{Y}_U$ to be 
$$\widetilde{Y}_U^\infty = B_1(\CO_{F,(p)})_+ \backslash \GL_2(\A_{F,\f}^{(p)}) /  U^p$$
Note that the natural surjection  $\widetilde{Y}^\infty_U\to Y^\infty_U$ identifies $Y^\infty_U$
with the quotient of $\widetilde{Y}^\infty_U$ by the left action of $B(\CO_{F,(p)})_+/B_1(\CO_{F,(p)}) \cong \CO_{F,(p),+}^\times$.
Furthermore the subgroup $(\CO_F^\times \cap U)^2$ acts trivially on $\widetilde{Y}^\infty_U$,  but the quotient
$\CO_{F,(p),+}^\times/(\CO_F^\times \cap U)^2$ need not act freely; more precisely, the
stabilizer of the cusp $B_1(\CO_{F,(p)}) g U$ of $\widetilde{Y}_U$ is the group
$\det(gUg^{-1} \cap B(F))_+  \subset \CO_{F,+}^\times$, in which $(\CO_F^\times \cap U)^2$
has finite index.

We also have a natural bijection between $Y^\infty_U$ and the set of isomorphism classes of data
$(H, I, [\lambda], [\eta])$ where
\begin{itemize}
\item $H$ is a projective $\CO_F$-module of rank two;
\item $I$ is an invertible submodule of $H$ such that $J:=H/I$ is invertible;
\item  $[\lambda]$ is a prime-to-$p$ orientation of $\wedge^2_{\CO_F}H \cong I \otimes_{\CO_F}J = IJ$,
by which we mean an $\CO_{F,(p),+}^\times$-orbit of isomorphisms 
$$\lambda:  \wedge^2_{\CO_{F,(p)}} H_{(p)}  \stackrel{\sim}{\longrightarrow}
    \CO_{F,(p)}.$$
\item $[\eta]$ is a level $U^p$-structure on $H$, i.e., a $U^p$-orbit
of $\widehat{\CO}_F^{(p)}$-linear isomorphisms 
$$\eta :  (\widehat{\CO}_{F}^{(p)})^2  \stackrel{\sim}{\longrightarrow}  \widehat{\CO}_F^{(p)} \otimes_{\CO_F} H.$$
\end{itemize}
The bijection is defined by associating the data $(H_g, I_g, [\lambda_g], [\eta_g])$ to the coset
$B(\CO_{F,(p)})_+ g U^p$, where $H_g = \widehat{\CO}_F^2 g^{-1} \cap F^2$, $I_g$
is its intersection with the subspace $\{0\} \times F$, $\lambda_g$ is induced by the determinant, 
and $\eta_g$ is induced by right pre-multiplication by $g^{-1}$.  

Note that to give a prime-to-$p$ orientation of $\wedge^2_{\CO_F}H$ is equivalent to giving
an $F_+^\times$-orbit of isomorphisms $\wedge_F^2 (\QQ\otimes H) \stackrel{\sim}{\longrightarrow} F$,
but the integrality condition is imposed for consistency with the fact that we have a bijection between
$\widetilde{Y}^\infty_U$ and the set of isomorphism classes of data of the form $(H, I, \lambda, [\eta])$,
which is similarly defined, $\CO_{F,(p),+}^\times$-equivariant, and compatible in the obvious sense
with the bijection describing $Y^\infty_U$.  In particular if $U = U(1) = \GL_2(\widehat{\CO}_F)$,
then the map sending $(H, I, [\lambda], [\eta])$ to the pair $(\wedge^2_{\CO_F}H,I)$ defines a bijection
between $Y^\infty_U$ and $C_F^+ \times C_F$, where $C_F^{(+)}$ denotes the (strict) class group of $F$.
For each such cusp we choose a polyhedral cone decomposition as in \cite[Lemme~4.2]{rap} (with $U' = U_+$
in the notation there) for $(M^*\otimes \RR)_+ \cup \{0\}$, where 
\begin{equation}\label{eqn:M} \begin{array}{l}
M = \Hom_{\CO_F}(I,\gd^{-1} J) = \gd^{-1}I^{-1}J,\\
 M^* = \Hom(M,\ZZ) \cong \Hom_{\CO_F}(J,I) = J^{-1}I,\end{array}\end{equation}
  and the positivity is induced by the orientation of 
$I\otimes_{\CO_F} J$.  

Suppose now that $U = U(N)$ for some $N \ge 3$ (not divisible by $p$) and that $\CO$ contains the $N$th roots of unity.
The above choice of cone decomposition (for the image of each cusp of $\widetilde{Y}^\infty_U$ in $Y^\infty_{U(1)}$) yields
a toroidal compactification\footnote{Compactified in the sense that its (infinitely many) connected components are proper over $\CO$.}
$\widetilde{Y}_U \hookrightarrow \widetilde{Y}_U^\tor$ such that the set of (geometrically) connected components of its (reduced) closed
subscheme $\widetilde{Z}_U^\tor := \widetilde{Y}_U^\tor - \widetilde{Y}_U$ is identified with $\widetilde{Y}^\infty_U$. 
The construction of
$\widetilde{Y}_U^\tor$ identifies its completion along the component corresponding to a cusp $\widetilde{\calC}$ represented by
$(H,I,\lambda,[\eta])$ with the quotient of a formal scheme\footnote{The formal scheme depends on the chosen cone decomposition
$\{\sigma_\alpha^{\widetilde{\calC}}\}$ and is denoted $S_N(\{\sigma_\alpha^{\widetilde{\calC}}\})^\wedge$ in \cite[3.4.2]{chai}.}
$\widehat{S}_{\widetilde{\calC}}$ by an action of $V_N^2$, where $V_N = \ker(\CO_F^\times \to (\CO_F/N\CO_F)^\times)$.
Furthermore this extends to an action of $\CO_{F,+}^\times$ on $\widehat{S} = \widehat{S}_{\widetilde{\calC}}$, and we have
an isomorphism
\begin{equation}\label{eqn:bigcompletion} \Gamma(\widehat{S}, \CO_{\widehat{S}}) \cong \CO[[q^m]]_{m \in N^{-1}M_+ \cup \{0\}}\end{equation}
compatible with the obvious action of $\CO_{F,+}^\times$ on the target.  (The isomorphism depends on a choice of splitting of
the exact sequence
$$0 \to I \to H \to J \to 0$$
of $\CO_F$-modules; modifying the splitting by an element $\beta \in J^{-1}I \cong M^*$
alters it by composition with the automorphism defined by $q^m \mapsto \zeta_N^{-\beta(m)}q^m$ for
$m \in N^{-1}M$.)   We let $\xi: \wS \to \widetilde{Y}_U$ denote the natural morphism
of formal schemes, and we write $F_\wS$ for the field of fractions of $\Gamma(\wS,\CO_\wS)$
and $\mu_\nu$ for the automorphism of $\wS$ defined by $\nu \in \CO_{F,+}^\times$.

The construction of the toroidal compactification also extends the universal abelian scheme $A$ to a semi-abelian scheme
$A^\tor$ whose pull-back to $\widehat{S}$ is identified with that of
the Tate semi-abelian scheme\footnote{More precisely, the formal scheme $\widehat{S}$ has an open cover by affine
formal subschemes $\Spf R_\sigma$ (indexed by cones $\sigma$) such that $\Spec R_\sigma \times_{\widetilde{Y}_U^\tor} A^\tor$
is identified with the semi-abelian scheme $T_{I,J}$ over $\Spec R_\sigma$.  The compatibilities in the discussion that follows
are then systematically checked by verifying them over the open subschemes $\Spec R_\sigma^0 = 
\Spec R_\sigma \times_{\widetilde{Y}_U^\tor} \widetilde{Y}_U$.}  associated to a quotient of the form
\begin{equation} \label{eqn:Tate} T_{I,J} := (\gd^{-1}I \otimes \G_m) /\widetilde{q}^{\gd^{-1}J},\end{equation}
where $\widetilde{q}^\cdot :\gd^{-1}J \to F_\wS^\times \otimes \gd^{-1}I$ is the homomorphism corresponding to the 
tautological element under the canonical isomorphism
$$\Hom(M,F_\wS^\times) =   \Hom_{\CO_F}(\gd^{-1}J, \gd^{-1}I \otimes F_\wS^\times).$$
Similarly its dual $A^\vee$ extends to a Tate semi-abelian scheme $(A^\vee)^\tor$ whose
pull-back via $\xi$ is associated to
$T_{\gd J^{-1},\gd I^{-1}}$, with the isomorphism $\gc\gd\otimes_{\CO_F}A^\tor \to (A^\vee)^\tor$
defined by the quasi-polarization pulling back to the composite
$$\gc\gd\otimes_{\CO_F} T_{I,J}  \stackrel{\sim}{\longrightarrow}  
\gd(IJ)^{-1} \otimes_{\CO_F} T_{I,J} \stackrel{\sim}{\longrightarrow} 
T_{\gd J^{-1},\gd I^{-1}},$$
where $\gc = \{\,\alpha\in F\,|\,\alpha\lambda(IJ)\subset \CO_F\,\}$, the first morphism is the isomorphism induced by $\lambda$, and
the second is the canonical one. 

The subschemes $Z_\theta$ of $\widetilde{Y}_U$ (defined in \S\ref{sec:stratification} by the vanishing of the partial Hasse invariants $H_\theta$)
are closed in $\widetilde{Y}_U^\tor$, and we let $\widetilde{Y}_U^\ord$
(resp.~$\widetilde{Y}_U^\tord$) denote the complement of their union, i.e., the ordinary locus, in $\widetilde{Y}_U$
(resp.~$\widetilde{Y}_U^\tor$), and we use similar notation for the restrictions of $A^\tor$ and $(A^\vee)^\tor$.
Since the sheaf $\Lie(A^\tord/\widetilde{Y}_U^\tord)$  is locally free over
$\CO_F \otimes \CO_{\widetilde{Y}_U^\tord}$, the universal filtration $\CF_\tau^{(j)}$ on
$$(s_*\Omega_{A/\widetilde{Y}_U}^1)_\tau  \cong \Shom_{\CO_{\widetilde{Y}_U}}(\Lie(A/\widetilde{Y}_U)_\tau,\CO_{\widetilde{Y}_U})$$
extends canonically to one on $\Shom_{\CO_{\widetilde{Y}^\tor_U}}(\Lie(A^\tor/\widetilde{Y}^\tor_U)_\tau,\CO_{\widetilde{Y}^\tor_U})$
for each $\tau \in \Sigma_0$.  Furthermore its pull-back to $\wS$ is identified (in the notation of (\ref{eqn:st})) with
$$0 \subset  t_{\tau,1}(I^{-1}\otimes \CO_\wS)_\tau \subset t_{\tau,2}(I^{-1}\otimes \CO_\wS) \subset \cdots 
       \subset  t_{\tau,e_\gp}(I^{-1} \otimes \CO_\wS)_\tau = (I^{-1}\otimes \CO_\wS)_\tau$$
 under the canonical isomorphism
 $$\Shom_{\CO_{\wS}}(\Lie(T_{I,J}/\wS),\CO_\wS)  \cong \Hom(\gd^{-1}I, \CO_\wS) \cong I^{-1}\otimes\CO_\wS.$$
 We thus obtain extensions ${\widetilde{\CL}}_\theta^\tor$ of the line bundles ${\widetilde{\CL}}_\theta = {\widetilde{\CL}}_{\gp,i,j}$
 to $\widetilde{Y}_U^\tor$ whose pull-back to $\wS$
 is identified with $(I^{-1})_\theta \otimes_{\CO}\CO_\wS$, where $(I^{-1})_\theta$ is defined by (\ref{eqn:graded}).
  
 Similarly $\Lie((A^\vee)^\tord/\widetilde{Y}_U^\tord)$ is locally free over $\CO_F \otimes \CO_{\widetilde{Y}_U^\tord}$, but
 the line bundles ${\widetilde{\CM}}_\theta = {\widetilde{\CM}}_{\gp,i,j}$ over $\widetilde{Y}_U^\ord$ are canonically identified with
 $$(R^1s_*\CO_{A^\ord})_\tau[u-\theta(\varpi_\gp)]  \cong \Lie((A^\vee)^\ord/\widetilde{Y}_U^\ord)_\tau[u-\theta(\varpi_\gp)].$$
 It follows that each ${\widetilde{\CM}}_\theta$ extends to a line bundle ${\widetilde{\CM}}_\theta^\tor$ on $\widetilde{Y}_U^\tor$ such that
 the identification
  $$\Lie(T_{\gd J^{-1},\gd I^{-1}}/\wS) = J^{-1}\otimes\CO_\wS$$
induces an isomorphism
 $$\xi^*{\widetilde{\CM}}_\theta^\tor \cong
 (J^{-1}\otimes \CO_{\wS})_\tau[u-\theta(\varpi_\gp)]  \cong (\gd J^{-1} \otimes \CO_\wS)_\tau  \otimes_{\CO[u],\theta} \CO.$$
 We can thus identify the pull-back $\xi^*{\widetilde{\CN}}_\theta^\tor$ of the line bundle 
 ${\widetilde{\CN}}_\theta^\tor = {\widetilde{\CL}}_\theta^\tor \otimes_{\CO_{\widetilde{Y}_U^\tor}} {\widetilde{\CM}}_\theta^\tor$
 with $(\gd (IJ)^{-1})_\theta\otimes_{\CO} \CO_\wS$, which the polarization in turn identifies 
 with $(\gc\gd)_\theta \otimes_{\CO} \CO_\wS$ in the notation of (\ref{eqn:graded}).
 Finally it follows that the automorphic bundles ${\widetilde{\CA}}_{\bk,\bl}$ extend to line bundles ${\widetilde{\CA}}_{\bk,\bl}^\tor$ on $\widetilde{Y}_U^\tor$
 such that 
 \begin{equation} \label{eqn:Dline}
 \xi^*{\widetilde{\CA}}_{\bk,\bl}^\tor \cong  D_{\bk,\bl} \otimes_{\CO} \CO_\wS, \quad\mbox{where\,\,}
 D_{\bk,\bl} := \bigotimes_{\theta\in \Sigma} \left( (I^{-1})^{\otimes k_\theta}_\theta \otimes (\gd(IJ) ^{-1})^{\otimes l_\theta}_\theta\right)
 \end{equation}
 (the tensor products being over $\CO$).  We refer to this isomorphism as 
 the {\em canonical trivialization} of  $\xi^*{\widetilde{\CA}}_{\bk,\bl}^\tor$.
 
Next we consider the completion of $\widetilde{Y}_U^\tor$ along the component corresponding to the cusp $\widetilde{\calC}$
represented by $(H,I,\lambda,[\eta])$,
which we denote $(\widetilde{Y}_{U^\tor})_{\widetilde{\calC}}^\wedge$.  
We now describe the global sections of the completions of the line bundles ${\widetilde{\CA}}_{\bk,\bl}^\tor$
using the identification $(\widetilde{Y}_{U}^\tor)_{\widetilde{\calC}}^\wedge = \wS/V_N^2$ and taking invariants under the action of
$V_N^2$ on their trivializations over $\wS$.  Note firstly that
$\Gamma( (\widetilde{Y}_U^\tor)_{\widetilde{\calC}}^\wedge,  \CO_{(\widetilde{Y}_U^\tor)_{\widetilde{\calC}}^\wedge})
  = \Gamma(\wS,\CO_\wS)^{V_N^2}$ corresponds to
\begin{equation} \label{eqn:completion}
\left\{\, \left.\sum r_mq^m \in \CO[[q^m]]_{m \in N^{-1}M_+ \cup \{0\}} \,\right|\, r_{\alpha^2m} = r_m \,\,\forall \alpha \in V_N, m \in N^{-1}M_+\,\right\}
\end{equation}
  under the isomorphism of (\ref{eqn:bigcompletion}).
One then finds that
the descent data for $\xi^*A^\tor$ is provided by the isomorphisms $T_{I,J} \stackrel{\sim}{\longrightarrow} \mu_{\alpha^2}^*T_{I,J}$
induced by $\alpha \otimes 1$ on $\gd^{-1}I \otimes \G_m$, from which it follows that the descent data for $\xi^*{\widetilde{\CL}}_{\theta}^\tor$ is
provided on the trivialization by the isomorphisms
$$(I^{-1})_\theta \otimes_\CO \mu^*_{\alpha^2} \CO_{\wS}  \stackrel{\sim}{\longrightarrow} (I^{-1})_\theta \otimes_\CO \CO_{\wS}$$
induced by $\theta(\alpha)$ on $(I^{-1})_\theta$.  On the other hand the descent data for $\xi^*{\widetilde{\CM}}_\theta^\tor$ is similarly induced
on the canonical trivialization by $\theta(\alpha)^{-1}$, so that the resulting trivialization of $\xi^*{\widetilde{\CN}}_\theta^\tor$ descends to
$(\widetilde{Y}_{U^\tor})_{\widetilde{\calC}}^\wedge$ (in fact extending the one already defined over $\widetilde{Y}_U$ via the choice of generator
$t_{\tau,j}(f'(\varpi_\gp) \otimes 1)$ of $(\gc\gd)_\theta$).  Since
$\Gamma( (\widetilde{Y}_U^\tor)_{\widetilde{\calC}},  ({\widetilde{\CA}}^\tor_{\bk,\bl})_{\widetilde{\calC}}^\wedge) =
\Gamma(\wS,\xi^*{\widetilde{\CA}}^\tor_{\bk,\bl})^{V_N^2}$, we conclude the following:
\begin{proposition} \label{prop:toroidal}
Suppose that $U=U(N)$ and $\CO$ contains the $N$th roots of unity. 
Then the  isomorphism~(\ref{eqn:bigcompletion}) and the canonical trivialization~(\ref{eqn:Dline}) identify 
$\Gamma( (\widetilde{Y}_U^\tor)_{\widetilde{\calC}},  ({\widetilde{\CA}}^\tor_{\bk,\bl})_{\widetilde{\calC}}^\wedge)$ with
 $$\left \{\, \left.\sum_{m \in N^{-1}M_+ \cup \{0\}}  b \otimes r_mq^m \,\right|\, 
   r_{\alpha^2m} = \chi_{\bk}(\alpha) r_m\,\, \forall \alpha \in V_N, m \in N^{-1}M_+\,\right\},$$
   where $b$ is any choice of basis for $D_{\bk,\bl}$.
\end{proposition}

\subsection{Minimal compactifications} \label{sec:min}
We now recall the construction due to Chai~\cite{chai} of minimal compactifications of Hilbert modular varieties.
The presentation in~\cite{chai} is very concise with numerous typos, but a more detailed treatment of the construction can
be found in \cite{dim} in the case of $U_1(\gn)$ (with different conventions than ours), and of the descriptions of $q$-expansions
in that case in \cite{dem}.

We continue to assume for the moment that $U = U(N)$ for some sufficiently large $N$ prime to $p$.
The minimal compactification $\widetilde{Y}_U \hookrightarrow \widetilde{Y}_U^{\min}$ is then constructed as in \cite[\S4]{chai}
or \cite[\S8]{dim}.   More precisely, letting $\bt = \sum \be_\theta$  and taking the global sections of  $\oplus_{k \ge 0} {\widetilde{\CA}}^\tor_{k\bt,\bf{0}}$
over each component of $\widetilde{Y}_U^\tor$ yields a projective scheme over $\CO$ containing the corresponding component
of the Deligne-Pappas model as an open subscheme.  Gluing their ordinary loci to $\widetilde{Y}_U$ along $\widetilde{Y}_U^\ord$ yields a scheme
$\widetilde{Y}_U^{\min}$ and a proper morphism $\pi: \widetilde{Y}_U^\tor \to \widetilde{Y}_U^{\min}$ such that ${\widetilde{\iota}}: \widetilde{Y}_U \to \widetilde{Y}_U^{\min}$
is an open immersion.  Furthermore the (reduced) complement $\widetilde{Y}_U^{\min} - \widetilde{Y}_U$ is an infinite
disjoint union of copies of $\Spec\CO$ indexed by $\widetilde{Y}^\infty_U$, the preimage of each in $\widetilde{Y}_U^\tor$ being the 
corresponding connected component of $\widetilde{Z}_U^\tor$, and the scheme $\widetilde{Y}_U^{\min}$ is independent of the choice
of cone decomposition in the construction of $\widetilde{Y}_U^\tor$.  

Now recall that the Koecher Principle implies that 
$\CO_{\widetilde{Y}_U^{\min}} = \pi_*\CO_{\widetilde{Y}_U^\tor} = {\widetilde{\iota}}_*\CO_{\widetilde{Y}_U}$, so that 
$\CO^\wedge_{\widetilde{Y}_U^{\min},\widetilde{\calC}} =
\Gamma( (\widetilde{Y}_U^\tor)_{\widetilde{\calC}}^\wedge,  \CO_{(\widetilde{Y}_U^\tor)_{\widetilde{\calC}}^\wedge})$
is the ring described by (\ref{eqn:completion}), 
where we have written $\widetilde{\calC}$ for the corresponding point of $\widetilde{Y}_U^{\min}$.
Furthermore the argument of \cite[Prop.~4.9]{rap} 
shows that ${\widetilde{\iota}}_*{\widetilde{\CA}}_{\bk,\bl} = \pi_*{\widetilde{\CA}}_{\bk,\bl}^\tor$
(see the discussion following~\cite[Def.~6.10]{rap}, or view $\widetilde{Y}_U$ as a disjoint union of
PEL Shimura varieties and apply~\cite[Thm.~2.5]{lan2}), so the Theorem on Formal Functions gives that
$(\iota_*{\widetilde{\CA}}_{\bk,\bl})^\wedge_{\widetilde{\calC}} = 
\Gamma( (\widetilde{Y}_{U^\tor})_{\widetilde{\calC}}^\wedge,  ({\widetilde{\CA}}^\tor_{\bk,\bl})_{\widetilde{\calC}}^\wedge))$ is the
$\CO^\wedge_{\widetilde{Y}_U^{\min},\widetilde{\calC}}$-module described in Proposition~\ref{prop:toroidal}.
(Note that ${\widetilde{\iota}}_*{\widetilde{\CA}}_{\bk,\bl}$ is coherent,
but not necessarily invertible.)  

Similarly for any $\CO$-algebra $R$, we may identify
$({\widetilde{\iota}}_{R,*}{\widetilde{\CA}}_{\bk,\bl,R})^\wedge_{\widetilde{\calC}}$ with
\begin{equation} \label{eqn:AYtilde}
\left\{\,\left. \sum_{m \in N^{-1}M_+ \cup \{0\}} b \otimes r_mq^m  \,\right|\, 
   r_{\alpha^2m} = \chi_{\bk,R}(\alpha) r_m \,\,\forall \alpha \in V_N, m \in N^{-1}M_+\,\right\}
   \end{equation}
 as a module over 
$\CO^\wedge_{\widetilde{Y}_{U,R}^{\min},\widetilde{\calC}}$, which the Koecher Principle and (\ref{eqn:bigcompletion}) identify with
\begin{equation} \label{eqn:RYtilde}
\left\{\, \left.\sum r_mq^m \in R[[q^m]]_{m \in N^{-1}M_+ \cup \{0\}} \,\right|\, r_{\alpha^2m} = r_m\,\, \forall \alpha \in V_N, m \in N^{-1}M_+\,\right\},
\end{equation}
where ${\widetilde{\iota}}_R: \widetilde{Y}_{U,R} \to \widetilde{Y}_{U,R}^{\min}$ is the base-change of ${\widetilde{\iota}}$ to $R$, the completions
are at the fibre over $\widetilde{\calC}$, and $b$ is any basis for $D_{\bk,\bl}$.

The compatibility of the choices of polyhedral cone decompositions ensures that the
natural action of $\CO_{F,(p),+}$ on $\widetilde{Y}_U$ extends (uniquely) to one on $\widetilde{Y}_U^\tor$.
Furthermore the stabilizer of each component of $\widetilde{Z}_U^\tor$ is $V_{N,+}$, and the action of
$V_{N,+}$ on each completion $(\widetilde{Y}_{U^\tor})_{\widetilde{\calC}}^\wedge = \wS/V_N^2$ is induced
by an action of $V_{N,+}$ on $\wS$ such that the effect of $\nu \in V_{N,+}$ on global sections of
$\CO_{\wS}$ is induced by multiplication by $\nu^{-1}$ on $M$.  We see also that the canonical isomorphism $A \to \nu^*A$
extends to an isomorphism $A^\tor \to \nu^*A^\tor$ whose pull-back via $\xi$ is induced by the
identity on $\gd^{-1}I \otimes \G_m$, from which it follows that the action of $\nu$ is compatible
with the canonical trivialization of the line bundle ${\widetilde{\CL}}_\theta^\tor$ over $\wS$.  On the other
hand the induced isomorphisms $\nu^*{\widetilde{\CM}}_\theta^\tor \to \CM_\theta^\tor$ and $\nu^*{\widetilde{\CN}}_\theta^\tor
\to {\widetilde{\CN}}_\theta^\tor$ pull back to ones corresponding to multiplication by $\theta(\nu)$.
 
Since $\widetilde{Y}_U^{\min}$ is a disjoint union of projective schemes over $\CO$ on which
$\CO_{F,(p),+}^\times/V_N^2$ acts with finite stabilizers, the quotient scheme exists, and we define
this to be the minimal compactification $Y_U^{\min}$ of $Y_U$.  We thus obtain an open
immersion $\iota:Y_U \to Y_U^{\min}$ such that $Y_U^{\min}$ is projective over $\CO$
and the (reduced) complement of $Y_U$ is a disjoint union of copies of $\Spec\CO$ in
canonical bijection with the set of cusps $Y_U^\infty$.  Furthermore we again have that
$\iota_*\CO_{Y_U} = \CO_{Y_U^{\min}}$, and its completion $\CO^\wedge_{Y_U^{\min},\calC}$
at the cusp $\calC$ represented by $(H,I,[\lambda],[\eta])$ is identified under (\ref{eqn:completion}) with
$$\left\{\, \left.\sum r_mq^m \in \CO[[q^m]]_{m \in N^{-1}M_+ \cup \{0\}} \,\right|\, r_{\nu m} = r_m\,\, \forall \nu \in V_{N,+}, m \in N^{-1}M_+\,\right\}$$
(where the identification depends as in (\ref{eqn:bigcompletion}) on a choice of splitting of the exact sequence $0 \to I \to H \to J \to 0$).
Now suppose that $R$ is an $\CO$-algebra such that $\chi_{\bk+2\bl,R}$ is trivial on $V_N$,
so that the line bundle $\widetilde{\CA}_{\bk,\bl,R}$ descends to one over $Y_{U,R}$ which we denote by $\CA_{\bk,\bl,R}$.
We then see that $\iota_*\CA_{\bk,\bl,R}$ is a coherent sheaf on $Y_{U,R}^{\min}$ whose completion at the (base-change to $R$ of the)
cusp $\calC$ is identified under (\ref{eqn:AYtilde}) with
$$ \left\{\, \left.\sum_{m \in N^{-1}M_+ \cup \{0\}} b \otimes r_mq^m \,\right|\, 
   r_{\nu^{-1} m} = \chi_{\bl,R}(\nu)  r_m\,\, \forall \nu \in V_{N,+}, m \in N^{-1}M_+\,\right\}.$$
In particular $\iota_*\CA_{\bk,\bl,R}$ is a line bundle if $\chi_{\bl,R}$ is trivial on $V_{N,+}$.

Suppose now that $U'$ is any sufficiently small open compact subgroup of $\GL_2(\widehat{\CO}_F)$ containing
$\GL_2(\CO_{F,p})$.  One can then carry out a construction similar to the one above to obtain the minimal compactification,
or choose an $N$ prime to $p$ such that $U(N) \subset U'$, extend the natural (right) action of $U'/U(N)$ on $Y_{U(N)}$
to $Y_{U(N)}^{\min}$ and take the quotient; we do the latter (see \cite{dim} for the former in the case of $U' = U_1(\gn)$).  
Firstly our choice of polyhedral cone decompositions ensures that the natural right action of $U'/U$ on $\widetilde{Y}_{U}$
extends to $\widetilde{Y}_{U}^\tor$,  where $U = U(N)$ for some choice of $N$ as above.
Denoting the resulting automorphism of $\widetilde{Y}_{U}^\tor$ by
$\widetilde{\rho}_g$ for $g \in U'$, the canonical identification of the universal $A$ over $\widetilde{Y}_{U}$ with its pull-back extends
to an identification $A^\tor = \widetilde{\rho}_g^*A^\tor$, giving rise to canonical isomorphisms 
$\widetilde{\rho}_g^*{\widetilde{\CA}}^\tor_{k\bt,\bf{0}} = {\widetilde{\CA}}^\tor_{k\bt,\bf{0}}$, and hence to an action of $U'/U$
on $\widetilde{Y}_{U}^{\min}$ extending its action on $\widetilde{Y}_{U}$.  Moreover the action commutes with the natural
action of $\CO_{F,(p),+}$, so it descends to an action on $Y_{U}^{\min}$ extending the action on $Y_{U}$.
We denote the resulting automorphisms of $Y_U^{\min}$ by $\rho_g$,
and define $Y_{U'}^{\min}$ to be the quotient of $Y_{U}^{\min}$ by the action of $U'/U$ (which we will shortly
see is independent of the choice of $N$ in its definition).  

Identifying the set of components of $Y_{U}^{\min} - Y_{U}$ with $Y^\infty_{U}$, the resulting action of
$g \in U'$ is given by pre-composing $\eta$ with right multiplication by $g^{-1}$, so the set of components of $Y_{U'}^{\min} - Y_{U'}$
may be identified with $Y^\infty_{U'}$.  For each cusp $\calC' \in Y^\infty_{U'}$, the completion $\CO^\wedge_{Y_{U'}^{\min},\calC'}$
is identified with the subring of $\CO^\wedge_{Y_{U}^{\min},\calC}$ invariant under the stabilizer of $\calC$
in $U'/U$, where $\calC$ is any cusp of $Y_{U}$ in the preimage of $\calC'$.  Choose such an $(H,I,[\lambda],[\eta])$
representing $\calC$ and a splitting $\sigma: H \stackrel{\sim}{\longrightarrow} J \times I$, and let 
\begin{equation}\label{eqn:Isotropy} \Gamma_{\calC} = \left\{\, \left.\left(\begin{array}{cc}\alpha&\beta\\0&\delta\end{array}\right)\,\right|\,
  \alpha \in \CO_F^\times, \beta \in J^{-1}I, \delta \in \alpha\CO_{F,+}^\times\,\right\},\end{equation}
which we view as acting on $J\times I$ by right multiplication.  The stabilizer of $\calC$ is then the set of classes 
$Ug = gU \in U'/U$ such that 
$$g \equiv \eta^{-1}\sigma^{-1}\gamma\sigma\eta \bmod N\widehat{\CO}_F\mbox{\,\,for some\,\,}
\gamma \in \Gamma_{\calC},$$
and we let $\Gamma_{\calC,U'} = \Gamma_{\calC} \cap \sigma\eta U' \eta^{-1}\sigma^{-1}$.  Thus the stabilizer of $\calC$
is the image of the homomorphism $s: \Gamma_{\calC,U'} \longrightarrow U'/U$ defined by
$\gamma \mapsto \eta^{-1}\sigma^{-1}\gamma\sigma\eta U$.

We claim that if $g = s(\gamma)$, then $\rho_g^*$ on $\CO^\wedge_{Y_{U}^{\min},\calC} = H^0(\wS,\CO_{\wS})^{V_{N,+}}$
is induced by an automorphism $\psi_\gamma$
of $\widehat{S}$ whose effect on global sections corresponds to the map defined by
\begin{equation} \label{eqn:action} \psi_\gamma^*: q^m \mapsto \zeta_N^{-\beta(\alpha^{-1}Nm)}q^{\alpha^{-1}\delta m} \end{equation}
under (\ref{eqn:bigcompletion}) and the identification $M^* = J^{-1}I$ of (\ref{eqn:M}).
Indeed letting $\nu$ denote $\alpha\delta$ (as well as the automorphism of $\widetilde{Y}_U^\tor$ defined
by its effect on the universal polarization), we see that $\delta\otimes 1$ on $T_{I,J}$
defines an isomorphism $\xi^*\rho_g^*\nu^*A^\tor \stackrel{\sim}{\longrightarrow} \psi_\gamma^*\xi^* A^\tor$
compatible with all auxiliary data, from which one deduces that $\nu\circ\rho_g\circ\xi = \xi\circ \psi_\gamma$.
Note also that (\ref{eqn:action}) defines an action of $\Gamma_{\calC}$ on $\Gamma(\wS,\CO_{\wS})$
which factors through the surjection
$$\begin{array}{ccc} \Gamma_{\calC} &\longrightarrow & (J^{-1}I\otimes  \ZZ/N\ZZ)\rtimes \CO_{F,+}^\times \\
\left(\begin{array}{cc}\alpha&\beta\\0&\delta\end{array}\right)&\mapsto & (-\alpha^{-1}\beta,\alpha^{-1}\delta),\end{array}$$
and the latter group acts on $\Gamma(\wS,\CO_{\wS})^{V_{N,+}}$ via its quotient
$$(J^{-1}I\otimes  \ZZ/N\ZZ) \rtimes  (\CO_{F,+}^\times/V_{N,+}).$$
We conclude that $\CO^\wedge_{Y_{U'}^{\min},\calC'} = H^0(\wS,\CO_{\wS})^{\Gamma_{\calC,U'}}$ is identified with
$$\left\{\, \sum_{m \in N^{-1}M_+ \cup \{0\}}  r_mq^m \,\left|\, r_{\alpha^{-1}\delta m} = \zeta_N^{-\beta(\alpha^{-1}Nm)}
 r_m\,\, \forall m \in N^{-1}M_+, \smallmat{\alpha}{\beta}{0}{\delta} \in \Gamma_{\calC,U'}
 \,\right.\right\},$$
where we recall that the isomorphism may depend on the choice of the splitting $\sigma$ and that we view $\beta$
as an element of $M^*$.
We note in particular if $U' = U(N')$ for some $N'|N$, then the resulting description of $\CO^\wedge_{Y_{U'}^{\min},\calC'}$
coincides with the one previously obtained, from which it follows that the same holds for the scheme $Y_{U'}^{\min}$,
and hence that $Y_{U'}^{\min}$ is independent of the choice of $N$ in its definition (for any sufficiently small $U'$ containing
$\GL_2(\CO_{F,p})$.)  

Suppose now that $\bk, \bl \in \ZZ^\Sigma$ and $R$ is an $\CO$-algebra such that $\chi_{\bk+2\bl,R}$ is trivial on
$U' \cap \CO_F^\times$, and consider the automorphic bundle $\CA'_{\bk,\bl,R}$ on $Y_{U',R}$.  Letting 
$\iota_R'$ denote the open immersion of $Y_{U',R}$ in $Y_{U',R}^{\min}$, similar considerations to those above show that 
$\iota'_{R,*}\CA'_{\bk,\bl,R}$ is a coherent sheaf on $Y_{U,R}^{\min}$ whose completion at $\calC'$ is identified with
the $\CO^\wedge_{Y_{U',R}^{\min},\calC'}$-module of $\Gamma_{\calC,U'}$-invariants in
$({\widetilde{\iota}}_{R,*}{\widetilde{\CA}}_{\bk,\bl,R})^\wedge_{\calC}$.  Using that the isomorphism
$\xi^*A^\tor \stackrel{\sim}{\longrightarrow} \psi_\gamma^*\xi^* A^\tor$ is induced by $\delta \otimes 1$ on the
Tate semi-abelian scheme, we find that the resulting automorphism multiplies the canonical trivialization (\ref{eqn:Dline})
of $\xi^*\widetilde{\CA}^\tor_{\bk,\bl,R}$ by $\chi_{\bl,R} (\alpha)\chi_{\bk+\bl,R}(\delta)$.
We therefore conclude:
\begin{proposition} \label{prop:minimal}  If $\chi_{\bk+2\bl,R}$ is trivial on $U' \cap \CO_F^\times$, then 
$\iota'_{R,*}\CA'_{\bk,\bl,R}$ is a coherent sheaf on $Y_{U',R}^{\min}$ whose completion at (the fibre over)
$\calC'$ is identified by the Koecher Principle and Proposition~\ref{prop:toroidal} with
$$\left \{\, \sum_{m \in N^{-1}M_+ \cup \{0\}}  b \otimes r_mq^m \,\left|\, \begin{array}{c}
   r_{\alpha^{-1}\delta m} =  \zeta_N^{-\beta(\alpha^{-1}Nm)}
 \chi_{\bl,R} (\alpha)\chi_{\bk+\bl,R}(\delta) r_m\\ \mbox{\rm for all $m \in N^{-1}M_+,
 \smallmat\alpha\beta0\delta \in \Gamma_{\calC,U'}$} \end{array}
 \,\right.\right\}.$$
 \end{proposition}
 Note that the description of $\CO^\wedge_{Y_{U',R}^{\min},\calC'} = (\iota'_{R,*} \CO_{Y_{U',R}^{\min}})^\wedge_{\calC'}$
 may be viewed as a special case (with $\bk = \bl = \bf{0}$), as can the prior formula for $U = U(N)$.  Furthermore
 the identifications are compatible in the obvious senses with base changes $R \to R'$, inclusions
 $U'' \subset U'$ (provided the splittings $\sigma$ are chosen compatibly), and the natural algebra structure on
 $\bigoplus_{\bk,\bl} \CA'_{\bk,\bl,R}$ (taking the direct sum over $\bk,\bl$ as in the statement).

 Recall that the $q$-expansion Principle allows one to characterize Hilbert modular forms in terms
 of their $q$-expansions:
 \begin{proposition} \label{prop:qexp}
If $C \subset Y^\infty_{U'}$ is any set of cusps containing at least one on each component of $Y_{U'}$,
then the natural map
$$M_{\bk,\bl}(U';R) = H^0(Y_{U'}^{\min}, \iota'_{R,*}\CA'_{\bk,\bl,R}) \longrightarrow \bigoplus_{\calC' \in C}
    (\iota'_{R,*}{\CA'}_{\bk,\bl,R})^\wedge_{\calC'}$$
 is injective.
 \end{proposition}

 Note also that we may replace
 $D_{\bk,\bl} \otimes_{\CO}\cdot$ with  $D_{\bk,\bl,R} \otimes_{R}\cdot$ in the description of $q$-expansions over $R$.
 In particular if $R$ is an $\FF$-algebra, the identification
 $$\begin{array}{rl} (I^{-1})_\theta \otimes_{\CO} \FF& = (I^{-1}\otimes\CO)_\tau \otimes_{\CO[u],\theta} \FF  = 
 t_{\tau,j}(I_{\gp}^{-1}\otimes_{W(\CO_F/\gp),\tau}\CO)\otimes_{\CO[u],\theta} \FF\\
 &=  t_{\tau,j}(I_{\gp}^{-1}\otimes_{W(\CO_F/\gp),\tau}\CO)/(u, \gm_\CO)t_{\tau,j}(I_{\gp}^{-1}\otimes_{W(\CO_F/\gp),\tau}\CO)\\
 &= u^{e_\gp - j}(I_{\gp}^{-1}\otimes_{W(\CO_F/\gp),\tau}\FF)/u^{e_\gp - j+1}(I_{\gp}^{-1}\otimes_{W(\CO_F/\gp),\tau}\FF)\end{array}$$
 yields a canonical isomorphism
\begin{equation}\label{eqn:modp}
(I^{-1})_\theta \otimes_{\CO} R = (\gp^{e_\gp - j}I^{-1}/\gp^{e_\gp - j+1}I^{-1}) \otimes_{\CO_F/\gp,\tau} R.\end{equation}
The analogous formula holds for the factors $(\gd(IJ) ^{-1})_\theta$ appearing in the definition of $D_{\bk,\bl}$.

 The condition on the $q$-expansion coefficients in the description of the completions in Proposition~\ref{prop:minimal}
 simplifies for certain standard level structures and cusps, as in \cite[Prop.~9.1.2]{DS}.
  Suppose that $\gn$ is an ideal of $\CO_F$ such that  $\chi_{\bk+2\bl,R}$ is trivial on
 $V_{\gn} = \ker(\CO_F^\times \to (\CO_F/\gn)^\times)$.  Letting $U' = U(\gn)$, we have
 $$\Gamma_{\calC,U'} = \left\{\, \left.\left(\begin{array}{cc}\alpha&\beta\\0&\delta\end{array}\right)\,\right|\,
  \alpha \in V_\gn, \beta \in \gn M^*,  \delta \in \alpha V_{\gn,+}\,\right\}$$
for every cusp $\calC$ of $Y_U$.
Note that $m \in \gn^{-1} M$ if and only if $\beta(N m) \in N\ZZ$ for all $\beta \in \gn M^*$, and that
$\alpha,\delta \in V_\gn$ implies that $\chi_{\bl,R}(\alpha)\chi_{\bk+2\bl}(\delta) = 
\chi_{\bl,R}(\alpha\delta^{-1})$, so we see that
$$(\iota'_{R,*}{\CA'}_{\bk,\bl,R})^\wedge_{\calC'}
 \simeq \left\{ \left.\, \sum_{m\in (\mathfrak{n}^{-1}M)_+\cup \{0\}} b \otimes  r_mq^m \,\right|\,
   r_{\nu^{-1} m} = \chi_{\bl,R}(\nu) r_m \ \mbox{for all $\nu \in V_{\mathfrak{n},+}$}\,\right\}$$
for every cusp $\calC'$ of $Y_{U'}$.

Keep the same assumption on $\gn$, but now let $U' = U_1(\gn)$ and suppose that $\calC'$ is
a cusp of $Y_{U'}$ ``at $\infty$'' in the sense that $\eta(0,1) \in I + \gn \widehat{H}^{(p)}$.
We then find that
$$\Gamma_{\calC,U'} = \left\{\, \left.\left(\begin{array}{cc}\alpha&\beta\\0&\delta\end{array}\right)\,\right|\,
  \alpha \in \delta\CO_{F,+}^\times, \beta \in M^*,  \delta \in V_{\gn}\,\right\},$$
  and we similarly conclude that
$$(\iota'_{R,*}{\CA'}_{\bk,\bl,R})^\wedge_{\calC'}
 \simeq \left\{ \left.\, \sum_{m\in M_+\cup \{0\}}b \otimes  r_mq^m \,\right|\,
   r_{\nu^{-1} m} = \chi_{\bl,R}(\nu) r_m \ \mbox{for all $\nu \in \CO_{F,+}$}\,\right\}.$$
 We remark that every component of $Y_{U'}$ contains such cusps (in the obvious sense), and
 that in this case the isomorphism is independent of the choice of splitting $\sigma$.

\subsection{Kodaira--Spencer filtration}  \label{sec:KS2}
We next explain how the Kodaira--Spencer filtration on differentials extends to compactifications.

We maintain the notation from the preceding section.  In particular, we first assume $U = U(N)$ for
some $N$ prime to $p$ before deducing results for more general level structures.  
The construction of $\widetilde{Y}_U^\tor$ via torus embeddings then yields a canonical isomorphism 
\begin{equation}\label{eqn:logdiffs} \xi^*(\Omega^1_{\widetilde{Y}_U^\tor/\CO}(\log \widetilde{Z}_U^\tor)) \cong N^{-1}M \otimes \CO_{\wS}\end{equation}
for each cusp $\tcalC$ of $\widetilde{Y}_U$ under which the descent data relative to the quotient map 
$\wS= \wS_{\tcalC} \to (\widetilde{Y}_U^\tor)_\tcalC^\wedge$
corresponds to that induced by the obvious action of $V_N^2$ on $N^{-1}M$, and the completion
of the canonical derivation
$$d:  \CO_{\widetilde{Y}_U^\tor}   \to \Omega^1_{\widetilde{Y}_U^\tor/\CO}(\log \widetilde{Z}_U^\tor)$$
pulls back to a derivation $\CO_{\wS} \to N^{-1}M \otimes \CO_{\wS}$ whose effect on global sections
corresponds under (\ref{eqn:bigcompletion}) to the map defined by
$$\sum_{m \in N^{-1}M_+ \cup \{0\}} r_mq^m \mapsto \sum_{m \in N^{-1}M_+ \cup \{0\}} m \otimes r_m q^m.$$

Recall also that 
$$\begin{array}{rcl}
s_*\Omega^1_{A^\ord/\widetilde{Y}_U^\ord} &\cong &\Shom_{\CO_{\widetilde{Y}_U^\ord}}(\Lie(A^\ord/\widetilde{Y}_U^\ord), \CO_{\widetilde{Y}_U^\ord}) \\
\mbox{and}\quad R^1s_*\CO_{A^\ord} &\cong& \Lie((A^\ord)^\vee/\widetilde{Y}_U^\ord) \end{array}$$
are locally free sheaves of $\CO_F \otimes \CO_{\widetilde{Y}_U^\ord}$-modules over $\widetilde{Y}_U^\ord$,
and therefore so is
$$s_*\Omega^1_{A^\ord/\widetilde{Y}_U^\ord} \otimes_{\CO_F\otimes\CO_{\widetilde{Y}_U^\ord}} 
    s_*\Omega^1_{(A^\ord)^\vee/\widetilde{Y}_U^\ord}$$
Decomposing this sheaf over embeddings $\tau \in \Sigma_0$ and equipping it with the filtration
defined by the images of the endomorphisms $t_{\tau,j}$ defined by (\ref{eqn:st}), we see that the successive quotients
$$\begin{array}{l} t_j\left(s_*\Omega^1_{A^\ord/\widetilde{Y}_U^\ord} \otimes_{\CO_F\otimes\CO_{\widetilde{Y}_U^\ord}} 
    s_*\Omega^1_{(A^\ord)^\vee/\widetilde{Y}_U^\ord}\right)_\tau
    \otimes_{\CO[u],\theta} \CO \\  \qquad
    \cong \Shom_{\CO_S}(R^1s_*\CO_{A^\ord}[u-\theta(\varpi)], 
    t_j(s_*\Omega^1_{(A^\ord)^\vee/\widetilde{Y}_U^\ord})_ \tau\otimes_{\CO[u],\theta} \CO)\end{array}$$
(where $\tau = \tau_i$, $t_j = t_{\tau,j}$ and $\theta = \theta_{\gp,i,j}$) are canonically
identified with the automorphic bundles $\widetilde{\CA}_{2\be_\theta,-\be_\theta}$ over $\widetilde{Y}_U^\ord$.
Furthermore the proof of Theorem~\ref{thm:KSfil} shows that the natural map
\begin{equation}\label{eqn:KSord}
s_*\Omega^1_{A^\ord/\widetilde{Y}_U^\ord} \otimes_{\CO_F\otimes\CO_{\widetilde{Y}_U^\ord}} 
    s_*\Omega^1_{(A^\ord)^\vee/\widetilde{Y}_U^\ord}   \longrightarrow \Omega^1_{\widetilde{Y}_U^\ord/\CO}\end{equation}
arising from Grothendieck--Messing theory, or equivalently the Gauss--Manin connection
on $\CH^1_\dr(A^\ord/\widetilde{Y}_U^\ord)$ (see \cite[\S2.1.7]{lan}),  is an isomorphism.  
In particular the Kodaira--Spencer filtration on $\Omega^1_{\widetilde{Y}_U^\ord/\CO}$ corresponds
under (\ref{eqn:KSord}) to the one defined by the images of the endomorphisms $t_{\tau,j}$.
Furthermore (\ref{eqn:KSord}) extends over $\widetilde{Y}_U^\tord$ to an isomorphism
$$\begin{array}{c} \Shom_{\CO_{\widetilde{Y}_U^\tord}}(\gd \otimes_{\CO_F} \Lie(A^\tord/\widetilde{Y}_U^\tord) \otimes_{\CO_F\otimes \CO_{\widetilde{Y}_U^\tord}} 
       \Lie((A^\tord)^\vee/\widetilde{Y}_U^\tord),  \CO_{\widetilde{Y}_U^\tord})\\ \stackrel{\sim}{\longrightarrow}\,\,
    \Omega^1_{\widetilde{Y}_U^\tord/\CO}(\log \widetilde{Z}_U^\tor)\end{array}$$
whose pull-back via $\xi = \xi_\tcalC$ for each cusp $\tcalC$ of $\widetilde{Y}_U$ is compatible
with the canonical isomorphisms of the pull-back of each with $M \otimes \CO_{\wS} = N^{-1}M\otimes \CO_{\wS}$
(the latter via (\ref{eqn:logdiffs})).
Indeed the existence of the extension and the claimed compatibility follow from the
analogous well-known result after base-change to $\CC$.  We therefore conclude that the Kodaira--Spencer
filtration on $\Omega^1_{\widetilde{Y}_U/\CO}$ extends over $\widetilde{Y}_U^\tor$ in the form of a decomposition
$$\Omega^1_{\widetilde{Y}^\tor_U/\CO}(\log \widetilde{Z}_U^\tor) = \bigoplus_{\gp \in S_p} \bigoplus_{i \in \ZZ/f_\gp\ZZ}
\left( \Omega^1_{\widetilde{Y}^\tor_U/\CO}(\log \widetilde{Z}_U^\tor)\right)_{\gp,i},$$
together with an increasing filtration of length $e_\gp$ on each component $(\Omega^1_{\widetilde{Y}^\tor_U/\CO}(\log \widetilde{Z}_U^\tor))_{\gp,i}$,
and isomorphisms
$$\widetilde{\CA}^\tor_{2\be_\theta,-\be_\theta} \stackrel{\sim}{\longrightarrow} \gr^j \left(\Omega^1_{\widetilde{Y}^\tor_U/\CO}(\log \widetilde{Z}_U^\tor)\right)_{\gp,i}.$$
Furthermore for each cusp $\tcalC$ of $\widetilde{Y}_U$ and embeddings $\tau = \tau_{\gp,i}$ and $\theta = \theta_{\gp,i,j}$, the pull-back via 
$\xi = \xi_{\tcalC}$ of $\Fil^j (\Omega^1_{\widetilde{Y}^\tor_U/\CO}(\log \widetilde{Z}_U^\tor))_{\gp,i}$ corresponds
to $t_j(N^{-1}M\otimes\CO)_\tau \otimes_\CO \CO_{\wS} = t_j(\gd^{-1}I^{-1}J \otimes \CO)_\tau \otimes_\CO \CO_{\wS}$  under (\ref{eqn:logdiffs}), and the resulting isomorphism
$$\xi^*\widetilde{\CA}^\tor_{2\be_\theta,-\be_\theta} \stackrel{\sim}{\longrightarrow} \gr^j \left(\xi^*(\Omega^1_{\widetilde{Y}^\tor_U/\CO}(\log \widetilde{Z}_U^\tor))_{\gp,i}\right)
   \cong (\gd^{-1}I^{-1}J \otimes\CO)_\theta \otimes_\CO \CO_{\wS}$$
coincides with the canonical trivialization of (\ref{eqn:Dline}).

We now interpret this in the context of minimal compactifications.  First we note that
the argument of \cite[Prop.~4.9]{rap} yields a Koecher Principle for $\Omega^1_{\widetilde{Y}_U^\tor/\CO}(\log \widetilde{Z}_U^\tor)$,
so that 
$$\pi_*(\Omega^1_{\widetilde{Y}_U^\tor/\CO}(\log \widetilde{Z}_U^\tor))  \longrightarrow
     \widetilde{\iota}_*(\Omega^1_{\widetilde{Y}_U/\CO})$$
is an isomorphism of coherent sheaves on $\widetilde{Y}_U^{\min}$ whose completion at the
cusp $\tcalC$ is identified with $\left(M\otimes \CO[[q^m]]_{m \in N^{-1}M_+\cup\{0\}}\right)^{V_N^2}$
$$ = \left\{ \left.\, \sum_{m\in N^{-1}M_+} c_m \otimes q^m \,\right|\,
   c_{\nu m} = \nu c_m \,\,\forall \nu \in V_N^2, m \in N^{-1}M_+\,\right\}.$$
(Note that $c_0 \in M^{V_N^2} = 0$.)  Furthermore the completion at $\tcalC$
of the canonical derivation $\CO_{\widetilde{Y}_U^{\min}} \to \widetilde{\iota}_* \Omega^1_{\widetilde{Y}_U/\CO}$
is given by $\sum r_mq^m \mapsto \sum m \otimes r_m q^m$.

We see also from the description of the extension of the Kodaira--Spencer filtration to
$\Omega^1_{\widetilde{Y}_U^\tor/\CO}(\log \widetilde{Z}_U^\tor)$ in terms of $q$-expansions that
$$\pi_*(\Fil^j(\Omega^1_{\widetilde{Y}_U^\tor/\CO}(\log \widetilde{Z}_U^\tor))_{\gp,i})
    ={\widetilde{\iota}}_* (\Fil^j(\Omega^1_{\widetilde{Y}_U/\CO})_{\gp,i}),$$
with completion at $\tcalC$ given by
$\left(t_j(M\otimes \CO)_\tau \otimes_\CO \CO[[q^m]]_{m \in N^{-1}M_+\cup\{0\}}\right)^{V_N^2}$
$$ = \left\{ \left.\, \sum_{m\in N^{-1}M_+} c_m \otimes q^m \,\right|\,
   c_{\nu m} = (\nu \otimes 1)c_m\,\, \forall \nu \in V_N^2, m \in N^{-1}M_+\,\right\}.$$
Furthermore since $V_N^2$ acts
freely on $N^{-1}M_+$ (and $M_\theta^{V_N^2} = 0$), the morphisms
$$\widetilde{\iota}_* (\Fil^j(\Omega^1_{\widetilde{Y}_U/\CO})_{\gp,i})  \to \widetilde{\iota}_* \widetilde{\CA}_{2\be_\theta,-\be_\theta}$$
are surjective on completions at cusps, and hence surjective.  It follows that the graded pieces of
(the obvious extension of) the Kodaira--Spencer filtration
on $\widetilde{\iota}_*\Omega^1_{\widetilde{Y}_U/\CO}$ are canonically isomorphic to
$\widetilde{\iota}_* \widetilde{\CA}_{2\be_\theta,-\be_\theta}$.  

The constructions above are compatible with the natural actions of $\CO_{F,(p),+}^\times$
on $\widetilde{Y}_U^\tor$, $\widetilde{Y}_U^{\min}$ and $M$, so that the resulting descriptions carry
over to $Y_U^{\min}$.  More precisely, $\iota_*(\Omega^1_{Y_U/\CO})$
is a coherent sheaf on ${Y}_U^{\min}$ whose completion at the
cusp $\calC$ is identified with
$\left(M\otimes \CO[[q^m]]_{m \in N^{-1}M_+\cup\{0\}}\right)^{V_{N,+}}$
$$= \left\{ \left.\, \sum_{m\in N^{-1}M_+} c_m \otimes q^m \,\right|\,
   c_{\nu m} = \nu c_m\,\, \forall \nu \in V_{N,+}, m \in N^{-1}M_+\,\right\},$$
in terms of which the canonical derivation is $\sum r_mq^m \mapsto \sum m \otimes r_m q^m$.
Furthermore the completion at $\calC$ of 
$\Fil^j(\iota_*(\Omega^1_{Y_U/\CO})_{\gp,i}) : = \iota_* (\Fil^j(\Omega^1_{Y_U/\CO})_{\gp,i})$ is identified
with 
$\left(t_j(M\otimes \CO)_\tau \otimes_\CO \CO[[q^m]]_{m \in N^{-1}M_+\cup\{0\}}\right)^{V_{N,+}}$,
and the natural maps 
\begin{equation} \label{eqn:KScusps} \gr^j(\iota_* (\Omega^1_{Y_U/\CO})_{\gp,i}) 
\hookrightarrow \iota_* (\gr^j(\Omega^1_{Y_U/\CO})_{\gp,i}) \stackrel{\sim}{\longrightarrow}
\iota_*\CA_{2\be_\theta,-\be_\theta} \end{equation}
are isomorphisms whose completions at the cusps are induced by the surjections
$t_j(M\otimes \CO)_\tau \to M_\theta$.

Suppose now that $U'$ is an arbitrary sufficiently small open compact subgroup
of $\GL_2(\AA_{F,\f})$ of level prime to $p$, and choose $N$ so that $U(N) \subset U'$.
The constructions above are then also compatible with the natural actions of $U'$, so we
arrive at similar conclusions with minor modifications to the descriptions of completions
that result from taking invariants under $\Gamma_{\calC,U'}$.  We omit the details, but
we remark that letting $L$ (resp.~$V$) denote the kernel (resp.~image) of the homomorphism
$$\begin{array}{ccc} \Gamma_{\calC,U'}/(\CO_F^\times\cap U') & \to & \CO_{F,+}^\times \\
 \left(\begin{array}{cc}\alpha&\beta\\0&\delta\end{array}\right)\cdot(\CO_F^\times \cap U') & \mapsto & \alpha^{-1}\delta,\end{array}$$
the coefficients $r_m$ of $q$-expansions in $\CO_{Y_{U'}^{\min},\calC'}^\wedge$ are indexed by
$m \in (L^*)_+\cup\{0\}$ (where $L$ is identified with a finite index subgroup of $M$).  
Since $V$ acts freely on $L^*$ (twisting $q$-expansion coefficients
by a possibly non-trivial cocycle valued in $L\otimes \mu_N(\CO)$), we still obtain the isomorphism
of (\ref{eqn:KScusps}) with $U$ replaced by $U'$, identifying the graded pieces
of the Kodaira--Spencer filtration on $\iota_*\Omega^1_{Y_{U'}/\CO}$ with the sheaves
$\iota_*\CA_{2\be_\theta,-\be_\theta}$.

The description of the extension of the Kodaira--Spencer filtration over compactifications also
applies after base-change to an arbitrary $\CO$-algebra $R$, with one significant difference.
If $R$ is not flat over $\CO$, then the modules $M\otimes R$ (and their subquotients) may have
invariants under the action of the unit groups $V_{N,+}$ (or more generally the isotropy groups
$\Gamma_{\calC,U'}$), so that $q$-expansions of meromorphic differentials on $Y_{U,R}$ (and forms of
weight $(2\be_\theta,-\be_\theta)$) may have non-zero constant terms, and the morphism analogous to
(\ref{eqn:KScusps}) may fail to be an isomorphism.  (Note that in this case the relevant base-change
morphisms $(\iota_*\CF)_R \to \iota_{R,*}(\CF_R)$ fail to be surjective at the cusps.)

We can however simplify matters by placing ourselves in the situation when this fails in the extreme.  
Suppose then that $p^n R = 0$ for some $n > 0$, and that $N$ is sufficiently large that $\nu \equiv 1 \bmod p^n\CO_F$
for all $\nu \in V_{N,+}$.  Arguing exactly as above, we find that $\iota_{R,*}(\Omega^1_{Y_{U,R}/R})$
is now a vector bundle over ${Y}_{U,R}^{\min}$ whose completion at $\calC$ is identified with
$\left(M\otimes R[[q^m]]_{m \in N^{-1}M_+\cup\{0\}}\right)^{V_{N,+}} =
  M \otimes \CO_{Y^{\min}_{U,R},\calC}^\wedge$
  $$ =  \left\{ \left.\, \sum_{m\in N^{-1}M_+ \cup\{0\}} c_m \otimes q^m \,\right|\,
   c_{\nu m} = \nu c_m\,\, \forall \nu \in V_{N,+}, m \in N^{-1}M_+\,\right\},$$
with the canonical derivation given by $\sum r_mq^m \mapsto \sum m \otimes r_m q^m$.
Furthermore each
$$\Fil^j(\iota_{R,*}(\Omega^1_{Y_{U,R}/R})_{\gp,i}) = \iota_{R,*} (\Fil^j(\Omega^1_{Y_{U,R}/R})_{\gp,i})$$
is a sub-bundle whose completion at $\calC$ is identified
with 
$t_j(M\otimes \CO)_\tau \otimes_\CO \CO_{Y^{\min}_{U,R},\calC}^\wedge$
and the natural maps
$$ \gr^j(\iota_{R,*}(\Omega^1_{Y_{U,R}/R})_{\gp,i}) 
\hookrightarrow \iota_{R,*} (\gr^j(\Omega^1_{Y_{U,R}/R})_{\gp,i}) \stackrel{\sim}{\longrightarrow}
\iota_{R,*}\CA_{2\be_\theta,-\be_\theta,R} $$
are isomorphisms of line bundles over $Y_{U,R}^{\min}$
whose completions at the cusps are induced by the surjections
$t_j(M\otimes \CO)_\tau \to M_\theta$.  We remark also that
this carries over with $U$ replaced by arbitrary $U'$,
provided $U'$ is sufficiently small that (in additional to the usual hypotheses)
$\alpha \equiv \delta \mod p^n\CO_F$ for all 
$\left(\begin{array}{cc}\alpha&\beta\\0&\delta\end{array}\right)\in \Gamma_{\calC,U'}$
and cusps $\calC'$ of $Y_{U'}$ (the condition being independent of the choice
of $N$ and $\calC$ in the definition of $\Gamma_{\calC,U'}$).
 
\section{Operators on $q$-expansions}  \label{q2}
\subsection{Partial Hasse invariants}  \label{sec:Hasse2}
We next describe the effect of the various weight-shifting operators on $q$-expansions, beginning
with the simplest case of (multiplication by) partial Hasse invariants.
We will now only be working in the setting of $R = \FF$, and we will use
$\overline{\cdot}$ to denote base-changes from $\CO$ to $\FF$.
Since the formation of $q$-expansions is compatible in the obvious sense with pull-back under
the projections $\overline{Y}_U \to \overline{Y}_{U'}$,  it will
suffice to consider the case $U = U(N)$.  

Recall that in \S\ref{sec:Hasse} we defined the partial Hasse invariants as certain elements
 $$ H_\theta  \in M_{\bh_\theta,\bf{0}} (U; \FF) = H^0(\overline{Y}_{U} , \ol{\CA}_{\bh_\theta,\bf{0}}),$$
where $\bh_\theta: = n_\theta \be_{\sigma^{-1}\theta} - \be_\theta$, with $n_\theta = p$ if $j=1$ and $n_\theta = 1$ if $j > 1$.
 In particular if $j > 1$, then $H_\theta$ is defined by the morphism $u:{\CL}_{\tau,j} \to {\CL}_{\tau,j-1}$
 induced by $\varpi_\gp$ on the universal abelian variety over $\widetilde{Y}_{U,\FF}$, which evidently
 extends to the endomorphism $\varpi_\gp$ of $A^\tor$ over $S:= \widetilde{Y}_{U,\FF}^\tor$. 
 Since its pull-back via $\ol{\xi}$ is defined
 by $\varpi_\gp$ on $\ol{T}_{I,J}$, the resulting morphism of line bundles
 $$\ol{\xi}^*{{\CL}}_{\tau,j}^\tor   \longrightarrow \ol{\xi}^*{{\CL}}_{\tau,j-1}^\tor$$
 is compatible with their canonical trivializations, and more precisely with the morphism
 $$(I^{-1})_\theta \otimes_{\CO} \FF   \to (I^{-1})_{\sigma^{-1}\theta} \otimes_{\CO} \FF $$ 
  induced by $u = \varpi \otimes 1$ on $(I^{-1} \otimes \CO)_\tau$.    It follows that
  $H_\theta$ has constant $q$-expansion, where the constant corresponds to the
  basis element $\varpi_\gp$ under the identification
$$\ol{D}_{\bh_\theta,\bf{0}} \cong \gp \otimes_{\CO_F,\theta}  \FF.$$
 provided by (\ref{eqn:modp}).  
 
 For $j=1$ we use also that the morphism of line bundles
 $$\ol{\xi}^*{{\CL}}_{\tau,1}^\tor   \longrightarrow \ol{\xi}^*\phi_S^*({{\CL}}_{\phi^{-1}\circ\tau,e_\gp}^\tor)$$
induced by the Verschiebung $\phi_S^*\ol{T}_{I,J} \to \ol{T}_{I,J}$
is compatible via the canonical trivializations with the canonical isomorphism
$$(I^{-1} \otimes \FF)_\tau   \stackrel{\sim}{\longrightarrow} \phi_S^* (I^{-1} \otimes \FF)_{\phi^{-1}\circ\tau}.$$
So in this case we again find that $H_\theta$ has constant $q$-expansion, the constant now
corresponding to the basis element $\varpi_\gp^{1-e_\gp}$ under the identification
$$\ol{D}_{\bh_\theta,\bf{0}} \cong \gp^{1-e_\gp} \otimes_{\CO_F,\theta}  \FF$$
given by  (\ref{eqn:modp}).

The $q$-expansions of the canonical sections $G_\theta  \in M_{\bf{0},\bh_\theta} (U; \FF)$ may
be described similarly.  Indeed for $j > 1$, the composites
$$\begin{array}{cccccc}
&{\CL}_{\tau,j} & \longrightarrow &{\CP}_{\tau,j} & \stackrel{u}{\longrightarrow} & {\CL}_{\tau,j-1} \\
\mbox{and}&  {\CM}_{\tau,j-1}  &\longrightarrow &{\CP}_{\tau,j} & \stackrel{u}{\longrightarrow} & {\CM}_{\tau,j} \end{array}$$
on $\widetilde{Y}^\ord_{U,\FF}$ are isomorphisms whose tensor product defines
${\CN}_\theta \cong {\CN}_{\sigma^{-1}\theta}$.  Its unique extension to $\widetilde{Y}^\tord_{U,\FF}$,
and hence to $\widetilde{Y}^\tor_{U,\FF}$,
is therefore the isomorphism whose pull-back via $\ol{\xi}$ is the tensor product of the
isomorphisms defined on canonical trivializations by
$(I^{-1})_\theta \otimes_{\CO} \FF   \stackrel{u}{\to} (I^{-1})_{\sigma^{-1}\theta} \otimes_{\CO} \FF $
and by the identity on $(\gd J^{-1} \otimes \FF)_\tau/u (\gd J^{-1} \otimes \FF)_\tau$, and hence corresponds to
$$(\gd(IJ)^{-1})_\theta \otimes_{\CO} \FF   \stackrel{u}{\to} (\gd(IJ)^{-1})_{\sigma^{-1}\theta} \otimes_{\CO} \FF.$$
Therefore $G_\theta$ has constant $q$-expansion, with constant corresponding to the
  basis element $\varpi_\gp$ under the identification
$$\ol{D}_{\bf{0},\bh_\theta} \cong \gp \otimes_{\CO_F,\theta}  \FF$$
 provided by the analogue of (\ref{eqn:modp}) for $(\gd(IJ) ^{-1})_\theta$.

Similarly if $j=1$, we find that ${\CN}_{\tau,1} \cong \phi_S^*({\CN}_{\phi^{-1}\circ\tau,e_\gp})$
over $S = \widetilde{Y}^\ord_{U,\FF}$ is the tensor product of the isomorphism
$${\CL}_{\tau,1} \to \phi_S^* ({\CL}_{\phi^{-1}\circ\tau,e_\gp})$$
defining the Hasse invariant and the isomorphism
$$\phi_S^*({\CM}_{\phi^{-1}\circ\tau,e_\gp}) \to {\CM}_{\tau,1}$$
induced by $\Frob_A$. The extensions to $\widetilde{Y}^\tord_{U,\FF}$ are again compatible
with the canonical trivializations, now corresponding to maps whose tensor product is
the inverse of the isomorphism
$$((\gd(IJ)^{-1})_{\sigma^{-1}\theta} \otimes_{\CO} \FF)^{\otimes p}   \stackrel{\sim}{\longrightarrow} 
(\gd(IJ)^{-1})_{\theta} \otimes_{\CO} \FF$$
induced by $u^{e_\gp-1}$.  So again $G_\theta$ has constant $q$-expansion, 
with constant corresponding to $\varpi_\gp^{1-e_\gp}$ under the identification
$$\ol{D}_{\bf{0},\bh_\theta} \cong \gp^{1-e_\gp} \otimes_{\CO_F,\theta}  \FF$$
given by (\ref{eqn:modp}) for $(\gd(IJ) ^{-1})_\theta$.

\subsection{Partial $\Theta$-operators}  \label{sec:theta2}
We now compute the effect of $\Theta$-operators on $q$-expansions exactly as in \cite{DS}.
Recall from \S\ref{sec:theta} that for each $\tau_0  = \tau_{\gp,i} \in \Sigma_0$, the associated partial $\Theta$-operator is a map
$$\Theta_{\tau_0}: M_{\bk,\bl}(U;\F) \to M_{\bk',\bl'}(U;\F)$$
where  $\bk' = \bk + n_{\theta_0} \be_{\sigma^{-1}\theta_0} + \be_{\theta_0}$, $\bl' = \bl + \be_{\theta_0}$
and $\theta_0 = \theta_{\gp,i,e_\gp}$.  It is defined for all sufficiently small $U$ of level prime to $p$,
and is Hecke-equivariant.  In particular it is compatible with restriction for $U \subset U'$, so we may 
assume $U = U(N)$ for some $N$ sufficiently large that $\nu \equiv 1 \mod p\CO_F$ for all $\nu \in V_{N,+}$.

Recall from the proof of Theorem~\ref{thm:theta} that $\Theta_{\tau_0}$ is defined by
a morphism $\CA_{\bk,\bl,\FF} \to \CA_{\bk',\bl',\FF}$ given locally
on sections by formula (\ref{eqn:altdeftheta}).  Our assumptions on $U$ imply that
$$(\iota_{\FF,*}\CA_{\bk,\bl,\FF})_{\calC}^\wedge = 
D_{\bk,\bl,\FF} \otimes_{\FF}\cdot \CO_{Y_{U,\FF}^{\min},\calC}^\wedge$$
is free of rank one over $\CO_{Y_{U,\FF}^{\min},\calC}^\wedge$ for all weights $\bk,\bl$
and cusps $\calC$, so the completion at $\calC$ of $\iota_{\FF,*}\Theta_{\tau_0}$ is the map
$$(\iota_{\FF,*}\CA_{\bk,\bl,\FF})_{\calC}^\wedge  \to (\iota_{\FF,*}\CA_{\bk',\bl',\FF})_{\calC}^\wedge$$
defined by (\ref{eqn:altdeftheta}), where $y_\theta$
is any basis for $(\iota_{\FF,*}{\CL}_\theta)_{\calC}^\wedge$ and
$y^{\bk} = \prod_\theta y_\theta^{k_\theta}$.  In particular we may choose 
$y_\theta = b_\theta \otimes 1$ where $b_\theta$ is a basis for 
$D_{\be_\theta,\bf{0},\FF} = (I^{-1})_\theta \otimes_\CO\FF$.  The fact that $H_\theta$
has (non-zero) constant $q$-expansion at $\calC$ then means the same holds for the 
element $r_\theta \in \CO_{Y_{U,\FF}^{\min},\calC}^\wedge$, i.e., $r_\theta \in \FF^\times$.
We can even select the bases $b_\theta$ so that $r_\theta = 1$ for all $\theta$
by choosing any basis $b_\gp$ for $(I_\gp \otimes \FF)_{\gp,0}$ over $\FF[u]/u^{e_\gp}$
for each $\gp$, letting $b_{\gp,i}$ correspond to $(\phi^i)^*(b_\gp)$ under the canonical
isomorphism $(I_\gp \otimes \FF)_{\gp,i} \cong (\phi^i)^*(I_\gp \otimes \FF)_{\gp,0}$,
and defining $b_{\gp,i,j}$ to be the image of $u^{e_\gp-j}b_{\gp,i}$ in 
$(I^{-1})_\theta \otimes_\CO\FF$

Recall also that $g^{\bl} = \prod_\theta g_\theta^{l_\theta}$ in (\ref{eqn:altdeftheta}),
where each $g_\theta$ is a trivialization of ${\CN}_\theta$ over $Y_{U,\FF}$.
Therefore $g_\theta$ trivializes $\iota_{\FF,*}{\CN}_\theta$ over $Y^{\min}_{U,\FF}$,
from which it follows that $g_\theta =  c_\theta \otimes 1$ for some basis $c_\theta$
of $D_{\bf{0},\be_\theta,\FF} = (\gd(IJ)^{-1})_\theta \otimes_\CO\FF$. 
The formula (\ref{eqn:altdeftheta}) therefore takes the form
$$\Theta_{\tau_0} (b^{\bk} c^{\bl} \otimes \varphi_f ) =  H_{\theta_0} b^{\bk}c^{\bl}  \otimes KS_{\tau_0} (d\varphi_f) =
b^{\bk+\bh_{\theta_0}} c^{\bl} \otimes KS_{\tau_0} (d\varphi_f),$$
for $\varphi_f \in \CO_{Y_{U,\FF}^{\min},\calC}^\wedge$.
Finally the descriptions in \S\ref{sec:KS2}
of the canonical derivation, the Kodaira--Spencer filtration and the isomorphism
(\ref{eqn:KScusps}) in terms of $q$-expansions yield the formula
\begin{equation}  \label{eqn:Thetaq} \Theta_{\tau_0}\left(  \sum_{m\in N^{-1}M_+\cup\{0\}}(b^{\bk} c^{\bl} \otimes  r_m) q^m \right) = 
 \sum_{m\in N^{-1}M_+}  (b^{\bk+\bh_{\theta_0}} c^{\bl} \ol{\tau}_0(m)  \otimes  r_m) q^m \end{equation}
where $\ol{\tau}_0$ is the canonical projection
$$\begin{array}{rcl}N^{-1}M  &\to &(N^{-1}M \otimes \FF)_{\tau_0}  = (M\otimes \FF)_{\tau_0} \\
& \to &M\otimes_{\CO_F,\ol{\theta}_0} \FF = (\gd^{-1} I^{-1} J) \otimes_{\CO_F,\ol{\theta}_0} \FF 
= D_{2\be_{\theta_0},-\be_{\theta_0},\FF}\end{array}$$
(writing $\ol{\theta}_0$ for the composite of $\CO_F \stackrel{\theta_0}{\to} \CO \to \FF$).
As noted above, it follows that (\ref{eqn:Thetaq}) holds with $U = U(N)$ replaced by
any sufficiently small open compact $U'$ of level prime to $p$ and $\calC$ replaced by any cusp
of $Y_{U'}$.  In this case the $q$-expansions
are necessarily invariant under the natural action of $\Gamma_{\calC,U'}$ (whose compatibility
with (\ref{eqn:Thetaq}) is a consequence of the construction, but is straightforward to check directly).

We see immediately from (\ref{eqn:Thetaq}) that the operators $\Theta_\tau$ for varying $\tau$ commute.
We see also that 
$$\Theta_{\tau_0}^p(f) = \Theta_{\tau_1}(f) H_{\theta_0}^p  H_{\theta_1} G_{\theta_1}^{-1}
\prod_{j=1}^{e_\gp-1} \left(H_{\sigma^j\theta_0} G^{-1}_{\sigma^j\theta_0}\right),$$
where $\tau_1 = \tau_0\circ\phi = \tau_{\gp,i+1}$ and $\theta_1 = \sigma^{e_\gp}\theta_0 = \theta_{\gp,i+1,e_\gp}$.
Indeed this follows from (\ref{eqn:Thetaq}) together with the fact that the $q$-expansions of 
$$H_{\tau_1} = \prod_{j=1}^{e_\gp} H_{\sigma^j\theta_0} \quad\mbox{and}\quad
    G_{\tau_1} = \prod_{j=1}^{e_\gp} G_{\sigma^j\theta_0}$$
are constants given by the canonical isomorphisms 
$((I^{-1})_{\theta_0} \otimes_{\CO} \FF)^{\otimes p} \cong (I^{-1})_{\theta_1} \otimes_{\CO} \FF$
and $((\gd(IJ)^{-1})_{\theta_0} \otimes_{\CO} \FF)^{\otimes p} \cong (\gd(IJ)^{-1})_{\theta_1} \otimes_{\CO} \FF$,
so we get a commutative diagram
$$\xymatrix{  &  D_{2\be_{\theta_1},-\be_{\theta_1},\FF} \ar[dd]\\ 
N^{-1}M  \ar[ur]^{\ol{\tau}_1} \ar[dr]_{\ol{\tau_0}^{\otimes p}} & \\
& D_{2\be_{\theta_0},-\be_{\theta_0},\FF}^{\otimes p}}$$
where the downward arrow is multiplication by the (constant) $q$-expansion of $H_{\tau_1}^2G_{\tau_1}^{-1}$.

\subsection{Partial Frobenius operators}  \label{sec:V2}
Finally we compute the effect on $q$-expansions of the partial Frobenius operators $V_\gp$ defined in \S\ref{sec:V}.
We must first extend the partial Frobenius endomorphisms 
$\widetilde{\Phi}_\gp$ (and $\Phi_\gp$) defined in \S\ref{sec:Frobenius} to compactifications.
To that end let $h_\gp$ denote
the matrix $\smallmat100\delta$, viewed as an element of
$\GL_2(\A_{F,\f}^{(p)})$, where $\delta \in F^\times$ is as in the definition fo $\widetilde{\Phi}_\gp$.
We let $\widetilde{\Phi}_{\gp}^\infty:\widetilde{Y}^\infty_U \to \widetilde{Y}^\infty_U$ be the permutation defined on double
cosets by $B_1(\CO_{F,(p)})gU^p \mapsto B_1(\CO_{F,(p)})h_\gp^{-1}gU^p$, and similarly
let $\Phi_{\gp}^\infty$ denote the induced permutation of $Y^\infty_U$.  Then $\widetilde{\Phi}^\infty_\gp$ translates to the
map on corresponding data sending $(H,I,\lambda,[\eta])$ to $(H',I',\lambda',[\eta'])$,
where
\begin{itemize}
\item $I' = \gp^{-1}\otimes_{\CO_F}I = \gp^{-1}I$,
\item $H'$ is the push-out of $H$ with respect to $I \to I'$;
\item $\lambda' = \delta\otimes\lambda$ (identifying $\wedge^2_{\CO_F}H' = \gp^{-1}\otimes_{\CO_F} \wedge^2_{\CO_F}H$),
\item and $\eta' = \eta$ (identifying $\widehat{\CO}_F^{(p)} \otimes_{\CO_F} H'  = \widehat{\CO}_F^{(p)} \otimes_{\CO_F} H$).
\end{itemize}

Suppose now that $U = U(N)$ for some sufficiently large $N$.  One then
checks that the morphism $\widetilde{\Phi}_\gp: \widetilde{Y}_{U,\FF} \to \widetilde{Y}_{U,\FF}$ extends to
a morphism $\widetilde{\Phi}_\gp^\tor: \widetilde{Y}_{U,\FF}^{\tor '}\to \widetilde{Y}_{U,\FF}^\tor$ (where the
$'$ indicates the choice of cone decomposition need not be the same) under which the
component corresponding to a cusp $\tcalC$ represented by $(H,I,\lambda,[\eta])$ is sent to the one
corresponding to $\tcalC' = \widetilde{\Phi}_\gp^\infty(\tcalC)$ represented by  $(H',I',\lambda',[\eta'])$, and the resulting
map on completions pulls back to a morphism $\wS'_{\tcalC,\FF} \to \wS_{\tcalC',\FF}$
whose effect on global sections corresponds under the isomorphisms of (\ref{eqn:bigcompletion})
to the homomorphism induced by the canonical inclusion 
$$M' := \gd^{-1}(I')^{-1}J = \gp\gd^{-1}I^{-1}J \hookrightarrow \gd^{-1}I^{-1}J.$$
Furthermore the pull-back of $A'^\tor = (\widetilde{\Phi}^\tor_\gp)^*A^\tor$ to $\wS'_{\tcalC,\FF}$ is
identified with the Tate semi-abelian variety $T_{\gp^{-1}I,J}$, and the isomorphisms 
defined in \S\ref{sec:V} relating the line bundles $\widetilde{\Phi}_\gp^*\CL_\theta$
to $\CL_\theta$ or $\CL_{\sigma^{-1}\theta}$  extend to isomorphisms
$$\mbox{$(\widetilde{\Phi}^\tor_\gp)^* \CL^\tor_\theta \cong \CL^{\tor'}_\theta$ for $\theta\not\in \Sigma_\gp$,
and $(\widetilde{\Phi}^\tor_\gp)^* \CL_\theta^\tor \cong (\CL^{\tor'}_{\sigma^{-1}\theta})^{\otimes n_\theta}$ for $\theta\in \Sigma_\gp$}$$
over $\widetilde{Y}_{U,\FF}^{\tor'}$ whose pull-backs are compatible via their canonical trivializations with isomorphisms induced by
the canonical $(\CO_F\otimes \FF)_\tau$-equivariant maps
$$(pI^{-1} \otimes \FF) \longrightarrow (\gp I^{-1} \otimes \FF)_\tau \longrightarrow (I^{-1} \otimes \FF)_\tau.$$
More precisely if $\tau \not\in \Sigma_{\gp,0}$, then the second map is an isomorphism identifying
$(\gp I^{-1})_\theta \otimes_\CO \FF$ with $(I^{-1})_\theta \otimes_\CO \FF $, and if
$\tau = \tau_{\gp,i}$ and $\theta = \theta_{\gp,i,j}$, then this map also induces the desired isomorphisms
$$\begin{array}{rcl}
(\gp I^{-1})_\theta\otimes_\CO \FF  &=& u^{e_\gp - j}(\gp I^{-1} \otimes \FF)_\tau \otimes_{\FF[u]/u^{e_\gp}} \FF\\
  & \stackrel{\sim}{\longrightarrow}& u^{e_\gp - j+1}(I^{-1} \otimes \FF)_\tau \otimes_{\FF[u]/u^{e_\gp}} \FF
   = (I^{-1})_{\sigma^{-1}\theta}\otimes_\CO \FF\end{array}$$
for $j=2,\ldots,e_\gp$.  On the other hand if $\tau = \tau_{\gp,i}$ and $\theta = \theta_{\gp,i,1}$, then
the first map induces an isomorphism
$$(pI^{-1} \otimes  \FF)_\tau \otimes_{\FF[u]/u^{e_\gp}} \FF  \stackrel{\sim}{\longrightarrow} 
u^{e_\gp - j}(\gp I^{-1} \otimes \FF)_\tau \otimes_{\FF[u]/u^{e_\gp}} \FF = (\gp I^{-1})_\theta\otimes_\CO \FF$$
whose composite with the ones induced by
$$\phi^*(I^{-1}\otimes \FF)_{\phi^{-1}\circ\tau}\stackrel{\sim}{\longrightarrow}
(I^{-1}\otimes \FF)_\tau  \stackrel{p\otimes 1}{\longrightarrow}
(pI^{-1}\otimes \FF)_\tau$$
yields the desired isomorphism
$$((I^{-1})_{\sigma^{-1}\theta} \otimes_\CO \FF)^{\otimes p} \stackrel{\sim}{\longrightarrow} 
(I^{-1})_{\theta} \otimes_\CO \FF.$$
The relations between the line bundles $\widetilde{\Phi}_\gp^*\CN_\theta$
and $\CN_\theta$ or $\CN_{\sigma^{-1}\theta}$ extend similarly over $\widetilde{Y}_{U,\FF}^{\tor'}$, so for
$\bk'',\bl''$ as in (\ref{eqn:Phiweight}) we obtain
isomorphisms
$$(\widetilde{\Phi}^\tor_\gp)^* \widetilde{\CA}_{\bk,\bl,\FF}^\tor \cong \widetilde{\CA}_{\bk'',\bl'',\FF}^{\tor'}$$
whose pull-backs to $\wS'_{\tcalC,\FF}$ are compatible via their canonical trivializations with the isomorphisms
\begin{equation} \label{eqn:Phibres} \overline{D}'_{\bk,\bl}   \cong \overline{D}_{\bk'',\bl''}  \end{equation}
obtained as the tensor products of the ones just defined
(where $\overline{D}'_{\bk,\bl}$ is associated to the data for the cusp $\tcalC'$, and $\overline{D}_{\bk'',\bl''}$ 
to the data for $\tcalC$).

It follows from the above description of $\widetilde{\Phi}_\gp^\tor$ that $\widetilde{\Phi}_\gp$ extends to the morphism 
$\widetilde{\Phi}_\gp^{\min}:\widetilde{Y}_{U,\FF}^{\min}  \to \widetilde{Y}_{U,\FF}^{\min}$ restricting to $\widetilde{\Phi}^\infty_{\gp}$
on the set of cusps, with the induced maps on completed local rings\footnote{Note that we have implicitly chosen
different representatives $(H,I,\lambda,[\eta])$ for each cusp $\tcalC$ according to whether $\widetilde{Y}_{U,\FF}$
is viewed as the source or target of $\widetilde{\Phi}_\gp$, but the rings and modules arising from the two different
descriptions of completions at $\tcalC$ are canonically isomorphic.} being the restriction to $V_N^2$-invariants of
the canonical inclusion
\begin{equation} \label{eqn:Phat} \FF[[q^m]]_{m \in (N^{-1}M')_+\cup\{0\}}  \longrightarrow \FF[[q^m]]_{m \in (N^{-1}M)_+\cup \{0\}},\end{equation}
where $M = \gd^{-1}I^{-1}J$ and $M' = \gp M$.  Furthermore the commutativity of the diagram
$$\xymatrix{
 ({\widetilde{\iota}}_{\FF,*}{\widetilde{\CA}}_{\bk,\bl,\FF})^\wedge_{\tcalC'}
  \ar@{^{(}->}[d]
&  \Gamma( (\widetilde{Y}_{U,\FF}^\tor)_{\tcalC'}^\wedge,  ({\widetilde{\CA}}^\tor_{\bk,\bl,\FF})_{\tcalC'}^\wedge) 
    \ar@{^{(}->}[d]  \ar[l]_-{\sim}  \ar[r]^-{\sim}
& ( \overline{D}'_{\bk,\bl}  \otimes_{\FF} \Gamma(\wS_{\tcalC',\FF}, \CO_{\wS_{\tcalC',\FF}}) )^{V_N^2} \ar@{^{(}->}[d]
\\({\widetilde{\iota}}_{\FF,*}{\widetilde{\Phi}_\gp^*\widetilde{\CA}}_{\bk,\bl,\FF})^\wedge_\tcalC  \ar[d]^{\wr}
 &  \Gamma( (\widetilde{Y}_{U,\FF}^{\tor'})_{\tcalC}^\wedge,  (\widetilde{\Phi}_\gp^{\tor,*}({\widetilde{\CA}}^\tor_{\bk,\bl,\FF}))_{\tcalC}^\wedge)  
  \ar[d]^{\wr}  \ar[l]_-{\sim} \ar[r]^-{\sim}
 &   ( \overline{D}'_{\bk,\bl}  \otimes_{\FF} \Gamma(\wS'_{\tcalC,\FF}, \CO_{\wS'_{\tcalC,\FF}}) )^{V_N^2}   \ar[d]^{\wr}
\\   ({\widetilde{\iota}}_{\FF,*}{\widetilde{\CA}}_{\bk'',\bl'',\FF})^\wedge_\tcalC  
 &  \Gamma( (\widetilde{Y}_{U,\FF}^{\tor'})_{\tcalC}^\wedge,  ({\widetilde{\CA}}^{\tor'}_{\bk'',\bl'',\FF})_{\tcalC}^\wedge)
   \ar[l]_-{\sim} \ar[r]^-{\sim}   
 & ( \overline{D}_{\bk'',\bl''}  \otimes_{\FF} \Gamma(\wS'_{\tcalC,\FF}, \CO_{\wS'_{\tcalC,\FF}}) )^{V_N^2}
   }$$
(where the top vertical arrows are defined by pulling back via $\widetilde{\Phi}_\gp^{\min}$, $\widetilde{\Phi}_\gp^{\tor}$
and the map $\wS'_{\tcalC,\FF} \to \wS_{\tcalC',\FF}$) shows that the resulting map on $q$-expansions
is the restriction to $V_N^2$-invariants of the map
$$ \overline{D}'_{\bk,\bl}   \otimes_\FF \FF[[q^m]]_{m \in (N^{-1}M')_+\cup\{0\}}  \longrightarrow
\overline{D}_{\bk'',\bl''}   \otimes_{\FF} \FF[[q^m]]_{m \in (N^{-1}M)_+\cup \{0\}},$$
obtained as the tensor product of (\ref{eqn:Phibres}) and (\ref{eqn:Phat}).
(Note that the isomorphism (\ref{eqn:Phibres}) is $V_N^2$-equivariant, but we can also choose 
$N$ sufficiently large that the action is trivial.)

The constructions above are all compatible with the natural action of $\CO_{F,(p),+}^\times$, so the morphism
$\widetilde{\Phi}_\gp^{\min}$ induces a morphism $\Phi_{\gp}^{\min} : \ol{Y}_U^{\min} \to \ol{Y}_U^{\min}$ extending
$\Phi_\gp$ by the map $\Phi^\infty_{\gp}$ on cusps, and its effect on completed local rings is given by the
$V_{N,+}$-invariants of (\ref{eqn:Phat}).  Furthermore the map $V_\gp:  M_{\bk,\bl}(U;\FF)  \longrightarrow M_{\bk'',\bl''}(U,\FF)$
is described on $q$-expansions by taking the $V_{N,+}$-invariants of the tensor product of (\ref{eqn:Phibres})
and (\ref{eqn:Phat}).  (Note that $\Phi_\gp^{\min}$ is proper and quasi-finite, hence finite, but not necessarily
flat at the cusps.)

Similarly for any sufficiently small $U'$, we may choose $N$ so $U = U(N) \subset U'$
and take invariants under the natural action of $U'/U$, with which the above constructions are also easily
seen to be compatible.  We thus obtain the description of $V_\gp$ on $q$-expansions (under the
identifications of Proposition~\ref{prop:minimal}) as the resulting map
\begin{equation} \label{eqn:Vonq}
\left( \overline{D}'_{\bk,\bl}   \otimes_\FF \FF[[q^m]]_{m \in (N^{-1}M')_+\cup\{0\}}  \right)^{\Gamma_{\calC',U'}}
\longrightarrow
\left( \overline{D}_{\bk'',\bl''}   \otimes_\FF \FF[[q^m]]_{m \in (N^{-1}M)_+\cup\{0\}}  \right)^{\Gamma_{\calC,U'}}.
\end{equation}
(Note that the maps (\ref{eqn:Phibres}) and (\ref{eqn:Phat}) are in fact $\Gamma_{\calC}$-equivariant, where
$\Gamma_{\calC}$ is defined in (\ref{eqn:Isotropy}) and its action on the target is via the natural inclusion in $\Gamma_{\calC'}$.)

Finally we note that the effect of the operator  $V_\gp^0:  M_{\bk,\bl}(U;\FF)  \longrightarrow M_{\bk'',\bl}(U,\FF)$
on $q$-expansions has the same description, but with (\ref{eqn:Phibres}) replaced by its composite with the
isomorphism 
$$\overline{D}_{\bk'',\bl''}  \cong \overline{D}_{\bf{0},\bl-\bl''} \otimes_{\FF}  \overline{D}_{\bk'',\bl''} 
= \overline{D}_{\bk'',\bl}$$
given by choosing the basis element of $\overline{D}_{\bf{0},\bl-\bl''}$ to be the (constant) $q$-expansion
of $\prod_{\theta\in\Sigma} G_\theta^{-l_\theta}$.

\section{The kernel of $\Theta$}  \label{kernel}
\subsection{Determination of the kernel}  \label{sec:ker}
In this section we analyze the kernel of the partial $\Theta$-operator 
$\Theta_{\tau}: M_{\bk,\bl}(U;\F) \to M_{\bk',\bl'}(U;\F)$
for $\tau \in \Sigma_{\gp,0}$ and relate it to the 
image of a partial Frobenius operator.

We allow $U$ to be any sufficiently small open compact subgroup of $\GL_2(\A_{F,\f})$ of
level prime to $p$, and $(\bk,\bl)$ any weight such that $\chi_{\bk+2\bl,\FF}$ is trivial on
$U \cap \CO_F^\times$.   First note that by (\ref{eqn:Thetaq}) and the $q$-expansion Principle
(Proposition~\ref{prop:qexp}), the kernel of
$\Theta_\tau$ consists precisely of those $f$ whose $q$-expansions 
$\sum_{m\in N^{-1}M_+\cup\{0\}}(b^{\bk} c^{\bl} \otimes  r_m) q^m$
at all cusps,  or even a cusp on each connected component, satisfy
$r_m = 0$ for all $m \not\in \gp N^{-1}M$.  (Recall that $M = \gd^{-1}I^{-1}J$ if the cusp $\calC$
if represented by $(H,I,[\lambda],[\eta])$, and note that the condition depends only on the isomorphism
class of $(H,I,[\lambda],[\eta])$ and in particular is independent of $N$ prime to $p$ such that $U(N) \subset U$.)
Note that the condition is the same for all $\tau \in \Sigma_{\gp,0}$, so that
$$\ker(\Theta_\tau) = \ker(\Theta_{\tau'}) \quad\mbox{for all $\tau,\tau' \in \Sigma_{\gp,0}$}.$$
Note also that the condition is invariant under multiplication
by the Hasse invariants $H_\theta$ (and of course the forms $G_\theta$) for all $\theta$, so
that $\Theta_\tau(f) = 0$ if and only if 
$$\Theta_\tau(f \prod_{\theta\in \Sigma} G_\theta^{m_\theta} H_\theta^{n_\theta}) = 0$$
for all $\bm \in \ZZ^\Sigma$, $\bn \in \ZZ_{\ge 0}^\Sigma$, if and only if
$\Theta_\tau(f \prod_{\theta\in \Sigma} G_\theta^{m_\theta} H_\theta^{n_\theta}) = 0$
for some $\bm \in \ZZ^\Sigma$, $\bn \in \ZZ_{\ge 0}^\Sigma$.  (Alternatively note
that this follows from the fact the partial $\Theta$-operators commute with multiplication
by the $G_\theta$ and $H_\theta$, as can be seen directly from their definition.)

Suppose now that $\bk = \bk_{\min}(f)$, so that $f$ is not divisible by any partial Hasse invariants (see \S\ref{sec:stratification}).
Then if $f \in \ker(\Theta_\tau)$, and hence $f \in \ker(\Theta_{\tau_{\gp,i}})$ for all $i \in \ZZ/f_\gp\ZZ$,
then Theorem~\ref{thm:theta} implies that $p| k_{\theta_{\gp,i,e_\gp}}$ for all $i \in \ZZ/f_\gp\ZZ$.
 Therefore $\bk$ is of the form $\bk_0''$ for some $\bk_0$, where $\bk_0''$ is as in the definition of $V_\gp$
 in \S\ref{sec:V}, or equivalently $V_\gp^0$.  Furthermore it is immediate from the description of the effect of $V_\gp^0$
 on $q$-expansions in (\ref{eqn:Vonq}) that its image is contained in the kernel of $\Theta_\tau$.   We now use the
 method of \cite[Thm.~9.8.2]{DS} to prove the kernel is precisely the image of $V_\gp^0$.

\begin{theorem} \label{thm:kertheta}
Suppose that $f \in M_{\bk_0'',\bl}(U;\F)$ and $\tau \in \Sigma_{\gp,0}$.
If $\Theta_\tau(f) = 0$, then $f = V_\gp^0(g)$ for some 
$g \in  M_{\bk_0,\bl}(U;\F)$.
\end{theorem}
\begpf  Let $\ol{\iota}$ denote the embedding $\ol{Y}_U \hookrightarrow \ol{Y}_U^{\min}$,
and choose a set of cusps $\CS \subset Y^\infty_U$ consisting of precisely one on each connected component of $\ol{Y}_U$.
Note that since $\Phi_\gp^{\min}$ (defined in \S\ref{sec:V2}) is bijective on cusps as well as connected components, the set
$\CS' : =\Phi^\infty_{\gp}(\CS)$ also includes exactly one cusp on each connected component.

Recall from Proposition~\ref{prop:minimal} that the sheaves $\ol{\iota}_*\CA_{\bk_0,\bl,\FF}$ and $\ol{\iota}_*\CA_{\bk''_0,\bl,\FF}$ are
coherent, as is $\ol{\iota}_*\Phi_{\gp,*}\CA_{\bk_0'',\bl'',\FF} = \Phi_{\gp,*}^{\min}\ol{\iota}_*\CA_{\bk_0'',\bl,\FF}$
since $\Phi_{\gp}^{\min}$ is finite.  For each $\calC \in \CS$, let $\calC' = \Phi_{\gp}^\infty(\calC)$, so that
$\Phi_{\gp}^{{\min},*}$ defines a finite extension $\CO_{\ol{Y}_U^{\min},\calC'} \hookrightarrow
\CO_{\ol{Y}_U^{\min},\calC}$ of local rings.  We let $N_{\calC'} = (\ol{\iota}_*\CA_{\bk_0,\bl,\FF})_{\calC'}$
denote the stalk at $\calC'$ of $\ol{\iota}_*\CA_{\bk_0,\bl,\FF}$, and similarly let 
$N''_{\calC} = (\ol{\iota}_*\CA_{\bk_0,\bl,\FF})_{\calC} = (\ol{\iota}_*\Phi_{\gp,*}\CA_{\bk_0'',\bl,\FF} )_{\calC'}$.  
The stalk at $\calC'$ of $\ol{\iota}_*$ of the adjoint of
$\Phi_{\gp}^*\CA_{\bk_0,\bl,\FF} \stackrel{\sim}{\longrightarrow} \CA_{\bk_0'',\bl,\FF}$ then defines
an injective homomorphism $N_{\calC'} \to N''_{\calC}$ of finitely generated
$\CO_{\ol{Y}_U^{\min},\calC'}$-modules, extending $V_\gp^0$ to a map 
$$\bigoplus_{\calC'\in\CS'} N_{\calC'} \to \bigoplus_{\calC\in \CS} N''_{\calC}.$$
Similarly localizing at the generic points of $\ol{Y}_U$ (or equivalently $\ol{Y}_U^{\min}$) extends $V_\gp^0$
to a map $H^0(\ol{Y}_U, \CA_{\bk_0,\bl,\FF}\otimes_{\ol{Y}_U} \CF_U) \to H^0(\ol{Y}_U, \CA_{\bk''_0,\bl,\FF}\otimes_{\ol{Y}_U} \CF_U)$,
so we obtain a commutative diagram of injective maps
\begin{equation} \label{eqn:Vdiag}
\xymatrix{M_{\bk_0,\bl}(U;\FF)  \ar[r]\ar[d] & \bigoplus_{\calC'\in \CS'}  N_{\calC'}  \ar[r]\ar[d]  & H^0(\ol{Y}_U, \CA_{\bk_0,\bl,\FF} \otimes_{\ol{Y}_U} \CF_U) \ar[d] \\
M_{\bk''_0,\bl}(U;\FF)  \ar[r]& \bigoplus_{\calC\in \CS}  N''_{\calC}  \ar[r] & H^0(\ol{Y}_U, \CA_{\bk''_0,\bl,\FF} \otimes_{\ol{Y}_U} \CF_U).}\end{equation}
(Note that the horizontal maps, defined by localization, are injective since $\CS$ and $\CS'$ each contain a unique cusp
on each component of $\ol{Y}_U$.)   

Let $N_{\calC'}^\wedge$ denote the completion of $N_{\calC'}$ with respect to the maximal ideal of 
$\CO_{\ol{Y}_U^{\min},\calC'}$, and similarly let $N''^\wedge_{\calC}$ denote the completion of 
$N''_{\calC}$ with respect to the maximal ideal of $\CO_{\ol{Y}_U^{\min},\calC}$, or equivalently $\CO_{\ol{Y}_U^{\min},\calC'}$.
Note that the map $N_{\calC'}^\wedge \to N''^\wedge_{\calC}$ is the one described by (\ref{eqn:Vonq}), or more
precisely its variant for $V_\gp^0$.

Recall from (\ref{eqn:Thetaq}) that if $f \in \ker \Theta_\tau$, then for each $\calC\in \CS$, the $q$-expansion of $f$:
$$\sum_{m\in N^{-1}M_+\cup\{0\}}(b^{\bk_0''} c^{\bl} \otimes  r_m) q^m \in N''^\wedge_{\calC}$$
satisfies $r_m = 0$ for all $m \not\in \gp N^{-1}M = N^{-1}M'_+$, where $M' = \gp M$,
and is therefore in the image of $N_{\calC'}^\wedge$.  
Since the completion $\CO^\wedge_{\ol{Y}_U^{\min},\calC'}$ is faithfully flat over $\CO_{\ol{Y}_U^{\min},\calC'}$,
it follows that the image of $f$ in $\bigoplus_{\calC\in \CS} N''_{\calC}$ is of the form $V_\gp^0(g)$
for some $g \in \bigoplus_{\calC'\in \CS'}  N_{\calC'}$, and hence that its image in 
$H^0(\ol{Y}_U, \CA_{\bk''_0,\bl,\FF} \otimes_{\ol{Y}_U} \CF_U)$ is of the form $V_\gp^0(g)$
for some $g \in H^0(\ol{Y}_U, \CA_{\bk''_0,\bl,\FF} \otimes_{\ol{Y}_U} \CF_U)$.

It just remains to prove that $g \in M_{\bk_0,\bl}(U;\FF)$, or equivalently that $\ord_z(g) \ge 0$
for all prime divisors $z$ on $\ol{Y}_U$.  To that end, note that the operators 
$V_{\gp'}$ for all $\gp' \in S_p$, and $\epsilon_{\bk,\bl}$ for all $\bk,\bl\in \ZZ^\Sigma$ (see \S\ref{sec:V}))
similarly extend to maps on stalks at generic points satisfying (\ref{eqn:ppower}), so that
$V_\gp(g) = f \prod_{\theta\in\Sigma} G_\theta^{l_\theta}$ and
$$g^p  = \left(\epsilon_{p\bk,p\bl} \prod_{\gp'\in S_p}  V_\gq^{e_\gq} \right)   (g)  = 
\left(\delta_{p\bk,p\bl} \prod_{\gp' \neq \gp}  V_{\gp'}^{e_\gp} \right)   V_\gp^{e_\gp-1}(f).$$
Therefore $p\cdot \ord_z(g) \ge 0$, and hence $\ord_z(g) \ge 0$.
\epf

For the following corollary, recall that $\Xi^{\min}$ is defined by (\ref{eqn:Ximin}) and that the main
result of \cite{DK2} states that if $f$ is a non-zero form in $M_{\bk,\bl}(U;\F)$, then $\bk_{\min}(f) \in \Xi^{\min}$.
\begin{corollary} \label{cor:kertheta}
Suppose that $f \in M_{\bk,\bl}(U;\F)$ and $\tau \in \Sigma_{\gp,0}$.
Then $\Theta_\tau(f) = 0$ if and only if there exist $\bk_0 \in \Xi^{\min}$, $\bn \in \ZZ^\Sigma_{\ge 0}$ and
$g \in  M_{\bk_0,\bl}(U;\F)$ such that $\bk = \bk_0'' + \sum_\theta n_\theta \bh_\theta$ and
$$f = V_\gp^0(g) \prod_\theta H_\theta^{n_\theta}.$$
\end{corollary}
\begpf  We have already seen that if $f = V_\gp^0(g) \prod_\theta H_\theta^{n_\theta}$, then $\Theta_\tau(f) = 0$.

For the converse, note that we may assume $\bk = \bk_{\min}(f)$, so that $\bk \in \Xi^{\min}$ and
$\bk = \bk_0''$ for some $\bk_0 \in  \ZZ^\Sigma$.  Therefore
 the theorem implies that $f= V_\gp^0(g)$ for some $g \in  M_{\bk_0,\bl}(U;\F)$.
Finally it is immediate from the definitions of $\Xi^{\min}$ and $\bk_0''$
that $\bk_0 \in \Xi^{\min}$ if (and only if) $\bk_0'' \in \Xi^{\min}$.
\epf

\subsection{Forms of partial weight $0$}  \label{sec:ptwt0}
We now apply our results on partial $\Theta$-operators to prove a partial positivity result for minimal weights of
Hilbert modular forms.  Recall the main result of \cite{DK2} proves that minimal weights $\bk = \sum k_\theta \be_\theta$
of Hilbert modular forms necessarily lie in the cone $\Xi^{\min}$, and hence satisfy $k_\theta \ge 0$ for all $\theta$,
and that forms with $\bk = \bf{0}$ are easily described by Proposition~\ref{prop:wt0}.  We prove the following
restriction on possible minimal weights $\bk$ with $k_\theta = 0$ for some (but not all) $\theta \in \Sigma$.

\begin{theorem}  \label{thm:ptwt0}  Suppose that $\gp \in S_p$ is such that $f_\gp > 1$, or $e_\gp > 1$ and $p > 3$.
Suppose that $f \in M_{\bk,\bl}(U;\F)$ is non-zero and $\bk = \bk_{\min}(f)$.  If $k_\theta = 0$ for some $\theta \in \Sigma_\gp$,
then $\bk = \bf{0}$.
\end{theorem}
\begpf  Writing simply $f = f_\gp$ and $e=e_\gp$, note that the hypotheses mean that $ef > 1$ and $p^f > 3$.
Choose any $\tau = \tau_{\gp,i} \in \Sigma_{\gp,0}$ and let $\theta_0 = \theta_{\gp,i,e}$.  
We will first prove that $\Theta_\tau(f) = 0$.

Note that since $\bk \in \Xi^{\min}$ and $k_\theta = 0$ for some $\theta \in \Sigma_\gp$, we in fact have
$k_\theta = 0$ for all $\theta \in \Sigma_\gp$.  In particular $p|k_{\theta_0}$, so Theorem~\ref{thm:theta}
implies that $H_{\theta_0}|\Theta_\tau(f)$.  Therefore if $\Theta_\tau(f) \neq 0$, then 
$\bk_{\min}(\Theta_\tau(f)) \le_{\Ha} \bk + 2\be_{\theta_0}$, i.e.,
$$ \bk + 2\be_{\theta_0} - \sum_{\theta\in \Sigma} m_\theta \bh_\theta \in \Xi^{\min}$$
for some integers $m_\theta \ge 0$.  Letting $m_r  = m_{\sigma^r\theta_0}$ for $r=1,\ldots,ef$,
this implies that
\begin{equation}\label{eqn:chain} \begin{array}{rl}
      &   m_1 - m_2           \le m_2 - m_3              \le     \cdots \le m_{e-1} - m_e               \le m_e - pm_{e+1} \\
\le &p(m_{e+1} - m_{e+2}) \le p(m_{e+2} - m_{e+3}) \le \cdots \le p(m_{2e-1} - m_{2e}) \le p(m_{2e} - pm_{2e+1}) \\
       &   \vdots \\
\le &p^{f-1}(m_{(f-1)e+1} - m_{(f-1)e+2}) \le \cdots \le p^{f-1}(m_{ef-1} - m_{ef}) \le
    p^{f-1}(2+m_{ef} - pm_1) \\
\le &p^f(m_1-m_2),
\end{array}\end{equation}
(with the obvious collapsing here and in subsequent inequalities if $e$ or $f=1$).
In particular all the expressions in (\ref{eqn:chain}) are non-negative, so we have
$$m_1 \ge m_2 \ge \cdots m_e \ge pm_{e+1} \ge pm_{e+2} \ge \cdots \ge p^{f-1}m_{ef-1} \ge p^{f-1}m_{ef}$$
and $2 + m_{ef} - pm_1 \ge 0$, which implies that $(p^f-1)m_{ef} \le 2$.  Since $p^f > 3$, it follows that
$m_{ef} = 0$, so $pm_1 \le 2$, which implies that either $m_1 = 0$, or $m_1 =1$ and $p=2$.
If $m_1 = 0$, then $m_r = 0$ for all $r$, which contradicts the final inequality in (\ref{eqn:chain}).
On the other hand if $m_1 = 1$ and $p=2$, then all the expressions in (\ref{eqn:chain}) are zero,
which in turn implies that $m_1 = p^{f-1}m_{ef}$, which again yields a contradiction.

We have now shown that $\Theta_\tau(f) = 0$.  Note that $\bk'' = \bk$ since $k_\theta = 0$ for all $\theta \in \Sigma_\gp$,
so Theorem~\ref{thm:theta} implies that $f = V_\gp^0(f_1)$ for some $f_1 \in M_{\bk,\bl}(U;\FF)$.  We may therefore
iterate the above argument to conclude that $f_1 = V_\gp^0(f_2)$ for some $f_2 \in M_{\bk,\bl}(U;\FF)$, and by induction that
for all $n \ge 1$, we have $f = (V_\gp^0)^n(f_n)$ for some $f_n \in M_{\bk,\bl}(U;\FF)$.  It follows that for all $n \ge 1$,
the $q$-expansion of $f$ at every cusp of $\ol{Y}_U$ satisfies $r_m = 0$ for all $m \not\in \gp^n M$, so in fact
the $q$-expansion of $f$ at every cusp is constant.

To prove that $\bk = \bf{0}$,  recall that $\Xi^{\min}$ is contained in the cone spanned by the partial Hasse invariants, 
so $\bk = \sum_{\theta\in \Sigma} s_\theta \bh_\theta$ for some $s_\theta \in \QQ_{\ge 0}$.  Furthermore the denominators
are divisors of $M  = \lcm\{\, p^{f_\gq} - 1 \,|\, \gq \in S_p\,\}$, so that $M\bk = \sum m_\theta  \bh_\theta$ for some
$m_\theta \in \ZZ_{\ge 0}$.   Similarly $M\bl = \sum n_\theta  \bh_\theta$ for some $n_\theta \in \ZZ$. 
Since $f$ has constant $q$-expansions, so does $f^M$, and therefore 
$$f^M = h \prod_{\theta\in \Sigma} (H_\theta^{m_\theta} G_\theta^{n_\theta})$$
for some $h \in H^0(\ol{Y}_U, \CO_{\ol{Y}_U})$.   

For each $\theta\in\Sigma$, the assumption that $\bk = \bk_{\min}(f)$ means that $f$ is not divisible by $H_\theta$,
so $\ord_z(f)  = 0$ for some irreducible component $z$ of $Z_\theta$.  On the other hand we have 
$M\ord_z(f) = \ord_z(f^M) \ge m_\theta$, so $m_\theta = 0$.  As this holds for all $\theta \in \Sigma$, we 
conclude that $\bk = \bf{0}$.
\epf

\subsection{The kernel revisited}  \label{sec:exact}
Finally we present a cleaner, but less explicit, variant of Corollary~\ref{cor:kertheta} describing the kernels
of partial $\Theta$-operators.

We first record the effect of $V_\gp$ on the partial Hasse invariants $H_\theta$.  For each prime $\gp \in S_p$,
we let $\beta_\gp = p^{-1}\varpi_\gp^{e_\gp} \in \CO_{F,\gp}^\times$.  It is straightforward to check,
directly from the definition of $V_\gp$ or from the description (\ref{eqn:Vonq})  of its effect on $q$-expansions 
(and those of the $H_\theta$ in \S\ref{sec:Hasse2}), that
if $\theta\not\in \Sigma_\gp$ then $V_\gp(H_\theta) = H_\theta$ , but if $\theta = \theta_{\gp,i,j}$ then
$$V_\gp(H_\theta) = \left\{\begin{array}{ll}
H_{\sigma^{-1}\theta}^p, & \mbox{if $e_\gp = 1$;}\\
\overline{\theta}(\beta_\gp)^{-1} H_{\sigma^{-1}\theta}^p, & \mbox{if $e_\gp > 1$ and $j=1$;}\\
\overline{\theta}(\beta_\gp) H_{\sigma^{-1}\theta}, & \mbox{if $e_\gp > 1$ and $j=2$;}\\
H_{\sigma^{-1}\theta},& \mbox{otherwise.}\end{array}\right.$$
Therefore we define the modified partial Hasse invariant to be
$H'_\theta = \overline{\theta}(\beta_{\gp'})H_\theta$ if $\theta = \theta_{\gp',i,1}$
for some $\gp' \in S_p$ and $i \in \ZZ/f_{\gp'}\ZZ$, 
and $H'_\theta = H_\theta$ otherwise, so that
$$V_\gp(H'_\theta) =  \left\{ \begin{array}{ll}  H'^{n_\theta}_{\sigma^{-1}\theta} & \mbox{if $\theta \in \Sigma_{\gp,0}$;}\\
H'_\theta & \mbox{if $\theta \not\in \Sigma_{\gp,0}$.}\end{array}\right.$$
Similarly letting $G'_\theta =  \overline{\theta}(\beta_{\gp'})G_\theta$ if $\theta = \theta_{\gp',i,1}$
and $G'_\theta = G_\theta$ otherwise, we have $V_\gp(G'_\theta) = G'^{n_\theta}_{\sigma^{-1}\theta}$
if $\theta \in \Sigma_{\gp}$ and $V_\gp(G'_\theta) = G_\theta$ if $\theta\not\in \Sigma_\gp$.

Now for any sufficiently small $U$ of level prime to $p$, consider the $\FF$-algebra
$$M_{\tot}(U;\FF) = \bigoplus_{\bk,\bl \in (\ZZ^\Sigma)^2} M_{\bk,\bl}(U;\FF)$$
of Hilbert modular forms of all weights and level $U$ (where we let $M_{\bk,\bl}(U;\FF) =0$
if $\chi_{\bk+2\bl,\FF}$ is non-trivial on $U\cap \CO_F^\times$).    We may then consider
$V_\gp$ (resp.~$\Theta_\tau$) as an $\FF$-algebra homomorphism (resp.~$\FF$-linear derivation)
$M_{\tot}(U;\FF) \to M_{\tot}(U;\FF)$ for any $\gp \in S_p$ and $\tau \in \Sigma_{\gp,0}$.
Furthermore letting $\gI$ denote the ideal 
$\langle \,H_\theta' - 1, G_\theta'-1\,\rangle_{\theta \in \Sigma}$ in $M_{\tot}(U;\FF)$
and $R_U = M^{\tot}(U;\FF)/\gI$,  we see that $V_\gp(\gI) \subset \gI$ and $\Theta_\tau(\gI) \subset \gI$, 
so we obtain an $\FF$-algebra homomorphism $V_\gp$ and derivation $\Theta_\tau$
such that the composite
$$R_U \stackrel{V_\gp}{\longrightarrow} R_U \stackrel{\Theta_\tau}{\longrightarrow} R_U$$
is zero for any $\gp \in S_p$, $\tau \in \Sigma_{\gp,0}$.

Let $\Lambda$ denote the subgroup $\bigoplus_{\theta \in \Sigma}  \ZZ \bh_\theta$ of
$\ZZ^\Sigma = \bigoplus_{\theta \in \Sigma}  \ZZ \be_\theta$, so $\Lambda$ is the image
of the image of the endomorphism of $\ZZ^\Sigma$ defined by $\sum_\theta m_\theta \be_\theta
\mapsto \sum_\theta m_\theta \bh_\theta$.  Writing $\bh_\theta = \sum_{\theta'} n_{\theta,\theta'} \be_{\theta'}$,
it is straightforward to check that the matrix $(n_{\theta,\theta'})$ has determinant
$\prod_{\gp \in S_p}(p^{f_\gp} - 1)$, so this is the index of $\Lambda$ in  $\ZZ^\Sigma$.
On the other hand, let $\Psi$ denote the group of characters
$\psi: (\CO_F/\gq)^\times = \bigoplus_{\gp \in S_p} (\CO_F/\gp)^\times \longrightarrow \FF^\times,$
and consider the surjective homomorphism $\varrho: \ZZ^\Sigma \to \Psi$ defined by 
$$\varrho(\bk)  = \prod_{\theta\in \Sigma} \overline{\theta}^{k_\theta} = \bigoplus_{\gp \in S_p} (  \prod_{\theta \in \Sigma_\gp} \overline{\theta}^{k_\theta} ).$$
Note that $\varrho(\bh_\theta)$ is trivial for all $\theta \in \Sigma$, so $\Lambda \subset \ker(\varrho)$.
Since $\ZZ^\Sigma/\Lambda$ and $\Psi$ each have order $\prod_{\gp \in S_p}(p^{f_\gp} - 1)$, it follows that
$$\Lambda = \ker(\varrho) = \left\{\,\sum k_\theta \be_\theta\,\,\left|\,\, \sum_{i=1}^{f_\gp} \sum_{j=1}^{e_\gp} k_{\theta_{\gp,i,j}} p^i \equiv 1 \bmod (p^{f_\gp} - 1)\,\,
\forall \gp \in S_p\,\right.\right\}.$$

\begin{remark}  Recall that the $G_\theta$, and hence $G'_\theta$, are invertible in $M^\tot(U;\FF)$,
so if $\bl' - \bl = \sum_\theta m_\theta \bh_\theta \in \Lambda$, then multiplication by 
$\prod_{\theta} G'^{m_\theta}_\theta$ defines an isomorphism $M_{\bk,\bl}(U;\FF) \stackrel{\sim}{\to} M_{\bk,\bl'}(U;\FF)$.
We may therefore write $R_U$ as the quotient of $\displaystyle\bigoplus_{\bk\in\ZZ^\Sigma,\psi \in \Psi}  M_{\bk,\psi}(U;\FF)$
by the ideal $\langle\,H'_\theta - 1 \,\rangle_{\theta\in \Sigma}$, where $M_{\bk,\varrho(\bl)}(U;\FF)$ is canonically isomorphic
to $M_{\bk,\bl}(U;\FF)$ for each $\bl \in \ZZ^\Sigma$.  Furthermore the main result of \cite{DK2} immediately implies the natural map
$$\bigoplus_{\bk \in  \Xi^{\min}, \psi \in \Psi}  M_{\bk,\psi}(U;\FF) \longrightarrow R_U$$
is surjective, so we may also replace $\ZZ^\Sigma$ by the submonoid $\Xi^{\min} + \Lambda$
as the index set for $\bk$ in the definition of $R_U$.
\end{remark}

We will now describe the ideal $\gI$ in terms of $q$-expansions.  For each cusp $\calC \in Y_U^\infty$ we choose a representative
$(H,I,[\lambda],[\eta])$, and for each $\bk,\bl \in \ZZ^\Sigma$, we let $\overline{D}_{\bk,\bl}^\calC$ denote the one-dimensional
vector space $\FF\otimes_{\CO} D_{\bk,\bl}$ over $\FF$ (where $D_{\bk,\bl}$ is defined by (\ref{eqn:Dline})).  We then let
$\overline{D}_\tot^{\calC}$ denote $\bigoplus_{\bk,\bl \in \ZZ^\Sigma}  \overline{D}_{\bk,\bl}$ with its natural $\FF$-algebra structure.
For each $\theta \in \Sigma$, let $c_\theta \in \overline{D}^{\calC}_{\bh_\theta,\bf{0}}$ denote the (constant) $q$-expansion
of $H_\theta'$ at $\calC$.
(Recall that the $q$-expansions of $H_\theta$ were explicitly described in \S\ref{sec:Hasse2}, from which one gets a
description of $c_\theta$ by multiplying by $\overline{\theta}(b_\gp)$ if $\theta = \theta_{\gp,i,1}$.)  Similarly let
$d_\theta \in \overline{D}_{\bf{0},\bh_\theta}^\calC$ denote the (constant) $q$-expansion of $G'_\theta$ at $\calC$, and define 
$\gI^{\calC}$ to be the ideal $\langle\, c_\theta-1,d_\theta-1\,\rangle_{\theta \in \Sigma}$ of $\overline{D}_\tot^{\calC}$.
We may then view the quotient $\overline{D}_\tot^{\calC} / \gI^\calC$ as the space of $\Lambda^2$-coinvariants of
the free $\FF[(\ZZ^\Sigma)^2]$-module $\overline{D}_\tot^{\calC}$ and decompose 
$$\overline{D}_\tot^{\calC} / \gI^\calC  = \bigoplus_{\chi,\psi \in \Psi}  \overline{D}^\calC_{\chi,\psi},$$
so that the natural projection map $\overline{D}_{\bk,\bl}^{\calC} \to \overline{D}^\calC_{\varrho(\bk),\varrho(\bl)}$
is an isomorphism for all $\bk,\bl \in \ZZ^\Sigma$.  Now observe that the collection of $q$-expansion maps
$$M_{\bk,\bl}(U;\FF) \to \overline{D}_{\bk,\bl}^{\calC} \otimes_\FF \FF[[q^m]]_{m \in N^{-1}M^\calC_+ \cup \{0\}}$$
(where $M^\calC = \gd^{-1}I^{-1}J$ and $U(N) \subset U$) induces an $\FF$-algebra homomorphism
$$\overline{q}: M_\tot(U;\FF)  \longrightarrow \bigoplus_{\calC \in Y_U^\infty} \bigoplus_{\chi,\psi\in \Psi} 
\overline{D}_{\chi,\psi}^{\calC} \otimes_\FF \FF[[q^m]]_{m \in N^{-1}M^\calC_+ \cup \{0\}}.$$
\begin{lemma} \label{lem:qbar} The kernel of $\overline{q}$ is $\gI$.
\end{lemma}
\begpf  The inclusion $\gI \subset \ker(\overline{q})$ is clear from the definitions.   

Suppose then that $\overline{q}(f) = 0$ and write $f = \sum_{\bk,\bl \in W} f_{\bk,\bl}$ for some finite subset 
$W$ of  $\ZZ^\Sigma$.  For each $\chi \in \Psi$, choose $\bk_\chi \in \varrho^{-1}(\chi)$ sufficiently large
that $\bk \le_{\Ha} \bk_\chi$ for all $\bk \in \varrho^{-1}(\chi) \cap W$.
Thus for each $\bk \in W$, there is a unique $\bm_{\bk} = \sum_{\theta} m_{\bk,\theta} \be_\theta \in \ZZ^\Sigma_{\ge0}$ such that 
$\bk_{\varrho(\bk)} = \bk + \sum_{\theta} m_{\bk,\theta} \bh_\theta$.  Now note that
$$ g :=  \sum_{\bk,\bl \in W} H'^{\bm_{\bk}} G'^{\bm_{\bl}} f_{\bk,\bl}  \in \bigoplus_{\chi,\psi \in \Psi} M_{\bk_\chi,\bk_\psi}(U;\FF)$$
and that $f-g \in \gI$ (where $H'^{\bm_\bk} = \prod_\theta H_\theta'^{m_{\bk,\theta}}$ and 
$G'^{\bm_\bl} = \prod_\theta G_\theta'^{m_{\bl,\theta}}$).
Since $\gI \subset \ker(\overline{q})$, it follows that $\overline{q}(g) = 0$.
However $\overline{q}$ restricts to the $q$-expansion map
$$\bigoplus_{\chi,\psi \in \Psi} M_{\bk_\chi,\bk_\psi}(U;\FF) \longrightarrow 
\bigoplus_{\calC \in Y_U^\infty} \bigoplus_{\chi,\psi \in \Psi} 
\overline{D}_{\chi,\psi}^{\calC} \otimes_\FF \FF[[q^m]]_{m \in N^{-1}M^\calC_+ \cup \{0\}},$$
which is injective by Proposition~\ref{prop:qexp}, so $g = 0$, and hence $f = f-g \in \gI$.
\epf

We also extract the following observation from the proof of the lemma:
\begin{lemma} \label{lem:qbar2} If $W \subset (\ZZ^\Sigma)^2$ is such that $(\varrho,\varrho):W \to \Psi^2$ is injective, then
$$ \gI \cap \bigoplus_{(\bk,\bl) \in W} M_{\bk,\bl}(U;\FF)  = 0.$$
\end{lemma}

We are now ready to interpret the description of the kernel of the partial $\Theta$-operator in terms of
the algebra $R_U$.
\begin{theorem} \label{thm:exact} If $\gp \in S_p$, $\tau \in \Sigma_{\gp,0}$ and $U$ is any sufficiently small open
compact subgroup of $\GL_2(\AA_{F,\f})$ containing $\GL_2(\CO_{F,p})$, then the sequence
$$0 \longrightarrow R_U \stackrel{V_\gp}{\longrightarrow} R_U  \stackrel{\Theta_\tau}{\longrightarrow} R_U$$
is exact.
\end{theorem}
\begpf  The injectivity of $V_\gp: R_U \to R_U$ follows\footnote{Alternatively, one can
appeal to Lemma~\ref{lem:qbar2} instead of Lemma~\ref{lem:qbar} and argue similarly to the forthcoming
proof of the other exactness assertion.} from Lemma~\ref{lem:qbar} and the commutativity of the diagram
$$\begin{array}{ccc}  M_\tot(U;\FF)&  \stackrel{\overline{q}'}{\longrightarrow}& \displaystyle{\bigoplus_{\calC' \in Y_U^\infty} \bigoplus_{\chi,\psi\in \Psi}}
\overline{D}_{\chi,\psi}^{\calC'} \otimes_\FF \FF[[q^m]]_{m \in N^{-1}M^{\calC'}_+ \cup \{0\}}\\
{\scriptstyle{V_\gp}}\downarrow && \downarrow\\
M_\tot(U;\FF)&  \stackrel{\overline{q}}{\longrightarrow}&\displaystyle{ \bigoplus_{\calC \in Y_U^\infty} \bigoplus_{\chi,\psi\in \Psi}}
\overline{D}_{\chi,\psi}^{\calC} \otimes_\FF \FF[[q^m]]_{m \in N^{-1}M^{\calC}_+ \cup \{0\}},\end{array}$$
where the right downward arrow is the direct sum over $(\calC,\chi,\psi)$ of the tensor product of
the isomorphism $\overline{D}_{\chi,\psi}^{\calC'} \to \overline{D}_{\chi,\psi}^{\calC}$ induced by
(\ref{eqn:Phibres}) and the injective $\FF$-algebra homomorphism (\ref{eqn:Phat}).  (Note that
$(\varrho(\bk),\varrho(\bl)) = (\varrho(\bk''),\varrho(\bl''))$, and recall that the data $(H,I,[\lambda],[\eta])$
representing each cusp of $Y_U^\infty$ were implicitly chosen differently for the source and target of
$\Phi_\gp$ to simplify the resulting description of $V_\gp$ on $q$-expansions in \S\ref{sec:V2}.)

Since $\Theta_\tau\circ V_\gp = 0$, it just remains to prove that $\ker(\Theta_\tau) \subset \im(V_\gp)$.
Suppose then that $f \in M_\tot(U;\FF)$ is such that $\Theta_\tau(f) \in \gI$.  As in the proof of Lemma~\ref{lem:qbar},
write $f = \sum_{\bk,\bl \in W} f_{\bk,\bl}$ for some finite subset $W$ of  $\ZZ^\Sigma$, and for each $\chi \in \Psi$ choose
$\bk_\chi \in \varrho^{-1}(\chi)$ so that $\bk \le_{\Ha} \bk_\chi$ for all $\bk \in \varrho^{-1}(\chi) \cap W$, and consider
$$ g :=  \sum_{\bk,\bl \in W} H'^{\bm_{\bk}} G'^{\bm_{\bl}} f_{\bk,\bl}  \in \bigoplus_{\chi,\psi \in \Psi} M_{\bk_\chi,\bk_\psi}(U;\FF).$$
Since $f-g \in \gI$ and $\Theta_\tau(f) \in \gI$, we have $\Theta_\tau(g) \in \gI$.  Note however that
$$\Theta_\tau(g) \in  \bigoplus_{\chi,\psi \in \Psi} M_{\bk_\chi',\bl_\psi'}(U;\FF)$$
where $\bk_\chi' = \bk_\chi + n_{\theta_0}\be_{\sigma^{-1}\theta_0} + \be_{\theta_0}$ and $\bl_\psi' = \bk_\psi - \be_{\theta_0}$,
so Lemma~\ref{lem:qbar2} implies that $\Theta_\tau(g) = 0$.   Writing $g = \sum_{\chi,\psi} g_{\chi,\psi}$ with
$g_{\chi,\psi} \in M_{\bk_\chi,\bk_\psi}$, this means that $\Theta_\tau(g_{\chi,\psi}) = 0$ for all $\chi,\psi\in \Psi$.

Corollary~\ref{cor:kertheta} now implies that for each $\chi,\psi\in \Psi$, we have
$$g_{\chi,\psi} = G'^{\bm}H'^{\bn} V_\gp(h_{\chi,\psi})$$
for some $h_{\chi,\psi} \in M_{\bk_0,\bl_0}(U;\FF)$, where $\bk_0,\bl_0,\bm\in \ZZ^\Sigma$ and $\bn \in \ZZ^\Sigma_{\ge0}$
depend on $\chi$ and $\psi$.  It follows that $g_{\chi,\psi} -  V_\gp(h_{\chi,\psi}) \in \gI$, so setting $h = \sum_{\chi,\psi} h_{\chi,\psi}$,
we conclude that $f - V_\gp(h) = (f-g) + (g-V_\gp(h)) \in \gI$, as required.
\epf

Finally consider $R = \varinjlim_U R_U$, where the direct limit is over all sufficiently small open compact subgroups $U$
containing $\GL_2(\CO_{F,p})$.  (Note that this is the same as the quotient of
$$M_\tot(\FF) := \varinjlim_U M_\tot(U;\FF) \cong \bigoplus_{\bk,\bl \in \ZZ^\Sigma}  \left(\varinjlim_U M_{\bk,\bl}(U;\FF) \right)
  = \bigoplus_{\bk,\bl \in \ZZ^\Sigma}  M_{\bk,\bl}(\FF)$$
by the ideal $\langle \,H_\theta' - 1, G_\theta'-1\,\rangle_{\theta \in \Sigma}$.)
Since the maps $V_\gp$ and $\Theta_\tau$ are Hecke-equivariant in the obvious sense,
 we obtain an $\FF$-algebra endomorphism $V_{\gp}$ of $R$, and an
$\FF$-linear derivation $\Theta_\tau$ on $R$, each of which is $\GL_2(\A_{F,\f}^{(p)})$-equivariant, and
such that the sequence
$$0 \longrightarrow R \stackrel{V_\gp}{\longrightarrow} R \stackrel{\Theta_\tau}{\longrightarrow} R$$
is exact.

\bibliographystyle{abbrv} 

\bibliography{MRrefs_Theta}

 \end{document}